\documentclass[11pt]{amsart}

\usepackage{geometry}
\usepackage{amsmath}
\usepackage{amssymb}
\usepackage{amsthm}
\usepackage{amsfonts, dsfont}
\usepackage{paralist}
\usepackage{graphics}
\usepackage{epsfig} 
\usepackage{graphicx}  
\usepackage{epstopdf}
\usepackage{epstopdf}
\usepackage{verbatim}
\usepackage{mathrsfs}
\usepackage{mathtools}
\usepackage{pstricks}
\usepackage{relsize}
\usepackage{tikz}
\usetikzlibrary{matrix}
\usepackage{subcaption}
\usepackage{pgfplots}
\usepackage{fixltx2e}
\usepackage{enumitem}
\usepackage{upgreek}
\usepackage{url}
\usepackage[colorlinks=true]{hyperref}
\usepackage{cleveref}
\hypersetup{urlcolor=blue, linkcolor=blue, citecolor=red}
\usepackage{bm}
\usepackage{appendix}
\usepackage{algorithm2e}
\usepackage[protrusion=false,final]{microtype}

\allowdisplaybreaks

\newcommand{\D}[1]{\mbox{\rm #1}} 
\newcommand{\dd}{\D{d}}

\DeclareMathOperator*{\argmin}{argmin}

\usepackage[most]{tcolorbox}
\tcbset{
    mydefinitionstyle/.style={
        enhanced,
        breakable,
        colback=black!3,
        colframe=black!35,
        boxrule=0.35pt,
        arc=0.6mm,
        left=7pt,
        right=7pt,
        top=6pt,
        bottom=6pt,
        grow sidewards by=0.8cm
    }
}
\newenvironment{boxeddefinition}
    {\begin{tcolorbox}[mydefinitionstyle]\begin{definition}}
    {\end{definition}\end{tcolorbox}}

\usepackage[abbrev]{amsrefs}

\numberwithin{equation}{section}

\newtheorem{theorem}{Theorem}[section]
\newtheorem{lemma}[theorem]{Lemma}
\newtheorem{corollary}[theorem]{Corollary}
\newtheorem{proposition}[theorem]{Proposition}
\newtheorem{remark}[theorem]{Remark}
\newtheorem{definition}[theorem]{Definition}

\theoremstyle{definition}

\definecolor{ForestGreen}{RGB}{34,139,34}
\definecolor{ao(english)}{rgb}{0.0, 0.5, 0.0}

\begin{document}

\title[Distributional and Temporal Interactions in Games and Control]
{From Distributional to Temporal Interactions in McKean--Vlasov Games and Control}
\thanks{
}

\author{Hicham Kouhkouh}
\address{Hicham Kouhkouh \newline
Department of Mathematics and Scientific Computing, NAWI, University of Graz,\newline 
8010, Graz, Austria
}
\email{\texttt{hicham.kouhkouh@uni-graz.at}}

\begin{abstract}
The paper examines the relation between time and population in stochastic games and control. The idea is to substitute the interaction described via the population law with the interaction described via the occupation measure of a single self-interacting trajectory. Thus we study how McKean--Vlasov fixed point games and Mean Field Control problems compare with their resulting self-interacting counterparts. Assuming appropriate dissipativity and regularity conditions, we find that the ergodic costs, value functions and optimal feedback sets are identical, and the fixed point equilibria are equivalent. In addition, we obtain uniform bounds on both trajectories and costs for finite horizon problems. These bounds allow proving asymptotic equivalence of the normalized value functions, and enable the transfer of equilibrium strategies and optimal feedbacks from one representation to another, with approximation errors that vanish as the horizon length increases.
\end{abstract}



\maketitle

\vspace*{-2em}

\setcounter{tocdepth}{1}
\tableofcontents

\vspace*{-2.8em}

\section{Introduction}

The Birkhoff Ergodic Theorem \cite{birkhoff1931proof} establishes a direct relationship between ensemble average and time average. Under suitable assumptions, the statistical information contained in the invariant distribution can be obtained from the long-term observation of a single path. In the context of games and control, one may therefore ask if the information about an interacting population can be substituted with the information obtained over time along the path of a single agent.

Motivated by the work of K. Du, Y. Jiang, and J. Li \cite{du2023empirical} and F.-Y. Wang \cite{wang2018distribution}, we are interested in exploring how to replace the population law of McKean-Vlasov dynamics with the empirical (occupation) measure of a single self-interacting trajectory. To do so, we will examine the relationship between the two representations in both ergodic and finite-horizon environments, for both games and control problems. This type of replacement has also been studied in \cite{huang2024self} within the framework of Consensus-Based Optimization. 

We begin with an ergodic McKean--Vlasov fixed-point game. In this model, the dynamics depends upon the endogenous population law, and the running cost is based on a predetermined population measure. This game is compared to its self-interacting counterpart, which we call the \emph{One-Player Game} (OPG). We show that both formulations have the same ergodic cost for each admissible feedback. Therefore, their values, optimal-feedback sets, and equilibrium pairs coincide whenever they exist. Then we consider ergodic Mean-Field Control (MFC) where the endogenous law affects the running cost. This model is also compared to its self-interacting analogue called \emph{Path-Dependent Control}. With some additional assumptions regarding the regularity of the cost functions, the same equivalence can be shown. The proofs rely on dissipativity and moment estimates to establish that both models converge to a single invariant measure and thus produce a single ergodic cost. Our framework includes globally Lipschitz feedbacks as well as certain locally Lipschitz feedbacks. 

We also examine the behavior of such correspondences during finite periods of time. We obtain uniform estimates of trajectories for a class of globally Lipschitz feedback policies, which we apply to compare finite horizon cost functions and normalized value functions. We define the deterministic \emph{averaged experience measure} (the expected value of the random occupation measure) for fixed point games. This provides a way to express the consistency condition for the self-interacting game. We demonstrate that the existence of an exact equilibrium for either of the finite horizon games implies the existence of an approximate equilibrium for the other game. Moreover, the difference in performance between the two games decreases as the horizon increases. Similarly, the performance loss associated with transferring optimal and/or approximately optimal feedbacks from one formulation of MFC/Path-Dependent Control to the other decreases as the horizon increases. 

The paper is organized as follows. \textbf{Section \ref{sec:path}} develops the probabilistic estimates for law-dependent and self-interacting dynamics, while \textbf{Section \ref{sec: admissible controls}} introduces the admissible feedback classes. \textbf{Section \ref{sec:ergodic mfg}} and \textbf{Section \ref{sec:ergodic mkv}} establish the ergodic game and control equivalences, respectively. \textbf{Section \ref{sec:longtime}} presents the uniform trajectory estimates and finite-horizon approximation results.  
\textbf{Appendix \ref{app: occupation averaged measures}} contains supplementary measure-theoretic results concerning occupation measures and their expectations. 
\textbf{Appendix \ref{sec: proofs}} contains additional proofs.

\section{McKean--Vlasov Dynamics and Self-Interacting Processes}\label{sec:path}

In this section, we collect the approximation results for McKean--Vlasov dynamics and self-interacting processes that will be used throughout the paper. 

We begin with the heuristic underlying the approximation. For an ergodic Markov process $X$ with invariant probability measure $\mu$, the ergodic theorem suggests that its occupation measure\footnote{In \cite{du2023empirical}, this measure is referred to as the \textit{empirical measure}. We use the term \textit{occupation measure}, which is customary in control theory \cite{huang2025non}; see also \cite{huang2025non2,bardi2023eikonal,bardi2024long}. Indeed, for every Borel set $\mathcal O\subset\mathbb R^{d}$,
\(\;
    \mathcal E_t[X](\mathcal O)=\frac{1}{t}\int_0^t\mathbf 1_{\mathcal O}(X_r)\,\dd r
\)\;
is the proportion of time spent by the process in $\mathcal O$ over  $[0,t]$.}
\begin{equation}\label{occupation measure}
    \mathcal E_t[X]:=\frac{1}{t}\int_0^t\delta_{X_r}\,\dd r=\int_0^1\delta_{X_{ts}}\,\dd s,\qquad t>0,
\end{equation}
converges, under suitable ergodicity assumptions, to $\mu$ as $t\to+\infty$. On the other hand, we denote by $\mathcal L_t[X]:=\mathrm{Law}(X_t)$ the time-$t$ marginal law of $X$. If both $\mathcal L_t[X]$ and $\mathcal E_t[X]$ approach the same invariant measure as $t\to+\infty$, it is natural to expect that the two measure-valued quantities become asymptotically close.

This observation motivates the comparison between the McKean--Vlasov SDE
\begin{equation}\label{MkV SDE}
    \dd X_t=b\big(X_t,\mathcal L_t[X]\big)\,\dd t+\sigma\big(X_t,\mathcal L_t[X]\big)\,\dd W_t
\end{equation}
and the corresponding self-interacting SDE
\begin{equation}\label{path SDE}
    \dd Y_t = b(Y_t, \mathcal E_t[Y])\,\dd t + \sigma(Y_t, \mathcal E_t[Y])\,\dd W_t.
\end{equation}
The results provided below provide a quantitative estimate of how well the approximation holds under appropriate dissipativity and regularity conditions.

The probabilistic framework is the usual one: We are given a filtered probability space $(\Omega,\mathfrak F,\mathbb F,\mathbb P)$, with $\mathbb F=(\mathfrak F_t)_{t\geq 0}$, satisfying the standard assumptions, and is equipped with a $d$-dimensional $\mathbb F$-Brownian motion $W=(W_t)_{t\geq 0}$. The expected value associated with the probability measure $\mathbb P$ is denoted by $\mathbb E$. All stochastic processes mentioned  are defined in this setting and are adapted to the filtration $\mathbb F$.

For an initial distribution $m\in\mathscr P_2(\mathbb R^d)$, let $X_0$ be an $\mathfrak F_0$-measurable random variable with $\mathcal L_0[X]=\mathrm{Law}(X_0)=m$. Whenever \eqref{MkV SDE} is well posed, its solution induces the nonlinear semigroup $(P_t^*)_{t\geq0}$ on $\mathscr P_2(\mathbb R^d)$ defined by $P_t^*m:=\mathcal L_t[X]$. A probability measure $\mu\in\mathscr P_2(\mathbb R^d)$ is said to be invariant for \eqref{MkV SDE} if $P_t^*\mu=\mu$ for all $t\geq0$. When the initial condition is deterministic, we write $X_0=z$ for $z\in\mathbb R^d$, or equivalently $\mathcal L_0[X]=\delta_z$. In this case, the dependence of $\mathbb E$ on the initial point $z$ is implicit. The same convention applies to the self-interacting process $Y$. For deterministic initial conditions, we also set $\mathcal E_0[Y]:=\delta_z$.

Finally, we denote by $\|\cdot\|$ the Frobenius norm on matrices. For $r\geq1$, let
\begin{equation*}
    \mathscr P_r(\mathbb R^d):=\left\{\mu\in\mathscr P(\mathbb R^d):\mu(|\cdot|^r):=\int_{\mathbb R^d}|x|^r\,\mu(\dd x)<+\infty\right\}.
\end{equation*}
In particular, $\mathscr P_2(\mathbb R^d)$ is endowed with the $2$-Wasserstein distance $\mathcal W_2$.

As in \cite{du2023empirical}, the main assumption we need is the following.

\noindent\textbf{Assumptions (A).}
The coefficients $b$ and $\sigma$ are continuous on $\mathbb R^d\times\mathscr P_2(\mathbb R^d)$ and bounded on bounded subsets. Moreover, there exist constants $\mathfrak a>\mathfrak b\geq0$, $\mathfrak c>0$, $\rho>0$, and $K\geq0$ such that for all $x,y\in\mathbb R^d$ and $\mu,\nu\in\mathscr P_2(\mathbb R^d)$
\begin{equation}\label{eq: assumption 1}
    2\langle b(x,\mu),x\rangle+(1+\rho)\|\sigma(x,\mu)\|^2\leq-\mathfrak c|x|^2+K\big[1+\rho+\mu(|\cdot|^2)\big],
\end{equation}
\begin{equation}\label{eq: assumption 2}
    2\langle b(x,\mu)-b(y,\nu),x-y\rangle+\|\sigma(x,\mu)-\sigma(y,\nu)\|^2\leq-\mathfrak a|x-y|^2+\mathfrak b\mathcal W_2(\mu,\nu)^2,
\end{equation}

Condition \eqref{eq: assumption 1} is a coercivity condition ensuring stability of the dynamics. It is related to the classical Khasminskii recurrence condition. Condition \eqref{eq: assumption 2} is a dissipativity (confluence) condition. We emphasize that no non-degeneracy assumption is imposed on the diffusion coefficient $\sigma$.

Under Assumptions (A), the McKean--Vlasov SDE \eqref{MkV SDE} is well posed and admits a unique invariant probability measure; see \cite{du2023empirical} and the references therein. In particular, the dissipativity condition \eqref{eq: assumption 2} yields exponential contraction of the associated nonlinear semigroup. We first recall the convergence of occupation measures. 

\begin{theorem}[{\cite[Theorem 1.2]{du2023empirical}}]\label{thm: emp}
Let Assumptions (A) be satisfied, and let $\mu\in\mathscr P_2(\mathbb R^d)$ be the unique invariant probability measure of \eqref{MkV SDE}. Let $Y$ solve \eqref{path SDE} with $Y_0=z\in\mathbb R^d$. Then, for every
\begin{equation}\label{eq: condition epsilon}
    0<\varepsilon<\left(1-\frac{\mathfrak b}{\mathfrak a}\right)\wedge\frac{\rho}{(d+2)(\rho+2)},
\end{equation}
it holds
\begin{equation}\label{eq: conv path to mu}
    \mathbb E\left[\mathcal W_2\big(\mathcal E_t[Y],\mu\big)^2\right]=O(t^{-\varepsilon}),\qquad t\to+\infty.
\end{equation}
\end{theorem}

For the law-dependent dynamics, we recall the following exponential convergence estimate from \cite{wang2018distribution} which holds under an additional regularity condition. 

\begin{theorem}[{\cite[Theorem 3.1(2)]{wang2018distribution}}]\label{thm: conv law to mu}
Let Assumptions (A) be satisfied. Assume also that there exist constants $K_{\sigma},K_b\geq0$ such that
\begin{equation}\label{extra assumptions}
\begin{aligned}
    \|\sigma(x,\mu)-\sigma(y,\nu)\|^2&\leq K_{\sigma}\big(|x-y|^2+\mathcal W_2(\mu,\nu)^2\big),\\
    |b(0,\mu)|^2&\leq K_b\big[1+\mu(|\cdot|^2)\big],
\end{aligned}
\end{equation}
for all $x,y\in\mathbb R^d$ and $\mu,\nu\in\mathscr P_2(\mathbb R^d)$. Let $\mu\in\mathscr P_2(\mathbb R^d)$ be the unique invariant probability measure of \eqref{MkV SDE}. Then, for every $z\in\mathbb R^d$ and $t\geq0$,
\begin{equation}\label{eq: conv law to mu}
    \mathcal W_2\big(\mathcal L_t[X],\mu\big)^2\leq e^{-(\mathfrak a-\mathfrak b)t}\mathcal W_2(\delta_z,\mu)^2,
\end{equation}
where $X$ solves \eqref{MkV SDE} with $X_0=z$. In particular, $\mathcal W_2(\mathcal L_t[X],\mu)^2=O(t^{-\varepsilon})$ as $t\to+\infty$ for every $\varepsilon>0$.
\end{theorem}

Combining these convergence results with the trajectory comparisons established in \cite{du2023empirical}, we obtain the following estimate.

\begin{proposition}\label{prop: estimate trajectories}
Let Assumptions (A) and \eqref{extra assumptions} be satisfied, and let $X,Y$ solve \eqref{MkV SDE} and \eqref{path SDE}, respectively, with $X_0=Y_0=z\in\mathbb R^d$, driven by the same Brownian motion. Then, for every $\varepsilon$ satisfying \eqref{eq: condition epsilon}, there exists a constant $C_{\varepsilon,z}>0$ such that, for all $t\geq1$,
\begin{equation}\label{eq: estimate trajectories}
\begin{aligned}
    \mathbb E\big[|X_t-Y_t|^2\big]&\leq C_{\varepsilon,z}\,t^{-\varepsilon},\quad \text{and } \; \quad 
    \mathbb E\left[\mathcal W_2\big(\mathcal L_t[X],\mathcal E_t[Y]\big)^2\right]&\leq C_{\varepsilon,z}\,t^{-\varepsilon}.
\end{aligned}
\end{equation}
\end{proposition}

\begin{proof}
These estimates are obtained using the auxiliary frozen process
\[
    \dd Z_t=b(Z_t,\mu)\,\dd t+\sigma(Z_t,\mu)\,\dd W_t,\qquad Z_0=z,
\]
driven by the same Brownian motion as $X$ and $Y$.  The detailed estimates are provided in the supplementary material \ref{sec: proofs}.
\end{proof}

\begin{lemma}\label{lem: short time moments}
Let Assumptions (A) be satisfied, and let $X,Y$ be the solutions of \eqref{MkV SDE} and \eqref{path SDE}, respectively, with $X_0=Y_0=z\in\mathbb R^d$. Then, for every $T_0>0$, there exists a constant $C_{T_0,z}>0$ such that
\begin{equation}\label{eq: short time moments}
    \sup_{0\leq t\leq T_0}\mathbb E|X_t|^2+\sup_{0\leq t\leq T_0}\mathbb E|Y_t|^2\leq C_{T_0,z}.
\end{equation}
In particular, we have \;
\(
    \int_0^1\mathbb E\left[|X_t-Y_t|+\mathcal W_2\big(\mathcal L_t[X],\mathcal E_t[Y]\big)\right]\dd t<+\infty.
\)
\end{lemma}

\begin{proof}
The moment bound for $X$ follows from It\^o's formula, the Lyapunov inequality \eqref{eq: assumption 1}, and Gr\"onwall's inequality. For the self-interacting process $Y$, one has
\[
    \mathbb E\big[\mathcal E_t[Y](|\cdot|^2)\big]=\frac1t\int_0^t\mathbb E|Y_s|^2\,\dd s,\qquad t>0.
\]
We can then apply It\^o's formula and \eqref{eq: assumption 1}, with standard localization, and obtain
\[
    \mathbb E|Y_t|^2\leq|z|^2+C\int_0^t\left(1+\frac1s\int_0^s\mathbb E|Y_r|^2\,\dd r\right)\dd s.
\]
Let $M(t):=\sup_{0\leq s\leq t}\mathbb E|Y_s|^2$. Then we obtain
\(
    M(t)\leq|z|^2+C\int_0^t(1+M(s))\,\dd s,
\)\;
which together with Gr\"onwall inequality, yields \eqref{eq: short time moments}.

Finally, the triangle inequality for $\mathcal W_2$ applied through $\delta_0$ gives
\[
    \mathcal W_2\big(\mathcal L_t[X],\mathcal E_t[Y]\big)^2\leq2\mathbb E|X_t|^2+\frac2t\int_0^t|Y_s|^2\,\dd s.
\]
Taking expectations and using \eqref{eq: short time moments}, together with Cauchy--Schwarz for $|X_t-Y_t|$, proves the integrability on $[0,1]$. The detailed estimates are provided in the supplementary material \ref{sec: proofs}.
\end{proof}

In the remaining part of this paper, we intend to study how to extend the estimates developed previously to dynamics with a control. To begin with, we first provide a class of feedback controls which preserve the desired structural properties. Next, we employ the comparison results to examine McKean-Vlasov fixed point games and optimal control problems.

\section{Admissible Feedback Policies}\label{sec: admissible controls}

The estimates of Section \ref{sec:path} require that the coefficients are dissipative with suitable growth and regularity. As the conditions must be satisfied in the controlled case after composing with a feedback, we shall define an abstract class of admissible functions $\mathcal A_{\mathrm{adm}}$ by using the same characteristics required for the comparison of trajectories. Afterward, we shall introduce the specific class $\mathcal A_{\gamma}^{\vartheta}$ of globally Lipschitz feedback controls and establish structural conditions under which \(\; \mathcal A_{\gamma}^{\vartheta}\subseteq\mathcal A_{\mathrm{adm}}. \)\; The results presented in Section \ref{sec:path} will therefore be applicable to the laws and self-interacting systems produced by the same feedback. Lastly, we will demonstrate that global Lipschitz continuity for the feedback is not a requirement. In fact, as long as there exists a reasonable localization of the dependence of the data upon the control, even feedbacks with local Lipschitz continuity in the state variables may be included into $\mathcal A_{\mathrm{adm}}$.

Throughout this section, we will assume that $\mathds A\subset\mathbb R^m$ is a non-empty compact set of controls. The controlled coefficients are maps
\[
    b:\mathbb R^d\times\mathscr P_2(\mathbb R^d)\times\mathds A\longrightarrow\mathbb R^d,\qquad \sigma:\mathbb R^d\times\mathscr P_2(\mathbb R^d)\times\mathds A\longrightarrow\mathbb R^{d\times d},
\]
and the running cost is $f:\mathbb R^d\times\mathscr P_2(\mathbb R^d)\times\mathds A\longrightarrow\mathbb R$.

The proofs of the results in this section are based on direct estimates derived from the conditions of dissipativity, Lipschitz continuity and localization. Because the verification of these conditions is primarily technical, we have included the proof as part of our supplementary material \ref{sec: proofs}.

\subsection{Abstract admissibility and globally Lipschitz feedbacks}\label{subsec: abstract admissibility}

We begin with a regularity assumption on the running cost.

\noindent\textbf{Assumption (B).} 
There is $C_f>0$ constant so that for all $\mu\in\mathscr P_2(\mathbb R^d)$, $x,y\in\mathbb R^d$, and $a,\tilde a\in\mathds A$,
\begin{equation}\label{eq: assumption f}
    |f(x,\mu,a)-f(y,\mu,\tilde a)|\leq C_f\big(|x-y|+|a-\tilde a|\big).
\end{equation}

Thus, the running cost has a global Lipschitz condition on the state and control variables for all measures. No regularity condition is imposed on the measure variable at this point. The latter will only be required in Section \ref{sec:ergodic mkv}, where the measure appearing in the running cost is also generated by the controlled process.

We now define admissibility using both the feedback-composed coefficients and the running cost.

\begin{definition}[Admissible feedback policies]\label{def: admissible controls}
A Borel measurable feedback map $\alpha:\mathbb R^d\times\mathscr P_2(\mathbb R^d)\rightarrow\mathds A$, where $\mathscr P_2(\mathbb R^d)$ is endowed with the topology induced by $\mathcal W_2$, is said to be \emph{admissible} if the following conditions are satisfied.
\begin{enumerate}[label=\textnormal{(\roman*)}]
    \item The following feedback-composed coefficients satisfy Assumptions (A)
    \[
        b^\alpha(x,\mu):=b(x,\mu,\alpha(x,\mu)),\qquad \sigma^\alpha(x,\mu):=\sigma(x,\mu,\alpha(x,\mu)).
    \]

    \item The coefficients $(b^\alpha,\sigma^\alpha)$ satisfy the additional conditions \eqref{extra assumptions} required for the trajectory comparison results of Section \ref{sec:path}.

    \item $\exists\, C_\alpha>0$ a constant such that, for all $x,y\in\mathbb R^d$ and $\mu,\eta,\nu\in\mathscr P_2(\mathbb R^d)$,
    \begin{equation}\label{eq: admissible cost condition}
        \big|f(x,\mu,\alpha(x,\eta))-f(y,\mu,\alpha(y,\nu))\big|\leq C_\alpha\big(|x-y|+\mathcal W_2(\eta,\nu)\big).
    \end{equation}
\end{enumerate}
The collection of all such feedbacks is denoted by $\mathcal A_{\mathrm{adm}}$.
\end{definition}

Note that global Lipschitz continuity of $\alpha$ is not required for admissibility. Furthermore, the constant $C_\alpha$ in \eqref{eq: admissible cost condition}, and also all of the other constants in \textnormal{(i)}--\textnormal{(ii)}, can be different depending on which specific feedback $\alpha$ is considered. This policy-by-policy formulation is sufficient for the comparison arguments used in the sequel.

\begin{remark}
\label{rem: feedback measurability}
The admissible controls are feedback policies. For each $\alpha \in \mathcal A_{\mathrm{adm}}$, the associated  processes are: $\alpha^{\texttt L}_t = \alpha(X_t, \mathcal L_t[X])$ and $\alpha^{\texttt E}_t = \alpha(Y_t, \mathcal E_t[Y])$. As $\mathcal L_t[X]$ is a deterministic process and $\mathcal E_t[Y]$ depends only upon the history of the system prior to $t$, the two controls will be non-anticipative. Additionally, because the state and measure-valued processes are continuous and adapted and because $\alpha$ is a Borel measurable function, the controls are progressively measurable.
\end{remark}

We will now consider a more explicitly defined type of feedback control that satisfies the requirements for being admissible as in Definition \ref{def: admissible controls}. 

\begin{definition}[Globally Lipschitz feedbacks]\label{def: lip controls}
Let us assume that $\gamma,\vartheta\geq0$. Let $\mathcal A_\gamma^\vartheta$ denote the set of mappings $\alpha:\mathbb R^d\times\mathscr P_2(\mathbb R^d)\longrightarrow\mathds A$ that satisfy the following:
\begin{equation}\label{eq: lip feedback}
    \alpha(x,\mu)-\alpha(y,\nu)\leq\gamma|x-y|+\vartheta\,\mathcal W_2(\mu,\nu), \quad \forall\,x,y\in\mathbb R^d,\,\mu,\nu\in\mathscr P_2(\mathbb R^d).
\end{equation}
\end{definition}

Each $\alpha \in \mathcal A_\gamma^\vartheta$ is continuous and therefore Borel measurable on $\mathbb R^d \times \mathscr P_2(\mathbb R^d)$. As a result, the requirement that the control be measurable, as stated in Definition \ref{def: admissible controls}, will always be satisfied.

The constants $\gamma$ and $\vartheta$ are fixed \emph{a priori} and are common to all feedbacks in $\mathcal A_\gamma^\vartheta$.  The class $\mathcal A_\gamma^\vartheta$ is a convenient sufficient class rather than the definition of admissibility itself. Globally Lipschitz feedbacks of this type are commonly used in stochastic control and mean-field models; see, for example, \cite[\S 3]{bao2023ergodic}, \cite[\S 5.5, \S 6.5, \S 8.2]{bensoussan2013mean}, \cite[(H6)]{caines2006large}, \cite[Lemma 1 and \S 4]{carmona2013probabilistic}, \cite[p. 74 and Lemma 3.3, p. 137]{carmona2018probabilistic}, and \cite[\S 3.1]{pham2018bellman}. Related approximation questions for stationary Markov policies by Lipschitz policies have also been investigated in \cite{pradhan2024near}, using Borkar topology \cite{borkar1989topology}.

\subsection{Sufficient structural conditions for admissibility}\label{subsec: sufficient admissibility}

We now give conditions on the controlled coefficients under which the globally Lipschitz feedbacks introduced in Definition \ref{def: lip controls} are admissible in the sense of Definition \ref{def: admissible controls}, i.e. $\mathcal A_\gamma^\vartheta\subseteq\mathcal A_{\mathrm{adm}}$. The main point is to ensure that, after composition with any $\alpha\in\mathcal A_\gamma^\vartheta$, the resulting drift and diffusion retain the dissipativity, regularity, and growth properties required by the approximation results of Section \ref{sec:path}, with constants that can be chosen uniformly over the whole class $\mathcal A_\gamma^\vartheta$.

To make these conditions explicit, we consider controlled coefficients such that
\begin{equation}\label{eq: b sigma separable}
    b(x,\mu,a)=b_0(x,\mu)+b_1(x,a),\qquad \sigma(x,\mu,a)=\sigma_0(x,\mu)+\sigma_1(x,a),
\end{equation}
for $(x,\mu,a)\in\mathbb R^d\times\mathscr P_2(\mathbb R^d)\times\mathds A$. The decomposition \eqref{eq: b sigma separable} separates the mean-field component of the dynamics from its control-dependent component. We impose the following assumptions.

\noindent\textbf{Assumption (C).}
The coefficients appearing in \eqref{eq: b sigma separable} satisfy the following conditions.

\begin{enumerate}[label=\textnormal{(\roman*)}]
    \item The functions $b_0:\mathbb R^d\times\mathscr P_2(\mathbb R^d)\rightarrow\mathbb R^d$ and $\sigma_0:\mathbb R^d\times\mathscr P_2(\mathbb R^d)\rightarrow\mathbb R^{d\times d}$ are continuous and bounded on bounded subsets. Also, there exist  $\mathfrak a_0>\mathfrak b_0\geq0$, $\mathfrak c_0>0$, $\rho_0>0$, and $K_0\geq0$ such that,  $\forall\, x,y\in\mathbb R^d$ and $\forall\,\mu,\nu\in\mathscr P_2(\mathbb R^d)$,
    \begin{equation}\label{eq: D coercivity}
        2\langle b_0(x,\mu),x\rangle+2(1+\rho_0)\|\sigma_0(x,\mu)\|^2\leq-\mathfrak c_0|x|^2+K_0\big[1+\rho_0+\mu(|\cdot|^2)\big],
    \end{equation}
    \begin{equation}\label{eq: D dissipativity}
\begin{aligned}
     &   2\langle b_0(x,\mu)-b_0(y,\nu),x-y\rangle+2\|\sigma_0(x,\mu)-\sigma_0(y,\nu)\|^2 \\
	& \qquad \qquad \qquad \qquad\qquad \qquad  \leq-\mathfrak a_0|x-y|^2+\mathfrak b_0\mathcal W_2(\mu,\nu)^2.
\end{aligned}
    \end{equation}
    In addition, there exist constants $K_{\sigma,0},K_{b,0}\geq0$ such that
    \begin{equation}\label{extra assumptions b0 sigma0}
    \begin{aligned}
        \|\sigma_0(x,\mu)-\sigma_0(y,\nu)\|^2&\leq K_{\sigma,0}\big(|x-y|^2+\mathcal W_2(\mu,\nu)^2\big),\\
        |b_0(0,\mu)|^2&\leq K_{b,0}\big[1+\mu(|\cdot|^2)\big].
    \end{aligned}
    \end{equation}

    \item The coefficients $b_1:\mathbb R^d\times\mathds A\rightarrow\mathbb R^d$ and $\sigma_1:\mathbb R^d\times\mathds A\rightarrow\mathbb R^{d\times d}$ are globally Lipschitz: there exist $C_b,C_\sigma\geq0$ such that, for all $x,y\in\mathbb R^d$ and $a,\tilde a\in\mathds A$,
    \begin{equation}\label{eq: D control lip}
    \begin{aligned}
        |b_1(x,a)-b_1(y,\tilde a)|&\leq C_b\big(|x-y|+|a-\tilde a|\big),\\
        \|\sigma_1(x,a)-\sigma_1(y,\tilde a)\|&\leq C_\sigma\big(|x-y|+|a-\tilde a|\big).
    \end{aligned}
    \end{equation}
\end{enumerate}

Since $\mathds A$ is compact, define
\begin{equation}\label{eq: Mb Msigma}
    M_b:=\max_{a\in\mathds A}|b_1(0,a)|,\qquad M_\sigma:=\max_{a\in\mathds A}\|\sigma_1(0,a)\|.
\end{equation}
It follows from \eqref{eq: D control lip} that, for every $x\in\mathbb R^d$ and $a\in\mathds A$,
\begin{equation}\label{eq: b1 sigma1 growth}
    |b_1(x,a)|\leq C_b|x|+M_b,\qquad \|\sigma_1(x,a)\|\leq C_\sigma|x|+M_\sigma.
\end{equation}

Assumption (C) separates the strong stability requirements from the control-dependent perturbation. Conditions \eqref{eq: D coercivity}--\eqref{eq: D dissipativity} impose the dissipativity and confluence properties required in Section \ref{sec:path} on the uncontrolled component $(b_0,\sigma_0)$, while \eqref{eq: D control lip} imposes standard global Lipschitz regularity on $(b_1,\sigma_1)$. The additional conditions \eqref{extra assumptions b0 sigma0} are needed for the trajectory estimate of Proposition \ref{prop: estimate trajectories}. The next two lemmas show that these properties are preserved after composition with a feedback in $\mathcal A_\gamma^\vartheta$, provided the dissipativity of $(b_0,\sigma_0)$ dominates the Lipschitz perturbation generated by the control-dependent terms.

\begin{lemma}\label{lem: conditions ok}
Let $(b,\sigma)$ be given by \eqref{eq: b sigma separable}, and assume that Assumption (C) is satisfied. Fix $\gamma,\vartheta\geq0$ and assume that
\begin{equation}\label{eq: condition}
\left\{
\begin{aligned}
    \mathfrak a_0&>\mathfrak b_0+2C_b(1+\gamma+\vartheta)+4C_\sigma^2\big((1+\gamma)^2+\vartheta^2\big),\\
    \mathfrak c_0&>2C_b+4(1+\rho_0)C_\sigma^2+1.
\end{aligned}
\right.
\end{equation}
Then, for every $\alpha\in\mathcal A_\gamma^\vartheta$, the feedback-composed coefficients
\[
    b^\alpha(x,\mu):=b(x,\mu,\alpha(x,\mu)),\qquad \sigma^\alpha(x,\mu):=\sigma(x,\mu,\alpha(x,\mu))
\]
satisfy Assumptions (A). More precisely, Assumptions (A) hold with $\rho=\rho_0$ and with 
\begin{equation}\label{eq: resulting abc}
\begin{aligned}
    \mathfrak a&:=\mathfrak a_0-2C_b(1+\gamma)-C_b\vartheta-4C_\sigma^2(1+\gamma)^2,\\
    \mathfrak b&:=\mathfrak b_0+C_b\vartheta+4C_\sigma^2\vartheta^2, \; \text{ and }\;\; \mathfrak c:=\mathfrak c_0-2C_b-4(1+\rho_0)C_\sigma^2-1>0
\end{aligned}
\end{equation}
for which $\mathfrak a>\mathfrak b\geq0$. In particular, the corresponding McKean--Vlasov and self-interacting dynamics are well posed.
\end{lemma}

\begin{lemma}\label{lem: extra conditions feedback}
Let $(b,\sigma)$ be given by \eqref{eq: b sigma separable}, and assume that Assumption (C) is satisfied. Fix $\gamma,\vartheta\geq0$. Then, for every $\alpha\in\mathcal A_\gamma^\vartheta$, the feedback-composed coefficients satisfy the additional conditions \eqref{extra assumptions}. More precisely, there exist constants $\widetilde K_\sigma,\widetilde K_b\geq0$, independent of $\alpha\in\mathcal A_\gamma^\vartheta$, such that, for all $x,y\in\mathbb R^d$ and $\mu,\nu\in\mathscr P_2(\mathbb R^d)$,
\begin{equation}\label{eq: extra feedback estimates}
\begin{aligned}
    \|\sigma^\alpha(x,\mu)-\sigma^\alpha(y,\nu)\|^2&\leq\widetilde K_\sigma\big(|x-y|^2+\mathcal W_2(\mu,\nu)^2\big),\\
    |b^\alpha(0,\mu)|^2&\leq\widetilde K_b\big[1+\mu(|\cdot|^2)\big].
\end{aligned}
\end{equation}
\end{lemma}

The preceding two lemmas verify the coefficient conditions entering Definition \ref{def: admissible controls}. It remains to check the corresponding regularity of the feedback-composed running cost, which follows from Assumption (B) and the Lipschitz property of the feedback.

\begin{proposition}\label{prop: lip subclass admissible}
Let Assumptions (B) and (C) be satisfied. Fix $\gamma,\vartheta\geq0$. Assume that the compatibility condition \eqref{eq: condition} holds. Then $\mathcal A_\gamma^\vartheta\subseteq\mathcal A_{\mathrm{adm}}$. Moreover, the constants in Definition \ref{def: admissible controls} can be chosen uniformly over $\alpha\in\mathcal A_\gamma^\vartheta$.
\end{proposition}

\begin{corollary}\label{cor: feedback trajectory estimates}
Let Assumptions (B) and (C) be satisfied. Fix $\gamma,\vartheta\geq0$. Let the compatibility condition \eqref{eq: condition} hold. For $\alpha\in\mathcal A_\gamma^\vartheta$, let $X^\alpha$ and $Y^\alpha$ solve
\begin{equation*}
\begin{aligned}
    \dd X_t^\alpha&=b^\alpha\big(X_t^\alpha,\mathcal L_t[X^\alpha]\big)\,\dd t+\sigma^\alpha\big(X_t^\alpha,\mathcal L_t[X^\alpha]\big)\,\dd W_t,\\
    \dd Y_t^\alpha&=b^\alpha\big(Y_t^\alpha,\mathcal E_t[Y^\alpha]\big)\,\dd t+\sigma^\alpha\big(Y_t^\alpha,\mathcal E_t[Y^\alpha]\big)\,\dd W_t,
\end{aligned}
\end{equation*}
with $X_0^\alpha=Y_0^\alpha=z\in\mathbb R^d$, driven by the same Brownian motion. Then, for every
\[
    0<\varepsilon<\left(1-\frac{\mathfrak b}{\mathfrak a}\right)\wedge\frac{\rho_0}{(d+2)(\rho_0+2)},
\]
where $\mathfrak a$ and $\mathfrak b$ are defined in \eqref{eq: resulting abc}, and for every $\alpha\in\mathcal A_\gamma^\vartheta$, there exists a constant $C_{\varepsilon,z,\alpha}>0$ such that, for all $t\geq1$,
\begin{equation}\label{eq: feedback trajectory estimates}
\begin{aligned}
    \mathbb E\big[|X_t^\alpha-Y_t^\alpha|^2\big]&\leq C_{\varepsilon,z,\alpha}\,t^{-\varepsilon},\; \text{and } \; 
    \mathbb E\left[\mathcal W_2\big(\mathcal L_t[X^\alpha],\mathcal E_t[Y^\alpha]\big)^2\right]&\leq C_{\varepsilon,z,\alpha}\,t^{-\varepsilon}.
\end{aligned}
\end{equation}
\end{corollary}

\begin{remark}\label{rem: uniformity constants}
Proposition \ref{prop: lip subclass admissible} shows that the structural conditions for being admissible may be verified simultaneously for all $\alpha\in\mathcal A_\gamma^\vartheta$. In particular, the constants $\mathfrak a$, $\mathfrak b$, and $\rho_0$ can be fixed independently of $\alpha$, meaning that the admissible range of exponents $\varepsilon$ in Corollary \ref{cor: feedback trajectory estimates} is applicable to each element of the entire class.

We do not, however, state that the multiplicative constant $C_{\varepsilon,z,\alpha}$ in \eqref{eq: feedback trajectory estimates} is uniform for all $\alpha\in\mathcal A_\gamma^\vartheta$. Uniformity of the multiplicative constant is not necessary for the ergodic comparison results presented below as they apply to each policy individually and apply on the larger set $\mathcal A_{\mathrm{adm}}$. Uniform trajectory estimates over $\mathcal A_\gamma^\vartheta$ will be proved separately in Section \ref{sec:longtime} using additional assumptions for the quantitative finite horizon comparisons.
\end{remark}

\subsection{A localized class of admissible feedback policies}\label{subsec: local controls}

The class $\mathcal A_{\gamma}^{\vartheta}$ considered above provides a convenient and explicit subclass of $\mathcal A_{\mathrm{adm}}$. Its global Lipschitz requirement, however, is a sufficient condition rather than an intrinsic feature of the comparison arguments developed in this work. Indeed, what is ultimately required is that, after composition with the feedback, the controlled coefficients satisfy the dissipativity, regularity, and growth conditions entering Definition \ref{def: admissible controls}.

We now illustrate this point by exhibiting a class of feedbacks that are only locally Lipschitz in the state variable, but whose composition with the controlled coefficients nevertheless gives rise to globally regular coefficients. To obtain such a class under explicit conditions, we impose in this subsection a stronger structural assumption: the dependence of the dynamics and of the running cost on the control is localized in the state variable.

More precisely, let $\chi:\mathbb R^d\longrightarrow[0,1]$ be a globally Lipschitz function such that, for some $R>0$, $\chi(x)=0$ whenever $|x|\geq R$. 
We denote its Lipschitz constant by $L_\chi$.

In this subsection, we consider controlled coefficients of the form
\begin{equation}\label{eq: localized coefficients}
    b(x,\mu,a)=b_0(x,\mu)+\chi(x)\beta(a),\qquad \sigma(x,\mu,a)=\sigma_0(x,\mu)+\chi(x)\varsigma(a),
\end{equation}
where $\beta:\mathds A\longrightarrow\mathbb R^d$ and $\varsigma:\mathds A\longrightarrow\mathbb R^{d\times d}$ are Lipschitz continuous. We denote their Lipschitz constants by $L_\beta$ and $L_\varsigma$, respectively. Since $\mathds A$ is compact, we may also set
\begin{equation}\label{eq: localized bounds}
    M_\beta:=\max_{a\in\mathds A}|\beta(a)|<+\infty,\qquad M_\varsigma:=\max_{a\in\mathds A}\|\varsigma(a)\|<+\infty.
\end{equation}

We assume that $(b_0,\sigma_0)$ satisfy Assumption (C)(i), including the additional conditions \eqref{extra assumptions b0 sigma0}. Thus, outside the bounded region on which $\chi$ is nonzero, the dynamics coincide with the uncontrolled dissipative dynamics generated by $(b_0,\sigma_0)$.

We also assume  the dependence of the running cost on the control is localized:
\begin{equation}\label{eq: localized running cost}
    f(x,\mu,a)=f_0(x,\mu)+\chi(x)\ell(a),
\end{equation}
where $\ell:\mathds A\rightarrow\mathbb R$ is Lipschitz continuous withconstant $L_\ell$. Since $\mathds A$ is compact, we  set $M_\ell:=\max_{a\in\mathds A}|\ell(a)|<+\infty$. We further assume that there exists $L_{f,0}>0$ such that
\begin{equation}\label{eq: f0 lip}
    |f_0(x,\mu)-f_0(y,\mu)|\leq L_{f,0}|x-y|,
\end{equation}
for all $x,y\in\mathbb R^d$, uniformly with respect to $\mu\in\mathscr P_2(\mathbb R^d)$. These conditions imply Assumption (B). Indeed, for every $x,y\in\mathbb R^d$, $\mu\in\mathscr P_2(\mathbb R^d)$, and $a,\tilde a\in\mathds A$,
\[
    |f(x,\mu,a)-f(y,\mu,\tilde a)|\leq(L_{f,0}+L_\chi M_\ell)|x-y|+L_\ell|a-\tilde a|.
\]

The localization in \eqref{eq: localized coefficients}--\eqref{eq: localized running cost} makes it possible to replace global Lipschitz continuity of the feedback in the state variable by local Lipschitz continuity. Indeed, only the behavior of the feedback on a bounded neighborhood of $\operatorname{supp}(\chi)$ enters the controlled part of the coefficients and of the running cost. This leads to the following.

\begin{definition}[Localized locally Lipschitz feedbacks]\label{def: local feedbacks}
Fix $\vartheta\geq0$. We denote by $\mathcal A_{\mathrm{loc}}^\vartheta$ the class of Borel measurable feedback maps $\alpha:\mathbb R^d\times\mathscr P_2(\mathbb R^d)\longrightarrow\mathds A$ satisfying: 
\begin{enumerate}[label=\textnormal{(\roman*)}]
    \item The feedback $\alpha$ is locally Lipschitz in the state variable, uniformly with respect to the measure argument. More precisely, for every $r>0$, there exists a constant $\gamma_r\geq0$ such that, for every $\mu\in\mathscr P_2(\mathbb R^d)$,
    \begin{equation}\label{eq: local lip state}
        |\alpha(x,\mu)-\alpha(y,\mu)|\leq\gamma_r|x-y|,\qquad x,y\in B_r.
    \end{equation}

    \item The feedback $\alpha$ is globally Lipschitz with respect to the measure argument: 
    \begin{equation}\label{eq: local class measure lip}
        |\alpha(x,\mu)-\alpha(x,\nu)|\leq\vartheta\,\mathcal W_2(\mu,\nu), \qquad x\in \mathbb R^d, \; \mu,\nu\in\mathscr P_2(\mathbb R^d).
    \end{equation}
\end{enumerate}
\end{definition}

Combining \eqref{eq: local lip state} and \eqref{eq: local class measure lip}, we obtain, for every $r>0$,
\begin{equation}\label{eq: local joint lip}
    |\alpha(x,\mu)-\alpha(y,\nu)|\leq\gamma_r|x-y|+\vartheta\mathcal W_2(\mu,\nu), \qquad x,y\in B_r,\; \mu,\nu\in\mathscr P_2(\mathbb R^d).
\end{equation}
The important point is that no uniform bound on $\gamma_r$ as $r\to+\infty$ is imposed.

For $\alpha\in\mathcal A_{\mathrm{loc}}^\vartheta$, the control-dependent contributions to the feedback-composed coefficients are defined by
\begin{equation}\label{eq: localized composed terms}
    B^\alpha(x,\mu):=\chi(x)\beta\big(\alpha(x,\mu)\big),\qquad \Sigma^\alpha(x,\mu):=\chi(x)\varsigma\big(\alpha(x,\mu)\big).
\end{equation}
Then $b^\alpha(x,\mu)=b_0(x,\mu)+B^\alpha(x,\mu)$ and $\sigma^\alpha(x,\mu)=\sigma_0(x,\mu)+\Sigma^\alpha(x,\mu)$. The following lemma shows that localization compensates for the lack of global Lipschitz continuity of $\alpha$ in the state variable.

\begin{lemma}\label{lem: localized composed lip}
Let $\alpha\in\mathcal A_{\mathrm{loc}}^\vartheta$, and define $B^\alpha$ and $\Sigma^\alpha$ by \eqref{eq: localized composed terms}. Then there exist constants $\widehat C_{b,\alpha},\widehat C_{\sigma,\alpha}\geq0$ such that, for all $x,y\in\mathbb R^d$ and $\mu,\nu\in\mathscr P_2(\mathbb R^d)$,
\begin{equation}\label{eq: localized composed lip}
\begin{aligned}
    |B^\alpha(x,\mu)-B^\alpha(y,\nu)|&\leq\widehat C_{b,\alpha}\big(|x-y|+\mathcal W_2(\mu,\nu)\big),\\
    \|\Sigma^\alpha(x,\mu)-\Sigma^\alpha(y,\nu)\|&\leq\widehat C_{\sigma,\alpha}\big(|x-y|+\mathcal W_2(\mu,\nu)\big).
\end{aligned}
\end{equation}
\end{lemma}

The cutoff $\chi$ compensates for the lack of global Lipschitz continuity of $\alpha$. Indeed, when both state variables lie in a bounded neighborhood of $\operatorname{supp}(\chi)$, the estimates follow from the local Lipschitz continuity of $\alpha$ and the Lipschitz properties of $\chi$, $\beta$, and $\varsigma$. When at least one state variable lies outside this neighborhood, the vanishing of $\chi$ outside $B_R$, together with its global Lipschitz continuity and the boundedness of $\beta$ and $\varsigma$, yields the same estimates. Thus, only the local Lipschitz constant of $\alpha$ near $\operatorname{supp}(\chi)$ enters the global bounds. The  proof is in the supplementary material \ref{sec: proofs}.

\begin{proposition}\label{prop: local subclass admissible}
Assume that $(b_0,\sigma_0)$ satisfy Assumption (C)(i). Let $\alpha\in\mathcal A_{\mathrm{loc}}^\vartheta$, and let $\widehat C_{b,\alpha}$ and $\widehat C_{\sigma,\alpha}$ be the constants given in Lemma \ref{lem: localized composed lip}. 
Then $\alpha\in\mathcal A_{\mathrm{adm}}$ if the following conditions hold
\begin{equation}\label{eq: localized compatibility}
    \mathfrak a_0>\mathfrak b_0+4\widehat C_{b,\alpha}+8\widehat C_{\sigma,\alpha}^2,\qquad \mathfrak c_0>1.
\end{equation}
\end{proposition}

Note that Proposition \ref{prop: local subclass admissible} does not establish the inclusion $\mathcal A_{\mathrm{loc}}^\vartheta\subseteq\mathcal A_{\mathrm{adm}}$. Indeed, the constants $\widehat C_{b,\alpha}$ and $\widehat C_{\sigma,\alpha}$ depend on the particular feedback $\alpha$ through its local Lipschitz constant on a bounded neighborhood of $\operatorname{supp}(\chi)$. The compatibility condition \eqref{eq: localized compatibility} must therefore be verified for each feedback individually. The proposition thus identifies a subclass of $\mathcal A_{\mathrm{loc}}^\vartheta$ consisting of admissible policies, which may include feedbacks that are not globally Lipschitz in the state variable.

\begin{remark}\label{rem: local nonuniform constants}
An important difference between $\mathcal A_\gamma^\vartheta$ and $\mathcal A_{\mathrm{loc}}^\vartheta$ concerns the dependence of the constants entering the admissibility estimates. For the globally Lipschitz class $\mathcal A_\gamma^\vartheta$, the constants $\gamma$ and $\vartheta$ are fixed uniformly over the whole class. Consequently, under the assumptions of Proposition \ref{prop: lip subclass admissible}, the structural constants entering Definition \ref{def: admissible controls} can be chosen uniformly with respect to $\alpha\in\mathcal A_\gamma^\vartheta$. By contrast, for a policy $\alpha\in\mathcal A_{\mathrm{loc}}^\vartheta$ satisfying \eqref{eq: localized compatibility}, the constants $\widehat C_{b,\alpha}$, $\widehat C_{\sigma,\alpha}$, and $C_\alpha$ obtained above may depend on the particular feedback through its local Lipschitz constant on a bounded neighborhood of $\operatorname{supp}(\chi)$. No uniform bound over the entire localized class is required.
\end{remark}

In the following sections, the main ergodic comparison results are  stated on $\mathcal A_{\mathrm{adm}}$, while the two constructions $\mathcal A_\gamma^\vartheta$ and $\mathcal A_{\mathrm{loc}}^\vartheta$ provide explicit settings in which the admissibility requirements can be verified.

\section{Ergodic McKean--Vlasov Fixed-Point Games}\label{sec:ergodic mfg}

\subsection{From classical MFGs to a McKean--Vlasov fixed-point game}

Mean Field Games (MFG) were introduced by Lasry and Lions \cite{lasry2006jeux,lasry2006jeux2,lasry2007mean} and independently by Huang, Caines and Malham\'e \cite{caines2006large}. This approach allows to find Nash equilibria within a very large population where agents interact. By now, there is a lot of research on how MFGs relate to finite player games (e.g. \cite{lacker2016general,lacker2017limit,lacker2020convergence,fischer2017connection,cardaliaguet2019master}), and there is also a great interest in analyzing MFGs under ergodic conditions and over infinite horizons (e.g. \cite{arapostathis2017solutions,cardaliaguet2019long,cohen2023analysis,cohen2026asymptotic,carmona2024probabilistic}).

In a traditional MFG, the population measure $\mu$ is given at the time of the representative agent's optimization and cannot be changed by the agent's individual action. An equilibrium in the MFG occurs when the distribution of the optimal response matches the invariant distribution. However, in our setting, the dynamics will depend on the endogenously determined law $\mathcal L_t[X^\alpha]$, which will vary with the feedback, while the measure $\mu$ used in the running cost does not change during optimization. Therefore, we refer to our game as a \emph{McKean--Vlasov fixed-point game}. Also, our game differs from Mean-Field Control (MFC) because in MFC, the endogenous law appears in the objective function as well. Other related works include McKean--Vlasov control using dynamic programming \cite{bayraktar2018randomized}; infinite horizon games and control using FBSDEs \cite{bayraktar2023infinite}; and average cost MFC in discrete time \cite{bayraktar2024average}.

In contrast to the law $\mathcal L_t[X^\alpha]$ being replaced with the occupation measure $\mathcal E_t[Y^\alpha]$ of a self-interacting trajectory, which builds upon \cite{du2023empirical}, there are several works that study approximations to the invariant distribution of an uncontrolled McKean--Vlasov dynamics, such as \cite{chassagneux2024computing,chassagneux2025note}, and others that study potential monotone ergodic MFGs by updating beliefs and making repeated best response decisions, as in \cite{bai:hal-05610310}. However, unlike these works, we shall keep the feedback fixed and then apply Ces\`aro convergence. As a result, we obtain ergodic equivalences, and provide quantitative comparisons of the finite-horizon behaviors of the law-dependent and self-interacting formulations, both at the same finite-horizon, which is different from the infinite-horizon approximations studied in \cite{carmona2024probabilistic}.

Recall that the expression of the feedback-composed coefficients $b^\alpha$ and $\sigma^\alpha$ can be found in Definition \ref{def: admissible controls}.

\begin{boxeddefinition}[Ergodic McKean--Vlasov Fixed-Point Game]
\label{def: ergodic MKV FPG}
The ergodic McKean--Vlasov fixed-point game is defined through the following optimization and consistency procedure.

\textbf{Step 1.}
Given $\mu\in\mathscr P_2(\mathbb R^d)$, consider the control problem
\begin{equation}\label{ergodic mfg}
\tag{e-MKV-FPG}
\left[
\begin{aligned}
&\inf_{\alpha\in\mathcal A_{\mathrm{adm}}}\mathcal J_z^{\texttt L}(\alpha,\mu)\,:=\,\limsup_{T\to+\infty}\frac1T\mathbb E\left[\int_0^T f\big(X_t^\alpha,\mu,\alpha(X_t^\alpha,\mathcal L_t[X^\alpha])\big)\,\dd t\right],\\
&\text{subject to }\quad\dd X_t^\alpha=b^\alpha\big(X_t^\alpha,\mathcal L_t[X^\alpha]\big)\,\dd t+\sigma^\alpha\big(X_t^\alpha,\mathcal L_t[X^\alpha]\big)\,\dd W_t,\quad X_0^\alpha=z\in\mathbb R^d,
\end{aligned}
\right.
\end{equation}
Its value function and optimal-feedback set are
\[
v^{\texttt L}(z,\mu):=\inf_{\alpha\in\mathcal A_{\mathrm{adm}}}\mathcal J_z^{\texttt L}(\alpha,\mu),\qquad \mathcal O^{\texttt L}(\mu):=\argmin_{\alpha\in\mathcal A_{\mathrm{adm}}}\mathcal J_z^{\texttt L}(\alpha,\mu).
\]
Here $z$ is fixed; we shall show that the ergodic cost and optimal-feedback set are independent of $z$. The latter may be empty if the infimum is not attained.

\textbf{Step 2.}
Whenever $\mathcal O^{\texttt L}(\mu)\neq\varnothing$, choose $\alpha^*\in\mathcal O^{\texttt L}(\mu)$ and let $\mu^{\alpha^*}$ denote the unique invariant probability measure of the McKean--Vlasov dynamics generated by $\alpha^*$.

\textbf{Step 3.}
Replace $\mu$ with $\mu^{\alpha^*}$ and return to Step 1.

\medskip

The set-valued map $\Phi^{\texttt L}:\mathscr P_2(\mathbb R^d)\rightrightarrows\mathscr P_2(\mathbb R^d)$ is
\(\;
\Phi^{\texttt L}(\mu):=\{\mu^\alpha:\alpha\in\mathcal O^{\texttt L}(\mu)\},
\)\;
with $\Phi^{\texttt L}(\mu)=\varnothing$ whenever $\mathcal O^{\texttt L}(\mu)=\varnothing$.

An equilibrium is a measure $\mu^*\in\mathscr P_2(\mathbb R^d)$ satisfying $\mu^*\in\Phi^{\texttt L}(\mu^*)$, or equivalently $\mu^*=\mu^{\alpha^*}$ for some $\alpha^*\in\mathcal O^{\texttt L}(\mu^*)$.
\end{boxeddefinition}

The distinction from a classical MFG concerns not only the coefficients but also the feedback. Even if $b$ and $\sigma$ are independent of their measure argument, the implemented control $\alpha(X_t^\alpha,\mathcal L_t[X^\alpha])$ may still depend on the endogenous law. Restricting admissible policies to state-dependent feedbacks $\alpha=\alpha(x)$ removes this dependence and recovers the classical ergodic MFG best-response structure with mean-field interaction through the running cost; see \cite[Eq. (7.18) and the discussion on p. 637]{carmona2018probabilistic}. 
The present formulation nevertheless retains an MFG-type optimality--consistency structure, although convergence of the associated fixed-point iteration is not established.

\subsection{The self-interacting fixed-point game}

We introduce a self-interacting counterpart of \eqref{ergodic mfg}. Motivated by Theorem \ref{thm: emp}, we replace the endogenous law entering the controlled dynamics by the occupation measure
\[
    \mathcal E_t[Y^\alpha]=\frac1t\int_0^t\delta_{Y_s^\alpha}\,\dd s,\qquad t>0.
\]
The prescribed measure $\mu$ in the running cost remains fixed during optimization and is updated through the same type of consistency procedure.

\begin{boxeddefinition}[Ergodic One-Player Game]
\label{def: ergodic OPG}
The ergodic One-Player Game is defined as follows.

\textbf{Step 1.}
Given $\mu\in\mathscr P_2(\mathbb R^d)$, consider the control problem
\begin{equation}\label{one player game}
\tag{e-OPG}
\left[
\begin{aligned}
    &\inf_{\alpha\in\mathcal A_{\mathrm{adm}}}\mathcal J_z^{\texttt E}(\alpha,\mu) \,:=\, \limsup_{T\to+\infty}\frac1T\mathbb E\left[\int_0^T f\big(Y_t^\alpha,\mu,\alpha(Y_t^\alpha,\mathcal E_t[Y^\alpha])\big)\,\dd t\right],\\
    &\text{subject to }\quad \dd Y_t^\alpha=b^\alpha\big(Y_t^\alpha,\mathcal E_t[Y^\alpha]\big)\,\dd t+\sigma^\alpha\big(Y_t^\alpha,\mathcal E_t[Y^\alpha]\big)\,\dd W_t,\quad Y_0^\alpha=z\in\mathbb R^d,
\end{aligned}
\right.
\end{equation}
Its value function and optimal-feedback set are
\[
    v^{\texttt E}(z,\mu):=\inf_{\alpha\in\mathcal A_{\mathrm{adm}}}\mathcal J_z^{\texttt E}(\alpha,\mu),\qquad \mathcal O^{\texttt E}(\mu):=\argmin_{\alpha\in\mathcal A_{\mathrm{adm}}}\mathcal J_z^{\texttt E}(\alpha,\mu).
\]
The initial condition $z$ is fixed, and the optimal-feedback set may be empty.

\textbf{Step 2.}
Whenever $\mathcal O^{\texttt E}(\mu)\neq\varnothing$, choose $\alpha^*\in\mathcal O^{\texttt E}(\mu)$. Let $\mu^{\alpha^*}$ denote the unique invariant probability measure of the McKean--Vlasov dynamics generated by the same feedback. By Theorem \ref{thm: emp},
\begin{equation}\label{asymptotic measure}
    \lim_{t\to+\infty}\mathbb E\left[\mathcal W_2\big(\mathcal E_t[Y^{\alpha^*}],\mu^{\alpha^*}\big)^2\right]=0.
\end{equation}

\textbf{Step 3.}
Replace $\mu$ with $\mu^{\alpha^*}$ and return to Step 1.

\medskip

\noindent The corresponding set-valued map is
\(\;
    \Phi^{\texttt E}(\mu):=\{\mu^\alpha:\alpha\in\mathcal O^{\texttt E}(\mu)\},
\)\;
with $\Phi^{\texttt E}(\mu)=\varnothing$ whenever $\mathcal O^{\texttt E}(\mu)=\varnothing$.

An equilibrium is a measure $\mu^*\in\mathscr P_2(\mathbb R^d)$ satisfying $\mu^*\in\Phi^{\texttt E}(\mu^*)$, or equivalently $\mu^*=\mu^{\alpha^*}$ for some $\alpha^*\in\mathcal O^{\texttt E}(\mu^*)$.
\end{boxeddefinition}

Thus, the difference between \eqref{ergodic mfg} and \eqref{one player game} lies in replacing the endogenous distributional interaction by temporal self-interaction, while retaining the same prescribed measure in the running cost. The following lemma identifies a common invariant representation for their ergodic costs.

\begin{lemma}
\label{lem: ergodic cost independent z}
Let $\mu\in\mathscr P_2(\mathbb R^d)$ and $\alpha\in\mathcal A_{\mathrm{adm}}$ be fixed. For $z\in\mathbb R^d$, let $X^{\alpha,z}$ and $Y^{\alpha,z}$ denote, respectively, the law-dependent and self-interacting processes generated by $\alpha$ and starting from $z$. Then the corresponding ergodic costs are independent of the initial condition. More precisely, for every $z,z'\in\mathbb R^d$,
\begin{equation}\label{eq: ergodic cost independent z}
    \mathcal J_z^{\texttt L}(\alpha,\mu)=\mathcal J_{z'}^{\texttt L}(\alpha,\mu),\qquad \mathcal J_z^{\texttt E}(\alpha,\mu)=\mathcal J_{z'}^{\texttt E}(\alpha,\mu).
\end{equation}
In fact, if $\mu^\alpha$ denotes the unique invariant probability measure associated with the McKean--Vlasov dynamics generated by $\alpha$, then both ergodic costs are equal to
\begin{equation}\label{eq: common invariant cost}
    \int_{\mathbb R^d}f\big(x,\mu,\alpha(x,\mu^\alpha)\big)\,\mu^\alpha(\dd x).
\end{equation}
Moreover, the $\limsup$ in the definition of each ergodic cost is an actual limit.
\end{lemma}

\begin{proof}
Fix $\alpha\in\mathcal A_{\mathrm{adm}}$ and $\mu\in\mathscr P_2(\mathbb R^d)$, and let $\mu^\alpha$ denote the unique invariant probability measure of the feedback-composed McKean--Vlasov dynamics. Define
\[
g_{\alpha,\mu}(x):=f\big(x,\mu,\alpha(x,\mu^\alpha)\big),\qquad \ell_{\alpha,\mu}:=\int_{\mathbb R^d}g_{\alpha,\mu}(x)\,\mu^\alpha(\dd x).
\]
Using condition \textnormal{(iii)} in Definition \ref{def: admissible controls}, one has
\begin{equation}\label{eq: g alpha lip}
|g_{\alpha,\mu}(x)-g_{\alpha,\mu}(y)|\leq C_\alpha|x-y|,\qquad x,y\in\mathbb R^d.
\end{equation}
Moreover, since $\mu^\alpha\in\mathscr P_2(\mathbb R^d)$, the quantity $\ell_{\alpha,\mu}$ is finite.

\smallskip
\noindent\textit{\textbf{Step 1} (Law-dependent dynamics).}
Let us define $m_t^z:=\mathcal L_t[X^{\alpha,z}]$. By Theorem \ref{thm: conv law to mu}, we have
\begin{equation}\label{eq: law convergence alpha z}
\mathcal W_2\big(m_t^z,\mu^\alpha\big)\longrightarrow0,\qquad t\to+\infty.
\end{equation}
In particular, since $m_t^z$ is the law of $X_t^{\alpha,z}$, we also have
\[
\mathbb E\left[f\big(X_t^{\alpha,z},\mu,\alpha(X_t^{\alpha,z},m_t^z)\big)\right]=\int_{\mathbb R^d}f\big(x,\mu,\alpha(x,m_t^z)\big)\,m_t^z(\dd x).
\]
We can now add and subtract $\int_{\mathbb R^d}g_{\alpha,\mu}(x)\,m_t^z(\dd x)$, to obtain
\[
\begin{aligned}
&\left|\mathbb E\left[f\big(X_t^{\alpha,z},\mu,\alpha(X_t^{\alpha,z},m_t^z)\big)\right]-\ell_{\alpha,\mu}\right|\\
&\leq\left|\int_{\mathbb R^d}\!\left[f\big(x,\mu,\alpha(x,m_t^z)\big)\!-\!g_{\alpha,\mu}(x)\right]\!m_t^z(\dd x)\right|\!+\!\left|\int_{\mathbb R^d}\!g_{\alpha,\mu}(x)\,m_t^z(\dd x)\!-\!\!\!\int_{\mathbb R^d}g_{\alpha,\mu}(x)\mu^\alpha(\dd x)\right|
\end{aligned}
\]
Using condition \textnormal{(iii)} in Definition \ref{def: admissible controls}, the first term is bounded by $C_\alpha\mathcal W_2(m_t^z,\mu^\alpha)$. For the second term, we use \eqref{eq: g alpha lip} and the Kantorovich--Rubinstein duality which gives the bound $C_\alpha\mathcal W_1(m_t^z,\mu^\alpha)$. Consequently, and recalling $\mathcal W_1\leq\mathcal W_2$, one gets
\[
\left|\mathbb E\left[f\big(X_t^{\alpha,z},\mu,\alpha(X_t^{\alpha,z},m_t^z)\big)\right]-\ell_{\alpha,\mu}\right|\leq2C_\alpha\mathcal W_2\big(m_t^z,\mu^\alpha\big).
\]
Finally, the convergence in \eqref{eq: law convergence alpha z} yields
\begin{equation}\label{eq: instantaneous law cost limit}
\mathbb E\left[f\big(X_t^{\alpha,z},\mu,\alpha(X_t^{\alpha,z},m_t^z)\big)\right]\longrightarrow\ell_{\alpha,\mu},\qquad t\to+\infty.
\end{equation}

We next justify the interchange of expectation and time integration. Condition \textnormal{(iii)} in Definition \ref{def: admissible controls} gives
\[
|f(x,\mu,\alpha(x,\eta))|\leq|f(0,\mu,\alpha(0,\delta_0))|+C_\alpha\big(|x|+\mathcal W_2(\eta,\delta_0)\big).
\]
Since $\mu$ is fixed, the first term is finite. Moreover,
\[
\mathcal W_2\big(m_t^z,\delta_0\big)^2=\mathbb E|X_t^{\alpha,z}|^2,\qquad \mathcal W_2\big(\mathcal E_t[Y^{\alpha,z}],\delta_0\big)^2=\frac1t\int_0^t|Y_s^{\alpha,z}|^2\,\dd s.
\]
The finite-time moment bounds of Lemma \ref{lem: short time moments} ensure integrability of the running costs over bounded time intervals. Fubini's theorem therefore allows us to exchange expectation and time integration. Since the instantaneous expected cost converges to $\ell_{\alpha,\mu}$ by \eqref{eq: instantaneous law cost limit}, continuous-time Ces\`aro convergence implies that its time average converges to the same limit. Consequently, the ergodic cost is an actual limit and
\begin{equation}\label{eq: J L invariant representation}
\mathcal J_z^{\texttt L}(\alpha,\mu)=\ell_{\alpha,\mu}=\int_{\mathbb R^d}g_{\alpha,\mu}(x)\,\mu^\alpha(\dd x).
\end{equation}

\smallskip
\noindent\textit{\textbf{Step 2} (Self-interacting dynamics).}
Set $\eta_t^z:=\mathcal E_t[Y^{\alpha,z}]$. By Theorem \ref{thm: emp},
\begin{equation}\label{eq: empirical convergence alpha z}
\mathbb E\left[\mathcal W_2\big(\eta_t^z,\mu^\alpha\big)^2\right]\longrightarrow0,\qquad t\to+\infty.
\end{equation}
Moreover, Proposition \ref{prop: estimate trajectories} and \eqref{eq: law convergence alpha z} imply
\begin{equation}\label{eq: law Y convergence}
\mathcal W_2\big(\mathcal L_t[Y^{\alpha,z}],\mu^\alpha\big)\longrightarrow0,\qquad t\to+\infty.
\end{equation}
Indeed, the triangle inequality and the synchronous coupling yield
\[
\mathcal W_2\big(\mathcal L_t[Y^{\alpha,z}],\mu^\alpha\big)\leq\left(\mathbb E|Y_t^{\alpha,z}-X_t^{\alpha,z}|^2\right)^{1/2}+\mathcal W_2\big(m_t^z,\mu^\alpha\big),
\]
and both terms converge to zero.

We now argue as in the law-dependent case, this time adding and subtracting $\mathbb E[g_{\alpha,\mu}(Y_t^{\alpha,z})]$. By condition \textnormal{(iii)} in Definition \ref{def: admissible controls} and the Lipschitz estimate \eqref{eq: g alpha lip}, we obtain
\[
\begin{aligned}
&\left|\mathbb E\left[f\big(Y_t^{\alpha,z},\mu,\alpha(Y_t^{\alpha,z},\eta_t^z)\big)\right]-\ell_{\alpha,\mu}\right|\\
&\quad\leq\left|\mathbb E\left[f\big(Y_t^{\alpha,z},\mu,\alpha(Y_t^{\alpha,z},\eta_t^z)\big)-g_{\alpha,\mu}(Y_t^{\alpha,z})\right]\right|+\left|\mathbb E[g_{\alpha,\mu}(Y_t^{\alpha,z})]-\ell_{\alpha,\mu}\right|\\
&\quad\leq C_\alpha\mathbb E\left[\mathcal W_2\big(\eta_t^z,\mu^\alpha\big)\right]+C_\alpha\mathcal W_1\big(\mathcal L_t[Y^{\alpha,z}],\mu^\alpha\big)\\
&\quad\leq C_\alpha\left(\mathbb E\left[\mathcal W_2\big(\eta_t^z,\mu^\alpha\big)^2\right]\right)^{1/2}+C_\alpha\mathcal W_2\big(\mathcal L_t[Y^{\alpha,z}],\mu^\alpha\big).
\end{aligned}
\]
Here the first term is controlled by the Lipschitz dependence on the measure argument of the feedback, while the second follows from \eqref{eq: g alpha lip} and the Kantorovich--Rubinstein duality, as above.

Both terms converge to zero by \eqref{eq: empirical convergence alpha z} and \eqref{eq: law Y convergence}. Hence
\begin{equation}\label{eq: instantaneous E cost limit}
\mathbb E\left[f\big(Y_t^{\alpha,z},\mu,\alpha(Y_t^{\alpha,z},\eta_t^z)\big)\right]\longrightarrow\ell_{\alpha,\mu},\qquad t\to+\infty.
\end{equation}
The finite-time integrability established in \textit{Step 1} justifies Fubini's theorem. Applying continuous-time Ces\`aro convergence to \eqref{eq: instantaneous E cost limit}, we conclude that this ergodic cost is an actual limit and
\begin{equation}\label{eq: J E invariant representation}
\mathcal J_z^{\texttt E}(\alpha,\mu)=\ell_{\alpha,\mu}=\int_{\mathbb R^d}g_{\alpha,\mu}(x)\,\mu^\alpha(\dd x).
\end{equation}

\smallskip
\noindent\textit{\textbf{Step 3} (Conclusion).}
Combining \eqref{eq: J L invariant representation} and \eqref{eq: J E invariant representation}, we obtain
\[
\mathcal J_z^{\texttt L}(\alpha,\mu)=\mathcal J_z^{\texttt E}(\alpha,\mu)=\int_{\mathbb R^d}f\big(x,\mu,\alpha(x,\mu^\alpha)\big)\,\mu^\alpha(\dd x).
\]
Since the right-hand side is independent of $z$, this proves \eqref{eq: ergodic cost independent z} and \eqref{eq: common invariant cost}.
\end{proof}

In view of Lemma \ref{lem: ergodic cost independent z}, we henceforth suppress the initial condition from the notation for the ergodic costs and values, writing $\mathcal J^{\texttt L}(\alpha,\mu)$, $\mathcal J^{\texttt E}(\alpha,\mu)$, $v^{\texttt L}(\mu)$, and $v^{\texttt E}(\mu)$. The optimal-feedback sets are likewise independent of $z$.

\subsection{Equivalence of the two fixed-point games}

We now establish the main ergodic comparison result.

\begin{boxeddefinition}[Equilibrium pair]\label{definition}
For $\star\in\{\texttt L,\texttt E\}$, an equilibrium pair of the corresponding fixed-point game is a pair $(\mu^*,\alpha^*)\in\mathscr P_2(\mathbb R^d)\times\mathcal A_{\mathrm{adm}}$ satisfying
\[
    \alpha^*\in\mathcal O^\star(\mu^*),\qquad \mu^*=\mu^{\alpha^*}.
\]
For the One-Player Game, the consistency condition is equivalently expressed as convergence of $\mathcal E_t[Y^{\alpha^*}]$ to $\mu^*$ in the sense of \eqref{asymptotic measure}.
\end{boxeddefinition}

The value functions $v^{\texttt L}(\mu)$ and $v^{\texttt E}(\mu)$ are defined independently of whether an optimal feedback or an equilibrium exists.

\begin{theorem}\label{thm: mfg}
Let $\mu\in\mathscr P_2(\mathbb R^d)$ be fixed. Then, for every $\alpha\in\mathcal A_{\mathrm{adm}}$,
\begin{equation}\label{eq: equality poliyc and value}
    \mathcal J^{\texttt L}(\alpha,\mu)=\mathcal J^{\texttt E}(\alpha,\mu), \qquad \text{consequently} \qquad v^{\texttt L}(\mu)=v^{\texttt E}(\mu).
\end{equation}
Moreover, whenever the common value is finite, every $\varepsilon$-optimal feedback for one problem is $\varepsilon$-optimal for the other. The optimal-feedback sets also coincide:
\(
    \mathcal O^{\texttt L}(\mu)=\mathcal O^{\texttt E}(\mu),
\)\;
which may be empty if the common infimum is not attained.
\end{theorem}

\begin{proof}
By Lemma \ref{lem: ergodic cost independent z}, for every $\alpha\in\mathcal A_{\mathrm{adm}}$, we have 
\[
    \mathcal J^{\texttt L}(\alpha,\mu)=\mathcal J^{\texttt E}(\alpha,\mu)=\int_{\mathbb R^d}f\big(x,\mu,\alpha(x,\mu^\alpha)\big)\,\mu^\alpha(\dd x)
\]
which is the first part of \eqref{eq: equality poliyc and value}. Taking the infimum over $\mathcal A_{\mathrm{adm}}$ gives the second part of \eqref{eq: equality poliyc and value}.

If the common value is finite and $\alpha^\varepsilon$ is $\varepsilon$-optimal for either formulation, then the equality of costs and values shows that the same feedback is $\varepsilon$-optimal for the other. Finally, the two objective functionals coincide pointwise on $\mathcal A_{\mathrm{adm}}$, so their sets of minimizers coincide, whether or not the infimum is attained.
\end{proof}

\begin{corollary}\label{cor: lip class main result}
Let $(b,\sigma)$ be given by \eqref{eq: b sigma separable}. Assume that Assumptions (B) and (C) are satisfied, fix $\gamma,\vartheta\geq0$, and assume that the compatibility condition \eqref{eq: condition} holds. Then, for every $\alpha\in\mathcal A_\gamma^\vartheta$ and every $\mu\in\mathscr P_2(\mathbb R^d)$,
\[
    \mathcal J^{\texttt L}(\alpha,\mu)=\mathcal J^{\texttt E}(\alpha,\mu),
\qquad \text{consequently} \qquad 
    \inf_{\alpha\in\mathcal A_\gamma^\vartheta}\mathcal J^{\texttt L}(\alpha,\mu)=\inf_{\alpha\in\mathcal A_\gamma^\vartheta}\mathcal J^{\texttt E}(\alpha,\mu).
\]
Whenever the common value over $\mathcal A_\gamma^\vartheta$ is finite, $\varepsilon$-optimal feedbacks are transferred between the two problems within this class. Their optimal-feedback sets within $\mathcal A_\gamma^\vartheta$ also coincide, whether or not the common infimum is attained.
\end{corollary}

\begin{proof}
By Proposition \ref{prop: lip subclass admissible}, $\mathcal A_\gamma^\vartheta\subseteq\mathcal A_{\mathrm{adm}}$. The conclusions follow from Theorem \ref{thm: mfg} by restricting both optimization problems to $\mathcal A_\gamma^\vartheta$.
\end{proof}

\begin{corollary}\label{cor: mfg}
In the situation of Theorem \ref{thm: mfg}, the two fixed-point maps coincide:
\(\; \Phi^{\texttt L}(\mu)=\Phi^{\texttt E}(\mu),\; \mu\in\mathscr P_2(\mathbb R^d).\)\;
Consequently, the ergodic McKean--Vlasov fixed-point game admits an equilibrium if and only if the ergodic One-Player Game admits an equilibrium. More precisely, the two problems have the same equilibrium pairs $(\mu^*,\alpha^*)\in\mathscr P_2(\mathbb R^d)\times\mathcal A_{\mathrm{adm}}$, whenever such pairs exist. The corresponding equilibrium values also coincide.
\end{corollary}

\begin{proof}
By Theorem \ref{thm: mfg}, we have $\mathcal O^{\texttt L}(\mu)=\mathcal O^{\texttt E}(\mu)$ for every $\mu\in\mathscr P_2(\mathbb R^d)$. Since the two dynamics generated by a fixed feedback are associated with the same measure $\mu^\alpha$, we also have 
\(
    \Phi^{\texttt L}(\mu)=\{\mu^\alpha:\alpha\in\mathcal O^{\texttt L}(\mu)\}=\{\mu^\alpha:\alpha\in\mathcal O^{\texttt E}(\mu)\}=\Phi^{\texttt E}(\mu).
\) 
Thus the statements concerning equilibrium pairs and values follow.
\end{proof}

\section{Ergodic McKean--Vlasov Optimal Control}\label{sec:ergodic mkv}

We will now study the control version of the problems we discussed in Section \ref{sec:ergodic mfg}. In the fixed point representation that was developed in the previous section, the population measure which appears in the running cost function was kept constant at all times during the optimization problem and changed only when the consistency condition was applied. In contrast, in this section we allow the population distribution generated by the controlled process to appear in both the dynamics and the running cost function. Therefore, the population distribution can be updated dynamically as the control is changed. This leads to the ergodic control of a McKean-Vlasov process, usually called Mean-Field Control (MFC). We can then make a comparison between this ergodic control and a self-interacting path-dependent version of the same problem which is obtained by replacing the law of the process therein with its occupation measure.

Very recent studies related to MFC focus on discounted infinite-horizon McKean--Vlasov control using dynamic programming \cite{ruda2026infinite}, ergodic Mean-Field Control via Hamilton--Jacobi--Bellman equations in the Wasserstein space \cite{fuhrman2025ergodic}, and the convergence and turnpike properties of linear-quadratic Mean-Field Control problems \cite{bayraktar2027ergodicity}.

Throughout this section, all of our optimization problems will be within the admissible set $\mathcal A_{\mathrm{adm}}$ introduced in Definition \ref{def: admissible controls}.

\begin{boxeddefinition}[Ergodic Mean-Field Control]
\label{def: ergodic mfc}
The ergodic Mean-Field Control problem is
\begin{equation}\label{mean field control}
\tag{e-MFC}
\left[
\begin{aligned}
    &\inf_{\alpha\in\mathcal A_{\mathrm{adm}}}\mathcal J^{\texttt L}_{\mathrm C,z}(\alpha)\,:=\,\limsup_{T\to+\infty}\frac1T\mathbb E\left[\int_0^T f\big(X_t^\alpha,\mathcal L_t[X^\alpha],\alpha(X_t^\alpha,\mathcal L_t[X^\alpha])\big)\,\dd t\right],\\
    &\text{subject to }\quad \dd X_t^\alpha=b^\alpha\big(X_t^\alpha,\mathcal L_t[X^\alpha]\big)\,\dd t+\sigma^\alpha\big(X_t^\alpha,\mathcal L_t[X^\alpha]\big)\,\dd W_t,\; X_0^\alpha=z\in\mathbb R^d,
\end{aligned}
\right.
\end{equation}
As before, $b^\alpha(x,\eta):=b(x,\eta,\alpha(x,\eta))$ and $\sigma^\alpha(x,\eta):=\sigma(x,\eta,\alpha(x,\eta))$. We denote the value of this problem by $v^{\texttt L}_{\mathrm C}(z):=\inf_{\alpha\in\mathcal A_{\mathrm{adm}}}\mathcal J^{\texttt L}_{\mathrm C,z}(\alpha)$.
\end{boxeddefinition}

The subscript ``C'' refers to ``control''. Unlike the fixed-point game, MFC is a centralized optimization problem in which the controller accounts for the effect of the feedback on the distribution entering the objective; see, for instance, \cite{carmona2013control,bensoussan2013mean}.

\begin{boxeddefinition}[Ergodic Path-Dependent Control]
\label{def: ergodic pdc}
The ergodic Path-Dependent Control problem is
\begin{equation}\label{path dep control}
\tag{e-PDC}
\left[
\begin{aligned}
    &\inf_{\alpha\in\mathcal A_{\mathrm{adm}}}\mathcal J^{\texttt E}_{\mathrm C,z}(\alpha)\,:=\, \limsup_{T\to+\infty}\frac1T\mathbb E\left[\int_0^T f\big(Y_t^\alpha,\mathcal E_t[Y^\alpha],\alpha(Y_t^\alpha,\mathcal E_t[Y^\alpha])\big)\,\dd t\right],\\
    &\text{subject to }\quad \dd Y_t^\alpha=b^\alpha\big(Y_t^\alpha,\mathcal E_t[Y^\alpha]\big)\,\dd t+\sigma^\alpha\big(Y_t^\alpha,\mathcal E_t[Y^\alpha]\big)\,\dd W_t,\; Y_0^\alpha=z\in\mathbb R^d,
\end{aligned}
\right.
\end{equation}
We denote the value of this problem by $v^{\texttt E}_{\mathrm C}(z):=\inf_{\alpha\in\mathcal A_{\mathrm{adm}}}\mathcal J^{\texttt E}_{\mathrm C,z}(\alpha)$.
\end{boxeddefinition}

As in Section \ref{sec:ergodic mfg}, the same feedback map is used in both problems, but with different measure-valued arguments:
\(
    \alpha_t^{\texttt L}=\alpha\big(X_t^\alpha,\mathcal L_t[X^\alpha]\big),\text{ and } \alpha_t^{\texttt E}=\alpha\big(Y_t^\alpha,\mathcal E_t[Y^\alpha]\big).
\)
Unlike in the fixed-point-game setting, the measure entering the running cost is also endogenous. We therefore require an additional regularity condition.

\noindent\textbf{Assumption (B$'$).}
Assumption (B) is satisfied and, in addition, there exists a constant $\widetilde C_f>0$ such that, for every $x\in\mathbb R^d$, $a\in\mathds A$, and $\mu,\nu\in\mathscr P_2(\mathbb R^d)$,
\begin{equation}\label{eq: cost measure lipschitz}
    |f(x,\mu,a)-f(x,\nu,a)|\leq\widetilde C_f\,\mathcal W_2(\mu,\nu).
\end{equation}

This assumption allows us to compare running costs evaluated at the different endogenous measures of the two controlled dynamics.

\begin{lemma}\label{lem: mfc cost comparison}
Let Assumption \textnormal{(B$'$)} be satisfied, and let $\alpha\in\mathcal A_{\mathrm{adm}}$. Then there exists a constant $\widetilde C_\alpha>0$ such that, for every $x,y\in\mathbb R^d$ and $\eta,\nu\in\mathscr P_2(\mathbb R^d)$,
\begin{equation}\label{eq: mfc cost comparison}
    \big|f\big(x,\eta,\alpha(x,\eta)\big)-f\big(y,\nu,\alpha(y,\nu)\big)\big|\leq\widetilde C_\alpha\big(|x-y|+\mathcal W_2(\eta,\nu)\big).
\end{equation}
Consequently, if $X^\alpha$ and $Y^\alpha$ denote the law-dependent and self-interacting processes generated by the same feedback $\alpha$, then, for every $t\geq0$,
\begin{equation}\label{eq: mfc trajectory cost comparison}
\begin{aligned}
    &\big|f\big(X_t^\alpha,\mathcal L_t[X^\alpha],\alpha(X_t^\alpha,\mathcal L_t[X^\alpha])\big)-f\big(Y_t^\alpha,\mathcal E_t[Y^\alpha],\alpha(Y_t^\alpha,\mathcal E_t[Y^\alpha])\big)\big|\\
    &\qquad\leq\widetilde C_\alpha\big(|X_t^\alpha-Y_t^\alpha|+\mathcal W_2(\mathcal L_t[X^\alpha],\mathcal E_t[Y^\alpha])\big).
\end{aligned}
\end{equation}
\end{lemma}

\begin{proof}
By condition \textnormal{(iii)} in Definition \ref{def: admissible controls} and Assumption \textnormal{(B$'$)},
\[
\begin{aligned}
    &\big|f\big(x,\eta,\alpha(x,\eta)\big)-f\big(y,\nu,\alpha(y,\nu)\big)\big|\\
    &\quad\leq\big|f\big(x,\eta,\alpha(x,\eta)\big)-f\big(y,\eta,\alpha(y,\nu)\big)\big| +\big|f\big(y,\eta,\alpha(y,\nu)\big)-f\big(y,\nu,\alpha(y,\nu)\big)\big|\\
    &\quad\leq C_\alpha|x-y|+(C_\alpha+\widetilde C_f)\mathcal W_2(\eta,\nu).
\end{aligned}
\]
Thus \eqref{eq: mfc cost comparison} holds with $\widetilde C_\alpha:=C_\alpha+\widetilde C_f$. Applying this estimate to the two controlled processes gives \eqref{eq: mfc trajectory cost comparison}.
\end{proof}

\begin{lemma}[Invariant representation and independence of the initial condition]\label{lem: mfc ergodic cost independent z}
Let $\alpha\in\mathcal A_{\mathrm{adm}}$ and assume that Assumption \textnormal{(B$'$)} holds. For $z\in\mathbb R^d$, let $X^{\alpha,z}$ and $Y^{\alpha,z}$ denote, respectively, the law-dependent and self-interacting processes generated by $\alpha$ and starting from $z$. Then the corresponding ergodic control costs are independent of the initial condition. More precisely, for every $z,z'\in\mathbb R^d$,
\begin{equation}\label{eq: mfc cost independent z}
    \mathcal J^{\texttt L}_{\mathrm C,z}(\alpha)=\mathcal J^{\texttt L}_{\mathrm C,z'}(\alpha),\qquad \mathcal J^{\texttt E}_{\mathrm C,z}(\alpha)=\mathcal J^{\texttt E}_{\mathrm C,z'}(\alpha).
\end{equation}
In fact, if $\mu^\alpha$ denotes the unique invariant probability measure of the McKean--Vlasov dynamics generated by $\alpha$, then
\begin{equation}\label{eq: mfc invariant representation}
    \mathcal J^{\texttt L}_{\mathrm C,z}(\alpha)=\mathcal J^{\texttt E}_{\mathrm C,z}(\alpha)=\int_{\mathbb R^d}f\big(x,\mu^\alpha,\alpha(x,\mu^\alpha)\big)\,\mu^\alpha(\dd x).
\end{equation}
Moreover, the $\limsup$ defining each ergodic cost is an actual limit.
\end{lemma}

\begin{proof}
Fix $\alpha\in\mathcal A_{\mathrm{adm}}$, let $\mu^\alpha$ denote the unique invariant probability measure of the associated McKean--Vlasov dynamics, and define
\[
    G_\alpha(x,\eta):=f\big(x,\eta,\alpha(x,\eta)\big),\qquad \ell_\alpha:=\int_{\mathbb R^d}G_\alpha(x,\mu^\alpha)\,\mu^\alpha(\dd x).
\]
By Lemma \ref{lem: mfc cost comparison}, we have 
\begin{equation}\label{eq: G alpha lip}
    |G_\alpha(x,\eta)-G_\alpha(y,\nu)|\leq\widetilde C_\alpha\big(|x-y|+\mathcal W_2(\eta,\nu)\big).
\end{equation}
In particular, $\ell_\alpha$ is finite.

\smallskip
\noindent\textit{\textbf{Step 1} (A general estimate).} 
We first derive an estimate that applies to both dynamics. Let $U$ be a square-integrable $\mathbb R^d$-valued random variable, and let $\boldsymbol\eta$ be a $\mathscr P_2(\mathbb R^d)$-valued random variable satisfying $\mathbb E[\mathcal W_2(\boldsymbol\eta,\delta_0)^2]<+\infty$. 
Since $\mu^\alpha$ is deterministic, we have
\[
\mathbb E[G_\alpha(U,\mu^\alpha)]=\int_{\mathbb R^d}G_\alpha(x,\mu^\alpha)\,\mathcal L[U](\dd x).
\]
Adding and subtracting this quantity, we obtain
\[
\begin{aligned}
\left|\mathbb E\big[G_\alpha(U,\boldsymbol\eta)\big]-\ell_\alpha\right|
&\leq\mathbb E\left[\big|G_\alpha(U,\boldsymbol\eta)-G_\alpha(U,\mu^\alpha)\big|\right]\\
& \qquad \qquad +\left|\int_{\mathbb R^d}G_\alpha(x,\mu^\alpha)\,\mathcal L[U](\dd x)-\int_{\mathbb R^d}G_\alpha(x,\mu^\alpha)\,\mu^\alpha(\dd x)\right|\\
&\leq\widetilde C_\alpha\mathbb E\big[\mathcal W_2(\boldsymbol\eta,\mu^\alpha)\big]+\widetilde C_\alpha\mathcal W_1\big(\mathcal L[U],\mu^\alpha\big)\\
&\leq\widetilde C_\alpha\left(\mathbb E\big[\mathcal W_2(\boldsymbol\eta,\mu^\alpha)\big]+\mathcal W_2\big(\mathcal L[U],\mu^\alpha\big)\right).
\end{aligned}
\]
Indeed, applying \eqref{eq: G alpha lip} with $x=y=U$ gives $|G_\alpha(U,\boldsymbol\eta)-G_\alpha(U,\mu^\alpha)|\leq\widetilde C_\alpha\mathcal W_2(\boldsymbol\eta,\mu^\alpha)$ almost surely, which controls the first term. On the other hand, applying \eqref{eq: G alpha lip} with $\eta=\nu=\mu^\alpha$ shows that $x\mapsto G_\alpha(x,\mu^\alpha)$ is $\widetilde C_\alpha$-Lipschitz. The second term is therefore bounded by $\widetilde C_\alpha\mathcal W_1(\mathcal L[U],\mu^\alpha)$ by the Kantorovich--Rubinstein duality. The last inequality follows from $\mathcal W_1\leq\mathcal W_2$.

\smallskip
\noindent\textit{\textbf{Step 2} (Application).}
For the law-dependent dynamics, set $m_t^z:=\mathcal L_t[X^{\alpha,z}]$. By Theorem \ref{thm: conv law to mu},
\begin{equation}\label{eq: mfc law convergence}
    \mathcal W_2\big(m_t^z,\mu^\alpha\big)\longrightarrow0,\qquad t\to+\infty.
\end{equation}
Applying the general estimate of \textit{Step 1} with $U=X_t^{\alpha,z}$ and $\boldsymbol\eta=m_t^z$, which is deterministic, gives
\[
    \left|\mathbb E\left[G_\alpha\big(X_t^{\alpha,z},m_t^z\big)\right]-\ell_\alpha\right|\leq2\widetilde C_\alpha\mathcal W_2\big(m_t^z,\mu^\alpha\big)\longrightarrow0.
\]

For the self-interacting dynamics, set $\eta_t^z:=\mathcal E_t[Y^{\alpha,z}]$. By Theorem \ref{thm: emp},
\begin{equation}\label{eq: mfc empirical convergence}
    \mathbb E\left[\mathcal W_2\big(\eta_t^z,\mu^\alpha\big)^2\right]\longrightarrow0,\qquad t\to+\infty.
\end{equation}
Moreover, the synchronous coupling of $X^{\alpha,z}$ and $Y^{\alpha,z}$, together with the triangle inequality, yields
\[
\begin{aligned}
    \mathcal W_2\big(\mathcal L_t[Y^{\alpha,z}],\mu^\alpha\big)
    &\leq\mathcal W_2\big(\mathcal L_t[Y^{\alpha,z}],m_t^z\big)+\mathcal W_2\big(m_t^z,\mu^\alpha\big)\\
    &\leq\left(\mathbb E|Y_t^{\alpha,z}-X_t^{\alpha,z}|^2\right)^{\frac1 2}+\mathcal W_2\big(m_t^z,\mu^\alpha\big).
\end{aligned}
\]
The first term tends to zero by Proposition \ref{prop: estimate trajectories}, and the second by \eqref{eq: mfc law convergence}. Consequently, we have 
\begin{equation}\label{eq: mfc Y law convergence}
    \mathcal W_2\big(\mathcal L_t[Y^{\alpha,z}],\mu^\alpha\big)\longrightarrow0,\qquad t\to+\infty.
\end{equation}
Applying the general estimate of \textit{Step 1} with $U=Y_t^{\alpha,z}$ and $\boldsymbol\eta=\eta_t^z$, we obtain
\[
\begin{aligned}
&\left|\mathbb E\left[G_\alpha\big(Y_t^{\alpha,z},\eta_t^z\big)\right]-\ell_\alpha\right|\\
&\!\leq\!\widetilde C_\alpha\!\left(\mathbb E\!\!\left[\mathcal W_2\big(\eta_t^z,\mu^\alpha\big)\right]\!\!+\!\mathcal W_2\big(\mathcal L_t[Y^{\alpha,z}],\mu^\alpha\big)\!\right) 
\!\leq\!\widetilde C_\alpha\!\left(\!\mathbb E\!\!\left[\mathcal W_2\big(\eta_t^z,\mu^\alpha\big)^2\right]^{1/2}\!\!\!\!+\!\mathcal W_2\big(\mathcal L_t[Y^{\alpha,z}],\mu^\alpha\big)\!\!\right)
\end{aligned}
\]
where the last inequality follows from Cauchy--Schwarz. Both terms tend to zero by \eqref{eq: mfc empirical convergence} and \eqref{eq: mfc Y law convergence}, respectively. Consequently, we have
\[
\mathbb E\left[G_\alpha\big(Y_t^{\alpha,z},\eta_t^z\big)\right]\longrightarrow\ell_\alpha,\qquad t\to+\infty.
\]
Finally, \eqref{eq: G alpha lip} implies
\(\;
    |G_\alpha(x,\eta)|\leq|G_\alpha(0,\delta_0)|+\widetilde C_\alpha\big(|x|+\mathcal W_2(\eta,\delta_0)\big).
\) 
The finite-time moment bounds of Lemma \ref{lem: short time moments} ensure integrability of both running costs over every bounded time interval wich justifies Fubini's theorem. Then the continuous-time Ces\`aro convergence applied to the preceding instantaneous cost limits, as we did for Lemma \ref{lem: ergodic cost independent z}, gives
\begin{equation*}
    \mathcal J^{\texttt L}_{\mathrm C,z}(\alpha) = \mathcal J^{\texttt E}_{\mathrm C,z}(\alpha) =\ell_\alpha=\int_{\mathbb R^d}f\big(x,\mu^\alpha,\alpha(x,\mu^\alpha)\big)\,\mu^\alpha(\dd x).
\end{equation*}
Thus both $\limsup$ criteria are actual limits, and their common value is independent of $z$. Combining the two identities proves \eqref{eq: mfc cost independent z} and \eqref{eq: mfc invariant representation}.
\end{proof}

In view of Lemma \ref{lem: mfc ergodic cost independent z}, we henceforth suppress the dependence on the initial state and write $\mathcal J^{\texttt L}_{\mathrm C}(\alpha)$ and $\mathcal J^{\texttt E}_{\mathrm C}(\alpha)$, as well as $v^{\texttt L}_{\mathrm C}$ and $v^{\texttt E}_{\mathrm C}$.

The policy-by-policy equality of the ergodic control costs follows directly from the invariant representation \eqref{eq: mfc invariant representation}.

\begin{theorem}\label{thm: mfc}
Let Assumption \textnormal{(B$'$)} be satisfied. Then, for every $\alpha\in\mathcal A_{\mathrm{adm}}$, it holds
\begin{equation}\label{eq: mfc equality}
    \mathcal J^{\texttt L}_{\mathrm C}(\alpha)=\mathcal J^{\texttt E}_{\mathrm C}(\alpha), \qquad \text{consequently} \qquad 
    v^{\texttt L}_{\mathrm C}=v^{\texttt E}_{\mathrm C}.
\end{equation}
Moreover, whenever the common value is finite, every $\varepsilon$-optimal feedback for one problem is $\varepsilon$-optimal for the other. The optimal-feedback sets
\[
    \mathcal O^{\texttt L}_{\mathrm C}:=\argmin_{\alpha\in\mathcal A_{\mathrm{adm}}}\mathcal J^{\texttt L}_{\mathrm C}(\alpha),\qquad \mathcal O^{\texttt E}_{\mathrm C}:=\argmin_{\alpha\in\mathcal A_{\mathrm{adm}}}\mathcal J^{\texttt E}_{\mathrm C}(\alpha)
\]
also coincide, that is $\mathcal O^{\texttt L}_{\mathrm C}=\mathcal O^{\texttt E}_{\mathrm C}$. Both sets may be empty if the  infimum is not attained.
\end{theorem}

\begin{proof}
The two ergodic costs coincide over $\mathcal A_{\mathrm{adm}}$ thanks to Lemma \ref{lem: mfc ergodic cost independent z}, which proves the first part of \eqref{eq: mfc equality}. Then, since both problems are optimized over the same admissible class, one gets equality of their values, $\varepsilon$-optimal feedbacks whenever the common value is finite, and optimal-feedback sets.
\end{proof}

\section{Finite-Horizon Approximations for Games and Control}\label{sec:longtime}

In this section we examine the finite horizon versions of both control and game problems that we had previously studied in the infinite horizon (ergodic) context. As opposed to the ergodic formulation, the finite horizon optimal feedbacks may now depend on the time horizon. Therefore, the policy-by-policy estimates that we developed for the ergodic case cannot be used to compare either values or optimizers. This motivates the need for some estimates on the trajectories that should be uniform for a given family of feedbacks. With such estimates, we will be able to evaluate finite horizon costs, values and approximate equilibria for games and finite horizon costs, values and approximately optimal feedbacks for controls.

\subsection{Uniform estimates for the dynamics}

For the remainder of this section, we will use the following standard assumptions: we are optimizing over a single, globally Lipschitz class $\mathcal A_\gamma^\vartheta$, with $\gamma$ and $\vartheta$ being positive real numbers; the structural assumptions from Proposition \ref{prop: lip subclass admissible} are satisfied; and we can select the constants for Assumptions \textnormal{(A)} and \eqref{extra assumptions} independently of $\alpha \in \mathcal A_\gamma^\vartheta$.

We additionally assume
\begin{equation}\label{eq: longtime coercivity}
    \mathfrak c>K,
\end{equation}
where $\mathfrak c>0$ and $K\geq0$ are common constants satisfying
\begin{equation}\label{eq: uniform lyapunov longtime}
    2\langle b^\alpha(x,\eta),x\rangle+(1+\rho)\|\sigma^\alpha(x,\eta)\|^2\leq-\mathfrak c|x|^2+K\big[1+\rho+\eta(|\cdot|^2)\big],
\end{equation}
uniformly in $\alpha\in\mathcal A_\gamma^\vartheta$, $x\in\mathbb R^d$, and $\eta\in\mathscr P_2(\mathbb R^d)$.

Condition \eqref{eq: longtime coercivity} ensures uniform second-moment bounds for the invariant measures associated with the admissible feedbacks, which are needed to obtain uniform higher-moment and empirical-approximation estimates.

\begin{lemma}[Uniform moment bounds for the frozen dynamics]\label{lem: uniform frozen moments}
Assume the standing conditions of the present section, including \eqref{eq: longtime coercivity}, and fix $z\in\mathbb R^d$. For every $\alpha\in\mathcal A_\gamma^\vartheta$, let $\mu^\alpha$ denote the unique invariant probability measure of the McKean--Vlasov dynamics generated by $\alpha$, and let $Z^\alpha$ solve
\begin{equation}\label{eq: frozen alpha process}
    \dd Z_t^\alpha=b^\alpha(Z_t^\alpha,\mu^\alpha)\,\dd t+\sigma^\alpha(Z_t^\alpha,\mu^\alpha)\,\dd W_t,\qquad Z_0^\alpha=z.
\end{equation}
Then we have
\begin{equation}\label{uniform moment}
    \sup_{\alpha\in\mathcal A_\gamma^\vartheta}\mu^\alpha(|\cdot|^2)<+\infty, \qquad \text{and} \qquad \sup_{\alpha\in\mathcal A_\gamma^\vartheta}\sup_{t\geq0}\mathbb E\big[|Z_t^\alpha|^{2+\rho}\big]<+\infty.
\end{equation}
\end{lemma}

\begin{proof}
Fix $\alpha\in\mathcal A_\gamma^\vartheta$. Applying It\^o's formula to the corresponding McKean--Vlasov process $X^\alpha$ and using \eqref{eq: uniform lyapunov longtime}, we obtain
\begin{equation*}
    \frac{\dd}{\dd t}\mathbb E|X_t^\alpha|^2\leq-(\mathfrak c-K)\mathbb E|X_t^\alpha|^2+K(1+\rho).
\end{equation*}
Since $\mathfrak c>K$, Gr\"onwall's inequality yields
\begin{equation}\label{eq: uniform second moment X}
    \mathbb E|X_t^\alpha|^2\leq e^{-(\mathfrak c-K)t}|z|^2+\frac{K(1+\rho)}{\mathfrak c-K}\big(1-e^{-(\mathfrak c-K)t}\big).
\end{equation}
By Theorem \ref{thm: conv law to mu}, $\mathcal L_t[X^\alpha]\to\mu^\alpha$ in $\mathcal W_2$. Passing to the limit in \eqref{eq: uniform second moment X}, we obtain
\begin{equation}\label{eq: proof uniform invariant moment}
    \sup_{\alpha\in\mathcal A_\gamma^\vartheta}\mu^\alpha(|\cdot|^2)\leq\frac{K(1+\rho)}{\mathfrak c-K}<+\infty.
\end{equation}

For the higher-moment estimate, set $V(x):=(1+|x|^2)^{(2+\rho)/2}$ and let $\mathscr L^\alpha$ denote the generator of the frozen process $Z^\alpha$. Using \eqref{eq: uniform lyapunov longtime} and \eqref{eq: proof uniform invariant moment}, a standard generator calculation followed by Young's inequality gives constants $c,C>0$, independent of $\alpha$, such that $\;\mathscr L^\alpha V(x)\leq-cV(x)+C$.\; 
Applying the exponentially weighted It\^o formula, with localization, and then Fatou's lemma yields
\[
    \mathbb E[V(Z_t^\alpha)]\leq e^{-ct}V(z)+\frac Cc(1-e^{-ct}),\qquad t\geq0.
\]
Since $|x|^{2+\rho}\leq V(x)$, this proves the second part of  \eqref{uniform moment}. The detailed calculations are given in the supplementary material \ref{sec: proofs}.
\end{proof}

We also establish a uniform second-moment bound on $[0,1]$ for the law-dependent and self-interacting dynamics. This estimate will be used below to control the contribution of the initial time interval in occupation-measure and trajectory-comparison estimates. The bound for $X^\alpha$ follows from \eqref{eq: uniform second moment X}. For $Y^\alpha$, It\^o's formula and the Lyapunov inequality \eqref{eq: uniform lyapunov longtime}, after dropping the nonpositive coercive term, give
\[
    \mathbb E|Y_t^\alpha|^2\leq|z|^2+K(1+\rho)t+K\int_0^t\frac1s\int_0^s\mathbb E|Y_r^\alpha|^2\,\dd r\,\dd s.
\]
Writing $M_\alpha(t):=\sup_{0\leq s\leq t}\mathbb E|Y_s^\alpha|^2$, we obtain
\[
    M_\alpha(t)\leq|z|^2+K(1+\rho)t+K\int_0^tM_\alpha(s)\,\dd s.
\]
Gr\"onwall's inequality therefore gives
\begin{equation}\label{eq: uniform short time second moments finite horizon}
    \sup_{\alpha\in\mathcal A_\gamma^\vartheta}\sup_{0\leq t\leq1}\left\{\mathbb E|X_t^\alpha|^2+\mathbb E|Y_t^\alpha|^2\right\}\leq C_{z,\gamma,\vartheta}.
\end{equation}

The proofs of the next two lemmas follow the approach of \cite[Proposition 3.1 and Lemma 4.1]{du2023empirical}. The uniform moment bounds, together with the common contraction and dissipativity constants, allow us to carry out these arguments with constants independent of $\alpha\in\mathcal A_\gamma^\vartheta$. We specialize the sampling measure in \cite{du2023empirical} to Lebesgue measure and account explicitly for the mollifier scaling in the covariance estimate. 

\begin{lemma}[Uniform empirical approximation for frozen dynamics]\label{lem: uniform frozen empirical}
Assume the standing conditions of the present section. Fix $z\in\mathbb R^d$. Let $\kappa:=\frac{\rho}{(d+2)(\rho+2)}$. Then, for every $\varepsilon\in(0,1)$, there exists $C_{z,\gamma,\vartheta,\varepsilon}>0$  independent of $\alpha\in\mathcal A_\gamma^\vartheta$, such that
\begin{equation}\label{eq: uniform frozen empirical}
    \mathbb E\left[\mathcal W_2\big(\mathcal E_t[Z^\alpha],\mu^\alpha\big)^2\right]\leq C_{z,\gamma,\vartheta,\varepsilon}\,t^{-\kappa\varepsilon},\qquad t\geq1.
\end{equation}
\end{lemma}

\begin{proof}
We adapt the argument of \cite[Proposition 3.1]{du2023empirical} to the frozen diffusions $Z^\alpha$ and verify that its constants are uniform over $\alpha\in\mathcal A_\gamma^\vartheta$.

Assumption \textnormal{(A)}, with the measure fixed at $\mu^\alpha$, gives the synchronous contraction
\begin{equation}\label{eq: uniform frozen contraction}
\mathbb E|Z_t^{\alpha,x}-Z_t^{\alpha,y}|^2\leq e^{-\mathfrak a t}|x-y|^2,\qquad x,y\in\mathbb R^d,
\end{equation}
where $\mathfrak a>0$ independent of $\alpha$. Since $\mu^\alpha$ is invariant for the frozen diffusion, coupling $Z^{\alpha,z}$ with $Z^{\alpha,\xi}$, where $\xi\sim\mu^\alpha$ and both processes share the same noise, yields
\[
\mathcal W_2\big(\mathcal L_t[Z^\alpha],\mu^\alpha\big)^2\leq e^{-\mathfrak a t}\mathcal W_2(\delta_z,\mu^\alpha)^2\leq K_ze^{-\mathfrak a t},
\]
where $K_z:=2|z|^2+2\sup_{\alpha\in\mathcal A_\gamma^\vartheta}\mu^\alpha(|\cdot|^2)<\infty$ by the first part of \eqref{uniform moment}. The frozen $(2+\rho)$-moments are bounded uniformly in time and in $\alpha$ by Lemma \ref{lem: uniform frozen moments}.

For the ordinary occupation measure, we take Lebesgue sampling on $[0,1]$ in the proof of \cite[Proposition 3.1]{du2023empirical}. The mollification and averaged-law estimates in \cite[Eqs. (3.2)--(3.4)]{du2023empirical}, the density-coupling inequality in \cite[Lemma 3.3]{du2023empirical}, and the tail estimate in \cite[Eq. (3.5)]{du2023empirical} apply with the uniform bounds above. In the covariance step corresponding to \cite[Eqs. (3.6)--(3.7)]{du2023empirical}, we retain $\|\nabla\varphi_r\|_\infty=O(r^{-d-1})$ for a fixed rescaled mollifier and combine the semigroup contraction with a bounded-covariance estimate. As verified in the supplementary material \ref{sec: proofs}, this gives, for $0<r<1$, $R>r$, and $t>1$,
\begin{equation}\label{eq: uniform frozen empirical intermediate}
\mathbb E\left[\mathcal W_2\big(\mathcal E_t[Z^\alpha],\mu^\alpha\big)^2\right]\leq C\left[r^2+t^{-1}+(R-r)^{-\rho}+R^{d+2}r^{-d}\left(\frac{1+\log(1/r)}{t}\right)^{1/2}\right].
\end{equation}
Here $C$ depends only on $d$, $\rho$, $\mathfrak a$, $z$, the uniform invariant second-moment and frozen $(2+\rho)$-moment bounds, and the fixed mollifier; it is independent of $\alpha$, $t$, $r$, and $R$.

For $t>1$, choose
\[
r=t^{-\frac{\rho\varepsilon}{2(d+2)(\rho+2)}},\qquad R=2t^{\frac{\varepsilon}{(d+2)(\rho+2)}}.
\]
Then $r^2=t^{-\kappa\varepsilon}$, $(R-r)^{-\rho}\leq t^{-\kappa\varepsilon}$, and the last term in \eqref{eq: uniform frozen empirical intermediate} is bounded by $C_\varepsilon t^{-\kappa\varepsilon}$ because $\varepsilon<1$ and $t^{-(1-\varepsilon)/2}(1+\log t)^{1/2}$ is bounded on $[1,\infty)$. Also $t^{-1}\leq t^{-\kappa\varepsilon}$. The case $t=1$ follows from the same uniform second-moment bounds. This proves \eqref{eq: uniform frozen empirical}; the covariance adjustment and the dependence of constants are detailed in the supplementary material \ref{sec: proofs}.
\end{proof}

\begin{lemma}[Uniform trajectory comparison]\label{lem: uniform trajectory comparison}
Assume the standing conditions of the present section and fix $z\in\mathbb R^d$. Let $\kappa:=\frac{\rho}{(d+2)(\rho+2)}$. For every
\begin{equation}\label{eq: choice uniform comparison exponent}
    0<\bar\varepsilon<\left(1-\frac{\mathfrak b}{\mathfrak a}\right)\wedge\kappa,
\end{equation}
there exists a constant $C_{z,\gamma,\vartheta,\bar\varepsilon}>0$, independent of $\alpha\in\mathcal A_\gamma^\vartheta$, such that, for all $t\geq1$,
\begin{equation}\label{eq: uniform trajectory comparison}
\begin{aligned}
    \mathbb E\big[|X_t^\alpha-Y_t^\alpha|^2\big]&\leq C_{z,\gamma,\vartheta,\bar\varepsilon}\,t^{-\bar\varepsilon}, \text{and }
    \mathbb E\left[\mathcal W_2\big(\mathcal L_t[X^\alpha],\mathcal E_t[Y^\alpha]\big)^2\right]\leq C_{z,\gamma,\vartheta,\bar\varepsilon}\,t^{-\bar\varepsilon}.
\end{aligned}
\end{equation}
\end{lemma}

\begin{proof}
Fix $\bar\varepsilon$ as in \eqref{eq: choice uniform comparison exponent} and choose $\bar\varepsilon<q<\kappa$. For each $\alpha\in\mathcal A_\gamma^\vartheta$, let $\mu^\alpha$ be the invariant measure and couple $X^\alpha$, $Y^\alpha$, and the frozen diffusion $Z^\alpha$ from \eqref{eq: frozen alpha process} with the same Brownian motion and initial condition $z$. By Lemmas \ref{lem: uniform frozen moments} and \ref{lem: uniform frozen empirical}, including the uniform second-moment bound for $0<t\leq1$, there is a constant $F_q>0$, independent of $\alpha$, such that
\begin{equation}\label{eq: frozen empirical used comparison}
\mathbb E\left[\mathcal W_2\big(\mathcal E_t[Z^\alpha],\mu^\alpha\big)^2\right]\leq F_q(1\wedge t^{-q}),\qquad t>0.
\end{equation}
Also, $K_z:=\sup_{\alpha\in\mathcal A_\gamma^\vartheta}\mathcal W_2(\delta_z,\mu^\alpha)^2<\infty$ by the first part of \eqref{uniform moment}. If $\mathfrak b=0$, Assumption \textnormal{(A)} and the synchronous coupling give $X^\alpha=Y^\alpha=Z^\alpha$; the conclusion then follows from frozen contraction and \eqref{eq: frozen empirical used comparison}. Henceforth suppose $\mathfrak b>0$ and choose $\delta>0$ such that $\widehat{\mathfrak b}:=(1+\delta)\mathfrak b<\mathfrak a$ and $\bar\varepsilon<1-\widehat{\mathfrak b}/\mathfrak a$.

Set $h_{\texttt L}^\alpha(t):=\mathbb E|X_t^\alpha-Z_t^\alpha|^2$. By Assumption \textnormal{(A)}, synchronous coupling, and the frozen contraction \eqref{eq: uniform frozen contraction},
\[
\frac{\dd}{\dd t}h_{\texttt L}^\alpha(t)\leq-(\mathfrak a-\widehat{\mathfrak b})h_{\texttt L}^\alpha(t)+\mathfrak b(1+\delta^{-1})K_z e^{-\mathfrak a t}.
\]
Since $h_{\texttt L}^\alpha(0)=0$, integration gives \eqref{eq: uniform X Z estimate} with $\lambda:=\mathfrak a-\widehat{\mathfrak b}>0$ and a constant independent of $\alpha$:
\begin{equation}\label{eq: uniform X Z estimate}
 h_{\texttt L}^\alpha(t)\leq C_{z,\gamma,\vartheta}e^{-\lambda t},\qquad t\geq0.
\end{equation}
Consequently,
\begin{equation}\label{eq: uniform X Z polynomial}
 h_{\texttt L}^\alpha(t)\leq C_{z,\gamma,\vartheta,\bar\varepsilon}t^{-\bar\varepsilon},\qquad t\geq1.
\end{equation}

For $h_{\texttt E}^\alpha(t):=\mathbb E|Y_t^\alpha-Z_t^\alpha|^2$, the occupation-measure coupling
\[
\mathcal W_2\big(\mathcal E_t[Y^\alpha],\mathcal E_t[Z^\alpha]\big)^2\leq\frac1t\int_0^t|Y_s^\alpha-Z_s^\alpha|^2\,\dd s
\]
and Assumption \textnormal{(A)}, together with \eqref{eq: frozen empirical used comparison}, yield
\begin{equation}\label{eq: scalar inequality Y Z}
\frac{\dd}{\dd t}h_{\texttt E}^\alpha(t)\leq-\mathfrak a h_{\texttt E}^\alpha(t)+\widehat{\mathfrak b}\frac1t\int_0^t h_{\texttt E}^\alpha(s)\,\dd s+M_q(1\wedge t^{-q}),
\end{equation}
where $M_q:=\mathfrak b(1+\delta^{-1})F_q$ is independent of $\alpha$. Choose $\bar\varepsilon<r<1-\widehat{\mathfrak b}/\mathfrak a$ and set $r_0:=r\wedge q>\bar\varepsilon$. With Lebesgue time-sampling, $\int_0^1u^{-r}\,\dd u=(1-r)^{-1}<\mathfrak a/\widehat{\mathfrak b}$. Thus the polynomial supersolution argument of \cite[proof of Lemma 4.1, p. 2505]{du2023empirical}, applied with parameters $(\mathfrak a,\widehat{\mathfrak b},r,q,M_q)$ and comparison function $g(t)=t^{-r_0}$, gives
\begin{equation}\label{eq: uniform Y Z estimate}
 h_{\texttt E}^\alpha(t)\leq C_{z,\gamma,\vartheta,\bar\varepsilon}(1\wedge t^{-r_0}),\qquad t\geq0.
\end{equation}
The comparison constant is uniform in $\alpha$, since $\mathfrak a$, $\widehat{\mathfrak b}$, $r$, $q$, and $M_q$ are fixed; the bounded-time estimate is obtained directly from \eqref{eq: scalar inequality Y Z}. In particular,
\begin{equation}\label{eq: uniform Y Z desired exponent}
\mathbb E|Y_t^\alpha-Z_t^\alpha|^2\leq C_{z,\gamma,\vartheta,\bar\varepsilon}t^{-\bar\varepsilon},\qquad t\geq1.
\end{equation}
Combining \eqref{eq: uniform X Z polynomial} and \eqref{eq: uniform Y Z desired exponent} proves the first inequality in \eqref{eq: uniform trajectory comparison}.

For the measure-valued estimate, the triangle inequality, synchronous coupling, and frozen contraction give $\mathcal W_2(\mathcal L_t[X^\alpha],\mu^\alpha)^2\leq2h_{\texttt L}^\alpha(t)+2K_ze^{-\mathfrak a t}\leq Ct^{-\bar\varepsilon}$ for $t\geq1$. Moreover, since $r_0<1$, \eqref{eq: uniform Y Z estimate} and the occupation-measure coupling imply
\[
\mathbb E\left[\mathcal W_2\big(\mathcal E_t[Y^\alpha],\mathcal E_t[Z^\alpha]\big)^2\right]\leq\frac1t\int_0^t h_{\texttt E}^\alpha(s)\,\dd s\leq Ct^{-r_0},\qquad t\geq1.
\]
Together with \eqref{eq: frozen empirical used comparison} and $q,r_0>\bar\varepsilon$, this yields $\mathbb E[\mathcal W_2(\mathcal E_t[Y^\alpha],\mu^\alpha)^2]\leq Ct^{-\bar\varepsilon}$. The second inequality in \eqref{eq: uniform trajectory comparison} follows by the triangle inequality. The scalar comparison and the uniform dependence of its constants are detailed in the supplementary material \ref{sec: proofs}.
\end{proof}

\subsection{Finite-horizon costs and values for fixed-point games}\label{subsec: finite horizon cost comparison}

We shall use the uniform trajectory estimates developed in the previous subsection to tackle the finite horizon McKean--Vlasov fixed point game and its One-Player game (self-interacting) version. The given deterministic population flow affects only the running cost. However, the controlled dynamics continue to depend on their own endogenous laws or self-interactions. We recall that in the sequel, optimization is restricted to the given set $\mathcal A_\gamma^\vartheta$ under the assumptions made in the beginning of this section.

Let $T>0$ and let $\bm\mu=(\mu_t)_{0\leq t\leq T}\in C([0,T];\mathscr P_2(\mathbb R^d))$ be a family of deterministic flow.  The finite horizon problem differs from the ergodic problem, which assumes that the prescribed measure is constant. In the finite horizon problem, there is a time-dependent prescribed flow $(\mu_t)_{t\geq0}$. We will make assumptions on the flow to ensure that the running cost is integrable. Specifically, we will require that $f$ is jointly Borel measurable and will focus on flows which satisfy:
\begin{equation}\label{eq: finite horizon flow integrability}
    \int_0^T|f(0,\mu_t,a_0)|\,\dd t<+\infty,
\end{equation}
for some fixed $a_0\in\mathds A$. We denote this class by $\mathfrak M_T$.

By Assumption \textnormal{(B)}, this condition does not depend on the selection of $a_0$. If $\eta\mapsto f(0,\eta,a_0)$ is continuous in $\mathcal W_2$, then all continuous flows satisfy \eqref{eq: finite horizon flow integrability}. Furthermore, Assumption \textnormal{(B$'$)} implies that $\mathfrak M_T = C([0,T];\mathscr P_2(\mathbb R^d))$.

For $\alpha\in\mathcal A_\gamma^\vartheta$ and $\bm\mu=(\mu_t)_{0\leq t\leq T}\in\mathfrak M_T$, define
\begin{align}
    J_T^{\texttt L}(\alpha,\bm\mu)&:=\mathbb E\left[\int_0^T f\big(X_t^\alpha,\mu_t,\alpha(X_t^\alpha,\mathcal L_t[X^\alpha])\big)\,\dd t\right],\label{eq: finite horizon L cost}\\
    J_T^{\texttt E}(\alpha,\bm\mu)&:=\mathbb E\left[\int_0^T f\big(Y_t^\alpha,\mu_t,\alpha(Y_t^\alpha,\mathcal E_t[Y^\alpha])\big)\,\dd t\right].\label{eq: finite horizon E cost}
\end{align}
Here $X^\alpha$ and $Y^\alpha$ start from the same deterministic initial condition $z$ and are generated by the same feedback $\alpha$. 

By Assumption \textnormal{(B)},
\(\;
    |f(x,\mu,a)|\leq|f(0,\mu,a_0)|+C_f\big(|x|+|a-a_0|\big).
\)\;
Hence \eqref{eq: finite horizon flow integrability}, compactness of $\mathds A$, and the finite-time moment bounds ensure these costs are finite.

\begin{proposition}[Uniform finite-horizon cost comparison]\label{prop: uniform finite horizon cost comparison}
Assume the standing conditions of the present section. Fix $z\in\mathbb R^d$. Let $\bar\varepsilon>0$ satisfy \eqref{eq: choice uniform comparison exponent}. Then there exists a constant $C_{z,\gamma,\vartheta,\bar\varepsilon}>0$, independent of $T\geq1$, $\alpha\in\mathcal A_\gamma^\vartheta$, and $\bm\mu\in\mathfrak M_T$, such that
\begin{equation}\label{eq: uniform finite horizon cost comparison}
    \frac1T\left|J_T^{\texttt L}(\alpha,\bm\mu)-J_T^{\texttt E}(\alpha,\bm\mu)\right|\leq C_{z,\gamma,\vartheta,\bar\varepsilon}\big(T^{-1}+T^{-\bar\varepsilon/2}\big).
\end{equation}
\end{proposition}

\begin{proof}
By Assumption \textnormal{(B)} and the definition of $\mathcal A_\gamma^\vartheta$, there exists $C_{\gamma,\vartheta}>0$  a constant independent of $\alpha\in\mathcal A_\gamma^\vartheta$, such that
\begin{equation}\label{eq: uniform finite horizon cost lip}
    \big|f(x,\mu,\alpha(x,\eta))-f(y,\mu,\alpha(y,\nu))\big|\leq C_{\gamma,\vartheta}\big(|x-y|+\mathcal W_2(\eta,\nu)\big),
\end{equation}
for all $x,y\in\mathbb R^d$ and $\mu,\eta,\nu\in\mathscr P_2(\mathbb R^d)$. In particular, the constant is independent of the prescribed measure $\mu$.

Applying this estimate to the two controlled processes gives
\begin{equation}\label{eq: finite horizon instantaneous comparison}
\begin{aligned}
    &\big|f\big(X_t^\alpha,\mu_t,\alpha(X_t^\alpha,\mathcal L_t[X^\alpha])\big)-f\big(Y_t^\alpha,\mu_t,\alpha(Y_t^\alpha,\mathcal E_t[Y^\alpha])\big)\big|\\
    &\qquad\leq C_{\gamma,\vartheta}\big(|X_t^\alpha-Y_t^\alpha|+\mathcal W_2(\mathcal L_t[X^\alpha],\mathcal E_t[Y^\alpha])\big).
\end{aligned}
\end{equation}

For $t\geq1$, we use Lemma \ref{lem: uniform trajectory comparison}  which, together with Cauchy--Schwarz, yields
\[
    \sup_{\alpha\in\mathcal A_\gamma^\vartheta}\mathbb E\left[|X_t^\alpha-Y_t^\alpha|+\mathcal W_2\big(\mathcal L_t[X^\alpha],\mathcal E_t[Y^\alpha]\big)\right]\leq C_{z,\gamma,\vartheta,\bar\varepsilon}\,t^{-\bar\varepsilon/2}.
\]
Thus the expected absolute difference of the instantaneous running costs in \eqref{eq: finite horizon instantaneous comparison} is bounded by $C_{z,\gamma,\vartheta,\bar\varepsilon}t^{-\bar\varepsilon/2}$, uniformly over $\alpha\in\mathcal A_\gamma^\vartheta$ and $\bm\mu\in\mathfrak M_T$.

For $0<t\leq1$, we use the uniform short-time moment estimate \eqref{eq: uniform short time second moments finite horizon}, established in the preceding subsection. By the triangle inequality for $\mathcal W_2$, applied with the intermediate measure $\delta_0$, we have
\[
\begin{aligned}
\mathcal W_2\big(\mathcal L_t[X^\alpha],\mathcal E_t[Y^\alpha]\big)^2\!\leq\!2\mathcal W_2\big(\mathcal L_t[X^\alpha],\delta_0\big)^2\!\!+\!2\mathcal W_2\big(\mathcal E_t[Y^\alpha],\delta_0\big)^2\! \!=\!2\mathbb E|X_t^\alpha|^2\!+\!\frac2t\!\!\int_0^t \!|Y_s^\alpha|^2\dd s.
\end{aligned}
\]
Taking expectations and using \eqref{eq: uniform short time second moments finite horizon} yields
\[
\mathbb E\left[\mathcal W_2\big(\mathcal L_t[X^\alpha],\mathcal E_t[Y^\alpha]\big)^2\right]\leq2\mathbb E|X_t^\alpha|^2+\frac2t\int_0^t\mathbb E|Y_s^\alpha|^2\,\dd s\leq C_{z,\gamma,\vartheta}.
\]
Moreover, by Cauchy--Schwarz one has
\[
\mathbb E|X_t^\alpha-Y_t^\alpha|\leq\left(\mathbb E|X_t^\alpha-Y_t^\alpha|^2\right)^{1/2}\leq\left(2\mathbb E|X_t^\alpha|^2+2\mathbb E|Y_t^\alpha|^2\right)^{1/2}\leq C_{z,\gamma,\vartheta}.
\]
Similarly, Cauchy--Schwarz gives
\[
\mathbb E\left[\mathcal W_2\big(\mathcal L_t[X^\alpha],\mathcal E_t[Y^\alpha]\big)\right]\leq\left(\mathbb E\left[\mathcal W_2\big(\mathcal L_t[X^\alpha],\mathcal E_t[Y^\alpha]\big)^2\right]\right)^{1/2}\leq C_{z,\gamma,\vartheta}.
\]
Finally, these two estimates together yield
\[
\sup_{\alpha\in\mathcal A_\gamma^\vartheta}\sup_{0<t\leq1}\mathbb E\left[|X_t^\alpha-Y_t^\alpha|+\mathcal W_2\big(\mathcal L_t[X^\alpha],\mathcal E_t[Y^\alpha]\big)\right]\leq C_{z,\gamma,\vartheta}.
\]

Splitting the time integral over $[0,1]$ and $[1,T]$, we obtain, for $T\geq1$,
\[
    \frac1T\left|J_T^{\texttt L}(\alpha,\bm\mu)-J_T^{\texttt E}(\alpha,\bm\mu)\right|\leq\frac{C_{z,\gamma,\vartheta}}T+\frac{C_{z,\gamma,\vartheta,\bar\varepsilon}}T\int_1^Tt^{-\bar\varepsilon/2}\,\dd t.
\]
Since $0<\bar\varepsilon<1$, direct integration gives, for $T\geq1$,
\[
\frac1T\int_1^Tt^{-\bar\varepsilon/2}\,\dd t=\frac{T^{-\bar\varepsilon/2}-T^{-1}}{1-\bar\varepsilon/2}\leq\frac{T^{-\bar\varepsilon/2}}{1-\bar\varepsilon/2}
\]\
which concludes the proof.
\end{proof}

\begin{remark}\label{rem: uniformity finite horizon}
Unlike in the ergodic setting, finite-horizon optimal or approximately optimal feedbacks may depend on $T$. Policy-dependent trajectory estimates therefore do not suffice to compare finite-horizon values or optimizers. The uniform estimate \eqref{eq: uniform finite horizon cost comparison} allows feedbacks to vary with the horizon and motivates restricting the analysis to the fixed class $\mathcal A_\gamma^\vartheta$.
\end{remark}

\begin{corollary}[Finite-horizon value comparison]\label{cor: finite horizon value comparison}
Assume the standing conditions of the present section and fix $z\in\mathbb R^d$. Let $\bar\varepsilon>0$ satisfy \eqref{eq: choice uniform comparison exponent}. For $T\geq1$ and $\bm\mu\in\mathfrak M_T$, define the finite-horizon values restricted to $\mathcal A_\gamma^\vartheta$ by
\begin{equation}\label{eq: finite horizon restricted values}
    v_{T,\gamma,\vartheta}^{\texttt L}(z,\bm\mu):=\inf_{\alpha\in\mathcal A_\gamma^\vartheta}J_T^{\texttt L}(\alpha,\bm\mu),\qquad v_{T,\gamma,\vartheta}^{\texttt E}(z,\bm\mu):=\inf_{\alpha\in\mathcal A_\gamma^\vartheta}J_T^{\texttt E}(\alpha,\bm\mu).
\end{equation}
Then both values are finite, and there exists a constant $C_{z,\gamma,\vartheta,\bar\varepsilon}>0$, independent of $T$ and $\bm\mu$, such that
\begin{equation}\label{eq: finite horizon value comparison}
    \frac1T\left|v_{T,\gamma,\vartheta}^{\texttt L}(z,\bm\mu)-v_{T,\gamma,\vartheta}^{\texttt E}(z,\bm\mu)\right|\leq C_{z,\gamma,\vartheta,\bar\varepsilon}\big(T^{-1}+T^{-\bar\varepsilon/2}\big).
\end{equation}
\end{corollary}

\begin{proof}
For real-valued functions $F$ and $G$ on a common nonempty set $A$, with finite infima,
\begin{equation}\label{eq: infimum comparison elementary}
    \left|\inf_{a\in A}F(a)-\inf_{a\in A}G(a)\right|\leq\sup_{a\in A}|F(a)-G(a)|.
\end{equation}
By \eqref{eq: finite horizon flow integrability}, Assumption \textnormal{(B)}, compactness of $\mathds A$, and the uniform finite-time moment bounds, both costs are bounded in absolute value uniformly over $\mathcal A_\gamma^\vartheta$ for each fixed $T$ and $\bm\mu$. Thus the restricted values are finite, since the feedback class is nonempty.

Applying \eqref{eq: infimum comparison elementary} to $J_T^{\texttt L}(\cdot,\bm\mu)$ and $J_T^{\texttt E}(\cdot,\bm\mu)$, then using Proposition \ref{prop: uniform finite horizon cost comparison}, gives \eqref{eq: finite horizon value comparison}. The convergence is uniform because the bound is independent of $\bm\mu$.
\end{proof}

\subsection{Approximate equilibria for finite-horizon fixed-point games}\label{subsec: finite horizon equilibrium comparison}

We introduce finite horizon versions of the ergodic fixed point games discussed in Section \ref{sec:ergodic mfg}. In this case, the dynamics are still controlled, but the ergodic objectives are replaced by finite horizon costs as defined earlier, and the stationary consistency requirements are replaced by requirements on deterministic flow of probability measures. 

In the law dependent game, with $T>0$ and $\bm\mu=(\mu_t)_{0\leq t\leq T}\in\mathfrak M_T$, where $\mathfrak M_T$ is defined in subsection \ref{subsec: finite horizon cost comparison}, the optimal control uses the cost function $J_T^{\texttt L}(\alpha,\bm\mu)$ and the set of admissible controls $\mathcal A_{\gamma}^{\vartheta}$. For the consistency requirement, we need
\[
    \bm\mu=\big(\mathcal L_t[X^\alpha]\big)_{0\leq t\leq T}.
\]
In the self-interacting game, we use $J_T^{\texttt E}(\alpha,\bm\mu)$ in for the optimization. To address the fact that the occupation measure $\mathcal E_t[Y^\alpha]$ is random, we define a deterministic ``\textit{averaged experience measure}''
\begin{equation}\label{eq: averaged experience measure}
    \overline{\mathcal E}_t[Y^\alpha]:=\mathbb E[\mathcal E_t[Y^\alpha]]=\frac1t\int_0^t\mathcal L_s[Y^\alpha]\,\dd s,\qquad t>0,
\end{equation}
which represents the accumulated experience of a single trajectory (agent) averaged over realizations. Thus, for the consistency condition, we need
\[
    \bm\mu=\big(\overline{\mathcal E}_t[Y^\alpha]\big)_{0\leq t\leq T}.
\]

We set $\overline{\mathcal E}_0[Y^\alpha]:=\delta_z$. To verify continuity of the resulting flow in $\mathcal W_2$, observe that, for $s,t\in[0,T]$,
\[
\overline{\mathcal E}_t[Y^\alpha]=\int_0^1\mathcal L_{tu}[Y^\alpha]\,\dd u,\qquad \overline{\mathcal E}_s[Y^\alpha]=\int_0^1\mathcal L_{su}[Y^\alpha]\,\dd u,
\]
with the convention at $t=0$ or $s=0$. For each $u\in[0,1]$, the joint law $\pi_u:=\operatorname{Law}(Y_{tu}^\alpha,Y_{su}^\alpha)$ is a coupling of $\mathcal L_{tu}[Y^\alpha]$ and $\mathcal L_{su}[Y^\alpha]$. Hence $\pi:=\int_0^1\pi_u\,\dd u$ is a coupling of $\overline{\mathcal E}_t[Y^\alpha]$ and $\overline{\mathcal E}_s[Y^\alpha]$. By the definition of the Wasserstein distance and Fubini's theorem,
\begin{equation}\label{eq: averaged experience flow continuity bound}
\begin{aligned}
\mathcal W_2\big(\overline{\mathcal E}_t[Y^\alpha],\overline{\mathcal E}_s[Y^\alpha]\big)^2&\leq\int_{\mathbb R^d\times\mathbb R^d}|x-y|^2\,\pi(\dd x,\dd y)\;=\int_0^1\mathbb E|Y_{tu}^\alpha-Y_{su}^\alpha|^2\,\dd u.
\end{aligned}
\end{equation}
By the $L^2$-continuity of $Y^\alpha$ and its finite-time second-moment bounds of Lemma \ref{lem: short time moments}, dominated convergence proves continuity for $s>0$. At the origin,
\[
    \mathcal W_2\big(\overline{\mathcal E}_t[Y^\alpha],\delta_z\big)^2=\frac1t\int_0^t\mathbb E|Y_r^\alpha-z|^2\,\dd r\longrightarrow0,\qquad t\downarrow0.
\]
Consequently,
\(
    \big(\overline{\mathcal E}_t[Y^\alpha]\big)_{0\leq t\leq T}\in C\big([0,T];\mathscr P_2(\mathbb R^d)\big).
\) 
The marginal flow $(\mathcal L_t[X^\alpha])_{0\leq t\leq T}$ is likewise continuous in $\mathcal W_2$.

For $T>0$ and $\alpha\in\mathcal A_\gamma^\vartheta$, define the deterministic consistency maps
\begin{equation}\label{eq: finite horizon consistency maps}
    \Psi_T^{\texttt L}(\alpha):=\big(\mathcal L_t[X^\alpha]\big)_{0\leq t\leq T},\qquad \Psi_T^{\texttt E}(\alpha):=\big(\overline{\mathcal E}_t[Y^\alpha]\big)_{0\leq t\leq T}.
\end{equation}
Both maps take values in $C([0,T];\mathscr P_2(\mathbb R^d))$, but their images need not belong to $\mathfrak M_T$. We therefore impose admissibility of the prescribed flow explicitly in the equilibrium definitions.

For $\star\in\{\texttt L,\texttt E\}$ and $\bm\mu\in\mathfrak M_T$, define the finite-horizon best-response correspondence
\begin{equation}\label{eq: finite horizon best response maps}
    \operatorname{BR}_T^\star(\bm\mu):=\argmin_{\alpha\in\mathcal A_\gamma^\vartheta}J_T^\star(\alpha,\bm\mu).
\end{equation}
This set may be empty, since existence of an optimal feedback is not assumed. An \emph{exact equilibrium} of formulation $\star$ is a pair $(\alpha,\bm\mu)\in\mathcal A_\gamma^\vartheta\times\mathfrak M_T$ satisfying
\begin{equation}\label{eq: finite horizon exact equilibrium system}
    \alpha\in\operatorname{BR}_T^\star(\bm\mu),\qquad \bm\mu=\Psi_T^\star(\alpha).
\end{equation}

Unlike in the ergodic framework, the finite-horizon best-response correspondences and consistency maps need not coincide in both law-dependent and self-interacting settings. We instead compare them through uniform cost estimates and a time-averaged distance between flows. 
On $C([0,T];\mathscr P_2(\mathbb R^d))$, define
\begin{equation}\label{eq: finite horizon flow distance}
    \mathbf d_T(\bm\mu,\bm\nu):=\left(\frac1T\int_0^T\mathcal W_2(\mu_t,\nu_t)^2\,\dd t\right)^{1/2}.
\end{equation}

\begin{lemma}[Uniform comparison of the consistency flows]\label{lem: uniform consistency flow comparison}
Assume the standing conditions of the present section and fix $z\in\mathbb R^d$. Let $\bar\varepsilon>0$ satisfy \eqref{eq: choice uniform comparison exponent}. Then there exists a constant $C_{z,\gamma,\vartheta,\bar\varepsilon}>0$, independent of $T\geq1$ and $\alpha\in\mathcal A_\gamma^\vartheta$, such that
\begin{equation}\label{eq: uniform consistency flow comparison}
    \mathbf d_T\big(\Psi_T^{\texttt L}(\alpha),\Psi_T^{\texttt E}(\alpha)\big)\leq C_{z,\gamma,\vartheta,\bar\varepsilon}\big(T^{-1/2}+T^{-\bar\varepsilon/2}\big).
\end{equation}
In particular,
$\;\sup_{\alpha\in\mathcal A_\gamma^\vartheta}\mathbf d_T\big(\Psi_T^{\texttt L}(\alpha),\Psi_T^{\texttt E}(\alpha)\big)\longrightarrow0$, as $T\to+\infty.$
\end{lemma}

\begin{proof}
For $t>0$, we use the averaging inequality \eqref{eq: averaged experience Wasserstein consequence} from Appendix \ref{app: occupation averaged measures} which gives
\[
    \mathcal W_2\big(\mathcal L_t[X^\alpha],\overline{\mathcal E}_t[Y^\alpha]\big)^2\leq\mathbb E\left[\mathcal W_2\big(\mathcal L_t[X^\alpha],\mathcal E_t[Y^\alpha]\big)^2\right].
\]
We can then use Lemma \ref{lem: uniform trajectory comparison} and obtain
\[
    \mathcal W_2\big(\mathcal L_t[X^\alpha],\overline{\mathcal E}_t[Y^\alpha]\big)^2\leq C_{z,\gamma,\vartheta,\bar\varepsilon}\,t^{-\bar\varepsilon},\qquad\forall\, t\geq1, \qquad \forall\,\alpha\in\mathcal A_\gamma^\vartheta
\]

For $0<t\leq1$, we use the uniform short-time moment estimate \eqref{eq: uniform short time second moments finite horizon} which gives
\[
    \mathcal W_2\big(\mathcal L_t[X^\alpha],\overline{\mathcal E}_t[Y^\alpha]\big)^2
    \!\leq2\mathbb E|X_t^\alpha|^2\!+\!2\overline{\mathcal E}_t[Y^\alpha](|\cdot|^2)\!=\!2\mathbb E|X_t^\alpha|^2\!+\!\frac2t\!\int_0^t\mathbb E|Y_s^\alpha|^2\,\dd s\;
   \!\! \leq C_{z,\gamma,\vartheta}.
\]
Therefore, splitting the integral defining $\mathbf d_T$ over $[0,1]$ and $[1,T]$, we obtain, for $T\geq1$,
\[
\begin{aligned}
    \mathbf d_T\big(\Psi_T^{\texttt L}(\alpha),\Psi_T^{\texttt E}(\alpha)\big)^2
    &\leq\frac{C_{z,\gamma,\vartheta}}T+\frac{C_{z,\gamma,\vartheta,\bar\varepsilon}}T\int_1^Tt^{-\bar\varepsilon}\,\dd t\;
    \leq C_{z,\gamma,\vartheta,\bar\varepsilon}\big(T^{-1}+T^{-\bar\varepsilon}\big),
\end{aligned}
\]
where $0<\bar\varepsilon<1$. Taking square roots gives \eqref{eq: uniform consistency flow comparison}. The uniform convergence follows because the constant is independent of $\alpha$.
\end{proof}

We next introduce approximate equilibria. For $\star\in\{\texttt L,\texttt E\}$, $T\geq1$, and $\varepsilon,\delta\geq0$, a pair $(\alpha,\bm\mu)\in\mathcal A_\gamma^\vartheta\times\mathfrak M_T$ is an \emph{$(\varepsilon,\delta)$-equilibrium} of the finite-horizon $\star$-game if
\begin{equation}\label{eq: finite horizon approximate}
    \frac1T J_T^\star(\alpha,\bm\mu)\leq\inf_{\beta\in\mathcal A_\gamma^\vartheta}\frac1T J_T^\star(\beta,\bm\mu)+\varepsilon \quad \text{ and } \quad \mathbf d_T\big(\bm\mu,\Psi_T^\star(\alpha)\big)\leq\delta.
\end{equation} 
Thus, $\varepsilon$ measures the optimality error, while $\delta$ measures the consistency error. In particular, a $(0,0)$-equilibrium is an exact equilibrium within $\mathcal A_\gamma^\vartheta$.

\begin{theorem}[Finite-horizon approximate equilibrium comparison]\label{thm: finite horizon approximate equilibrium comparison}
Assume the standing conditions of the present section and fix $z\in\mathbb R^d$. Let $\bar\varepsilon>0$ satisfy \eqref{eq: choice uniform comparison exponent}. Then there exist constants $C_1,C_2>0$, depending only on $z,\gamma,\vartheta,\bar\varepsilon$ and the standing structural constants, but independent of $T\geq1$, such that:\\
\textbullet\, If $(\alpha_T,\bm\mu_T)\in\mathcal A_\gamma^\vartheta\times\mathfrak M_T$ is an exact equilibrium of the finite-horizon law-dependent game, then it is an $(\varepsilon_T,\delta_T)$-equilibrium of the finite-horizon self-interacting game, where
\begin{equation}\label{eq: finite horizon equilibrium errors}
    \varepsilon_T:=C_1\big(T^{-1}+T^{-\bar\varepsilon/2}\big),\qquad \delta_T:=C_2\big(T^{-1/2}+T^{-\bar\varepsilon/2}\big).
\end{equation}
\textbullet\, Conversely, every exact equilibrium $(\alpha_T,\bm\mu_T)\in\mathcal A_\gamma^\vartheta\times\mathfrak M_T$ of the finite-horizon self-interacting game is an $(\varepsilon_T,\delta_T)$-equilibrium of the finite-horizon law-dependent game, with the same error bounds. In particular, $\varepsilon_T\to0$ and $\delta_T\to0$ as $T\to+\infty$.
\end{theorem}

\begin{proof}
Suppose  $(\alpha_T,\bm\mu_T)$ is an exact equilibrium of the law-dependent game. Then
\begin{equation*}
    J_T^{\texttt L}(\alpha_T,\bm\mu_T)\leq J_T^{\texttt L}(\beta,\bm\mu_T)\quad\forall\, \beta\in\mathcal A_\gamma^\vartheta, \quad \text{ and } \quad \bm\mu_T=\Psi_T^{\texttt L}(\alpha_T).
\end{equation*}

Let $\eta_T:=C_{\mathrm{cost}}\big(T^{-1}+T^{-\bar\varepsilon/2}\big)$, where $C_{\mathrm{cost}}$ is the constant in Proposition \ref{prop: uniform finite horizon cost comparison}. Since $\bm\mu_T\in\mathfrak M_T$, the uniform cost estimate of that proposition applies to every feedback in $\mathcal A_\gamma^\vartheta$. 
Thus, for any $\beta\in\mathcal A_\gamma^\vartheta$,
\begin{align*}
    \frac1T J_T^{\texttt E}(\alpha_T,\bm\mu_T)
    &\leq\frac1T J_T^{\texttt L}(\alpha_T,\bm\mu_T)+\eta_T\;
    \leq\frac1T J_T^{\texttt L}(\beta,\bm\mu_T)+\eta_T\;
    \leq\frac1T J_T^{\texttt E}(\beta,\bm\mu_T)+2\eta_T.
\end{align*}
Taking the infimum over $\beta$ proves the first part of \eqref{eq: finite horizon approximate} for the self-interacting game with $\varepsilon_T=2\eta_T$.

To verify consistency, we use Lemma \ref{lem: uniform consistency flow comparison} which gives
\[
    \mathbf d_T\big(\bm\mu_T,\Psi_T^{\texttt E}(\alpha_T)\big)
    =\mathbf d_T\big(\Psi_T^{\texttt L}(\alpha_T),\Psi_T^{\texttt E}(\alpha_T)\big)
    \leq C_{\mathrm{cons}}\big(T^{-1/2}+T^{-\bar\varepsilon/2}\big),
\]
where $C_{\mathrm{cons}}$ is independent of $T$. Hence the second part of \eqref{eq: finite horizon approximate} holds with $\delta_T=C_{\mathrm{cons}}\big(T^{-1/2}+T^{-\bar\varepsilon/2}\big)$.

Taking $C_1:=2C_{\mathrm{cost}}$ and $C_2:=C_{\mathrm{cons}}$ proves the first assertion. The converse follows by interchanging the roles of $\texttt L$ and $\texttt E$, using the symmetry of the cost comparison and the distance $\mathbf d_T$.
\end{proof}

\subsection{Finite-horizon Mean-Field Control comparisons}\label{subsec: finite horizon mfc}

We will address the finite-horizon versions of the Mean-Field Control and Path-Dependent Control problems presented in Section \ref{sec:ergodic mkv}. As in the ergodic setting, the measure that appears in the running cost is produced by the controlled process itself and thus varies. Therefore, we shall assume in this subsection that Assumption \textnormal{(B$'$)} holds.

For $T>0$ and $\alpha\in\mathcal A_\gamma^\vartheta$, define
\begin{align}
    J_{\mathrm C,T}^{\texttt L}(\alpha)&:=\mathbb E\left[\int_0^T f\big(X_t^\alpha,\mathcal L_t[X^\alpha],\alpha(X_t^\alpha,\mathcal L_t[X^\alpha])\big)\,\dd t\right],\label{eq: finite horizon mfc L cost}\\
    J_{\mathrm C,T}^{\texttt E}(\alpha)&:=\mathbb E\left[\int_0^T f\big(Y_t^\alpha,\mathcal E_t[Y^\alpha],\alpha(Y_t^\alpha,\mathcal E_t[Y^\alpha])\big)\,\dd t\right].\label{eq: finite horizon mfc E cost}
\end{align}
Here $X^\alpha$ and $Y^\alpha$ start from the same deterministic initial condition $z$ and are generated by the same feedback $\alpha$.

In the first step, we compare the finite horizon costs uniformly over $\mathcal A_\gamma^\vartheta$. Then, we make a comparison of values and show how to approximate optimality of feedbacks transfers from one problem to the other.

\begin{proposition}[Uniform finite-horizon control-cost comparison]\label{prop: uniform finite horizon mfc cost comparison}
Assume the standing conditions of the present section and Assumption \textnormal{(B$'$)}, and fix $z\in\mathbb R^d$. Let $\bar\varepsilon>0$ satisfy \eqref{eq: choice uniform comparison exponent}. Then there exists a constant $C_{z,\gamma,\vartheta,\bar\varepsilon}>0$, independent of $T\geq1$ and $\alpha\in\mathcal A_\gamma^\vartheta$, such that
\begin{equation}\label{eq: uniform finite horizon mfc comparison}
    \frac1T\left|J_{\mathrm C,T}^{\texttt L}(\alpha)-J_{\mathrm C,T}^{\texttt E}(\alpha)\right|\leq C_{z,\gamma,\vartheta,\bar\varepsilon}\big(T^{-1}+T^{-\bar\varepsilon/2}\big).
\end{equation}
\end{proposition}

\begin{proof}
By Assumption \textnormal{(B$'$)} and the uniform estimate \eqref{eq: uniform finite horizon cost lip}, there exists a constant $\widetilde C_{\gamma,\vartheta}>0$, independent of $\alpha\in\mathcal A_\gamma^\vartheta$, such that
\begin{equation}\label{eq: uniform mfc running cost lip}
    \big|f(x,\eta,\alpha(x,\eta))-f(y,\nu,\alpha(y,\nu))\big|\leq\widetilde C_{\gamma,\vartheta}\big(|x-y|+\mathcal W_2(\eta,\nu)\big)
\end{equation}
for every $x,y\in\mathbb R^d$ and $\eta,\nu\in\mathscr P_2(\mathbb R^d)$. Indeed, adding and subtracting $f(y,\eta,\alpha(y,\nu))$ shows that one may take $\widetilde C_{\gamma,\vartheta}=C_{\gamma,\vartheta}+\widetilde C_f$, where $C_{\gamma,\vartheta}$ is the constant in \eqref{eq: uniform finite horizon cost lip}.

Applying \eqref{eq: uniform mfc running cost lip} to the two controlled processes gives
\begin{equation}\label{eq: finite horizon mfc instantaneous comparison}
\begin{aligned}
    &\big|f\big(X_t^\alpha,\mathcal L_t[X^\alpha],\alpha(X_t^\alpha,\mathcal L_t[X^\alpha])\big)-f\big(Y_t^\alpha,\mathcal E_t[Y^\alpha],\alpha(Y_t^\alpha,\mathcal E_t[Y^\alpha])\big)\big|\\
    &\qquad\leq\widetilde C_{\gamma,\vartheta}\big(|X_t^\alpha-Y_t^\alpha|+\mathcal W_2(\mathcal L_t[X^\alpha],\mathcal E_t[Y^\alpha])\big).
\end{aligned}
\end{equation}

The individual costs are integrable. Indeed, Assumptions \textnormal{(B)} and \textnormal{(B$'$)} imply, for any fixed $a_0\in\mathds A$,
\[
    |f(x,\eta,a)|\leq|f(0,\delta_0,a_0)|+C_f|x|+\widetilde C_f\mathcal W_2(\eta,\delta_0)+C_f|a-a_0|.
\]
Together with the compactness of $\mathds A$ and the finite-time moment bounds, this justifies Fubini's theorem for both costs.

Applying Lemma \ref{lem: uniform trajectory comparison} for $t\geq1$ and \eqref{eq: uniform short time second moments finite horizon} for $0<t\leq1$, we obtain, as in the proof of Proposition \ref{prop: uniform finite horizon cost comparison},
\[
    \frac1T\left|J_{\mathrm C,T}^{\texttt L}(\alpha)-J_{\mathrm C,T}^{\texttt E}(\alpha)\right|\leq\frac{C_{z,\gamma,\vartheta}}T+\frac{C_{z,\gamma,\vartheta,\bar\varepsilon}}T\int_1^Tt^{-\bar\varepsilon/2}\,\dd t.
\]
Since $0<\bar\varepsilon<1$, the right-hand side is bounded by $C_{z,\gamma,\vartheta,\bar\varepsilon}(T^{-1}+T^{-\bar\varepsilon/2})$, uniformly over $\alpha\in\mathcal A_\gamma^\vartheta$. This proves \eqref{eq: uniform finite horizon mfc comparison}.
\end{proof}

\begin{corollary}[Finite-horizon MFC value comparison]
\label{cor: finite horizon mfc value comparison}
In the situation of Proposition \ref{prop: uniform finite horizon mfc cost comparison}, define
\[
    v_{\mathrm C,T,\gamma,\vartheta}^{\texttt L}(z):=\inf_{\alpha\in\mathcal A_\gamma^\vartheta}J_{\mathrm C,T}^{\texttt L}(\alpha),\qquad v_{\mathrm C,T,\gamma,\vartheta}^{\texttt E}(z):=\inf_{\alpha\in\mathcal A_\gamma^\vartheta}J_{\mathrm C,T}^{\texttt E}(\alpha).
\]
Then both values are finite for every $T>0$, and there exists a constant $C_{z,\gamma,\vartheta,\bar\varepsilon}>0$, independent of $T\geq1$, such that
\begin{equation}\label{eq: finite horizon mfc value rate}
    \frac1T\left|v_{\mathrm C,T,\gamma,\vartheta}^{\texttt L}(z)-v_{\mathrm C,T,\gamma,\vartheta}^{\texttt E}(z)\right|\leq C_{z,\gamma,\vartheta,\bar\varepsilon}\big(T^{-1}+T^{-\bar\varepsilon/2}\big).
\end{equation}
\end{corollary}

\begin{proof}
Assumption \textnormal{(B$'$)}, compactness of $\mathds A$, and the uniform finite-time moment bounds imply that both costs are bounded in absolute value uniformly over $\mathcal A_\gamma^\vartheta$ for each fixed $T>0$. Since the feedback class is nonempty, both values are finite.

Applying \eqref{eq: infimum comparison elementary} and Proposition \ref{prop: uniform finite horizon mfc cost comparison} give \eqref{eq: finite horizon mfc value rate}, with the same constant as in the proposition. The asserted convergence follows immediately.
\end{proof}

The preceding results allow for a quantitative transfer of optimal and near-optimal feedback controls from one finite-horizon control problem to the other as we show next.

\begin{theorem}[Finite-horizon approximation of optimal feedbacks]\label{thm: finite horizon approximate minimizers}
Assume the standing conditions of the present section and Assumption \textnormal{(B$'$)}, and fix $z\in\mathbb R^d$. Let $\bar\varepsilon>0$ satisfy \eqref{eq: choice uniform comparison exponent}. Then there exists a constant $C_{z,\gamma,\vartheta,\bar\varepsilon}>0$, independent of $T\geq1$, such that the following holds.

If $\alpha_T\in\mathcal A_\gamma^\vartheta$ is an optimal feedback for the finite-horizon Mean-Field Control problem, then it is an $\varepsilon_T$-optimal feedback for the normalized finite-horizon Path-Dependent Control problem, where
\begin{equation}\label{eq: finite horizon approximate minimizer error}
    \varepsilon_T:=C_{z,\gamma,\vartheta,\bar\varepsilon}\big(T^{-1}+T^{-\bar\varepsilon/2}\big), \quad \text{that is} \quad  \frac1T J_{\mathrm C,T}^{\texttt E}(\alpha_T)\leq\frac1T v_{\mathrm C,T,\gamma,\vartheta}^{\texttt E}(z)+\varepsilon_T.
\end{equation}

Conversely, every optimal feedback for the finite-horizon Path-Dependent Control problem is an $\varepsilon_T$-optimal feedback for the normalized finite-horizon Mean-Field Control problem, with the same error bound.

More generally, if $\alpha_T$ is $\eta_T$-optimal for either normalized control problem, with $\eta_T\geq0$, then it is $(\eta_T+\varepsilon_T)$-optimal for the other. In particular, $\varepsilon_T\to0$ as $T\to+\infty$.
\end{theorem}

\begin{proof}
Let us define $r_T:=C_{\mathrm{cost}}\big(T^{-1}+T^{-\bar\varepsilon/2}\big)$ where 
$C_{\mathrm{cost}}>0$ is the constant in Proposition \ref{prop: uniform finite horizon mfc cost comparison}. 
The proposition and Corollary \ref{cor: finite horizon mfc value comparison} give
\[
    \sup_{\alpha\in\mathcal A_\gamma^\vartheta}\frac1T\left|J_{\mathrm C,T}^{\texttt L}(\alpha)-J_{\mathrm C,T}^{\texttt E}(\alpha)\right|\leq r_T,\qquad \frac1T\left|v_{\mathrm C,T,\gamma,\vartheta}^{\texttt L}(z)-v_{\mathrm C,T,\gamma,\vartheta}^{\texttt E}(z)\right|\leq r_T.
\]
Let  $\alpha_T\in\mathcal A_\gamma^\vartheta$ be $\eta_T$-optimal for the normalized Mean-Field Control problem. Then
\[
\begin{aligned}
    \frac1T J_{\mathrm C,T}^{\texttt E}(\alpha_T)
    \leq\frac1T J_{\mathrm C,T}^{\texttt L}(\alpha_T)\!+\!r_T
    &\leq\frac1T v_{\mathrm C,T,\gamma,\vartheta}^{\texttt L}(z)\!+\eta_T\!+\!r_T
    \leq\frac1T v_{\mathrm C,T,\gamma,\vartheta}^{\texttt E}(z)\!+\!\eta_T\!+\!2r_T.
\end{aligned}
\]
Thus $\alpha_T$ is $(\eta_T+2r_T)$-optimal for the normalized Path-Dependent Control problem. The converse follows by symmetry, while the assertions for exact optimal feedbacks correspond to $\eta_T=0$. Taking $C_{z,\gamma,\vartheta,\bar\varepsilon}:=2C_{\mathrm{cost}}$ in \eqref{eq: finite horizon approximate minimizer error} completes the proof.
\end{proof}

\begin{center}
    \textsc{Conclusion}
\end{center}

We established a population-time correspondence between McKean--Vlasov games and control problems and their self-interacting counterparts. With appropriate dissipativity and regularity, the ergodic formulations possess the same costs, value functions, and optimal feedback sets. Additionally, they exhibit the same equilibrium conditions. In the finte-hoizon setting, we also established quantitative transfers of approximate game equilibrium and optimal feedback controls. These results demonstrate that when a system exhibits sufficient dissipativity, information regarding a population may ultimately be substituted for experience gained over time by an individual on a single trajectory. This perspective also suggests further questions about individual learning and its relation to population-level behavior in socio-economic contexts.

\appendix

\section{On occupation and averaged measures}\label{app: occupation averaged measures}

\subsection{Occupation measures as Wasserstein-valued random variables}

In the first part, we record a continuity property of the occupation-measure map. In particular, this verifies that the occupation measures appearing in the self-interacting dynamics can be regarded as $\mathscr P_2(\mathbb R^d)$-valued random variables. The proofs of both lemmas in this appendix are given in the supplementary material \ref{sec: proofs}.

\begin{lemma}[Continuity of the occupation-measure map]\label{lem: continuity occupation measure map}
Fix $t>0$ and define
\[
    \mathfrak E_t:C([0,t];\mathbb R^d)\longrightarrow\mathscr P_2(\mathbb R^d),\qquad \mathfrak E_t(y):=\frac1t\int_0^t\delta_{y_s}\,\dd s.
\]
Then, for every $y,\widetilde y\in C([0,t];\mathbb R^d)$, it holds
\begin{equation}\label{eq: occupation measure map Lipschitz}
    \mathcal W_2\big(\mathfrak E_t(y),\mathfrak E_t(\widetilde y)\big)^2\leq\frac1t\int_0^t|y_s-\widetilde y_s|^2\,\dd s\leq\|y-\widetilde y\|_\infty^2.
\end{equation}
In particular, $\mathfrak E_t$ is $1$-Lipschitz from $C([0,t];\mathbb R^d)$, endowed with the uniform norm, into $\mathscr P_2(\mathbb R^d)$, endowed with the $\mathcal W_2$ distance.
\end{lemma}

Consequently, if $Y$ is an $\mathbb R^d$-valued stochastic process with continuous sample paths, its path restriction $(Y_s)_{0\leq s\leq t}$ can be regarded as a measurable $C([0,t];\mathbb R^d)$-valued random variable. Since $\mathfrak E_t$ is continuous,
\(\;
    \mathcal E_t[Y]=\mathfrak E_t\big((Y_s)_{0\leq s\leq t}\big)
\)\; 
is a $\mathscr P_2(\mathbb R^d)$-valued random variable. Moreover,
\begin{equation}\label{eq: occupation random measure second moment}
    \mathcal E_t[Y](|\cdot|^2)=\frac1t\int_0^t|Y_s|^2\,\dd s.
\end{equation}
So when $\sup_{0\leq s\leq t}\mathbb E|Y_s|^2<+\infty$, Tonelli gives
\(
    \mathbb E\big[\mathcal E_t[Y](|\cdot|^2)\big]=\frac1t\int_0^t\mathbb E|Y_s|^2\,\dd s<+\infty.
\)

\subsection{Wasserstein distance to an averaged random measure}

We  establish the averaging inequality used to pass from the random occupation measure of the self-interacting dynamics to its deterministic mean. The mean of a random probability measure is understood in the sense of averaging its evaluations on Borel sets. 

\begin{lemma}[Wasserstein distance to a mean random measure]\label{lem: Wasserstein averaged random measure}
Let $\boldsymbol\nu$ be a $\mathscr P_2(\mathbb R^d)$-valued random variable satisfying
\begin{equation}\label{eq: random measure second moment}
    \mathbb E\big[\boldsymbol\nu(|\cdot|^2)\big]<+\infty.
\end{equation}
Define its mean probability measure $\overline{\boldsymbol\nu}:=\mathbb E[\boldsymbol\nu]$ by
\begin{equation}\label{eq: mean random probability measure}
\overline{\boldsymbol\nu}(A):=\mathbb E[\boldsymbol\nu(A)],\qquad A\in\mathscr B(\mathbb R^d).
\end{equation}
Then $\overline{\boldsymbol\nu}\in\mathscr P_2(\mathbb R^d)$ and, for every $\mu\in\mathscr P_2(\mathbb R^d)$,
\begin{equation}\label{eq: Wasserstein averaged random measure}
    \mathcal W_2\big(\mu,\overline{\boldsymbol\nu}\big)^2\leq\mathbb E\big[\mathcal W_2(\mu,\boldsymbol\nu)^2\big].
\end{equation}
\end{lemma}

We finally specialize these results to the self-interacting dynamics. Fix $t>0$ and $\alpha\in\mathcal A_{\mathrm{adm}}$. By Definition \ref{def: admissible controls}, the feedback-composed coefficients
\[
    b^\alpha(x,\eta)=b(x,\eta,\alpha(x,\eta)),\qquad \sigma^\alpha(x,\eta)=\sigma(x,\eta,\alpha(x,\eta))
\]
satisfy Assumption \textnormal{(A)}. Hence the finite-time moment estimates for the corresponding self-interacting dynamics give
\(\;
    \sup_{0\leq s\leq t}\mathbb E|Y_s^\alpha|^2<+\infty.
\)\;
On the other hand, Lemma \ref{lem: continuity occupation measure map} shows that $\mathcal E_t[Y^\alpha]$ is a $\mathscr P_2(\mathbb R^d)$-valued random variable. Therefore, by \eqref{eq: occupation random measure second moment},  \vspace*{-0.5em}
\[
    \mathbb E\big[\mathcal E_t[Y^\alpha](|\cdot|^2)\big]=\frac1t\int_0^t\mathbb E|Y_s^\alpha|^2\,\dd s<+\infty,
\]
so the integrability assumption of Lemma \ref{lem: Wasserstein averaged random measure} is satisfied.

It remains to identify the mean measure $\mathbb E[\mathcal E_t[Y^\alpha]]$. Let $h:\mathbb R^d\to\mathbb R$ be bounded and Borel measurable. By the definition of the mean measure and Fubini's theorem,  \vspace*{-0.5em}
\begin{align*}
    & \int_{\mathbb R^d}h(x)\,\mathbb E[\mathcal E_t[Y^\alpha]](\dd x)
    =\mathbb E\left[\int_{\mathbb R^d}h(x)\,\mathcal E_t[Y^\alpha](\dd x)\right]\\
    &\qquad =\mathbb E\left[\frac1t\int_0^t h(Y_s^\alpha)\,\dd s\right]
    =\frac1t\int_0^t\mathbb E[h(Y_s^\alpha)]\,\dd s\;
    =\frac1t\int_0^t\int_{\mathbb R^d}h(x)\,\mathcal L_s[Y^\alpha](\dd x)\,\dd s.
\end{align*}
The last expression is the integral of $h$ with respect to the deterministic probability measure $\frac1t\int_0^t\mathcal L_s[Y^\alpha]\,\dd s$. Since the identity holds for every bounded Borel function $h$, the measures coincide. Consequently, we have  \vspace*{-0.5em}
\[
    \mathbb E[\mathcal E_t[Y^\alpha]]=\frac1t\int_0^t\mathcal L_s[Y^\alpha]\,\dd s=\overline{\mathcal E}_t[Y^\alpha].
\]

Lemma \ref{lem: Wasserstein averaged random measure}, applied with $\mu=\mathcal L_t[X^\alpha]$ and $\boldsymbol\nu=\mathcal E_t[Y^\alpha]$, therefore yields  \vspace*{-0.5em}
\begin{equation}\label{eq: averaged experience Wasserstein consequence}
    \mathcal W_2\big(\mathcal L_t[X^\alpha],\overline{\mathcal E}_t[Y^\alpha]\big)^2\leq\mathbb E\left[\mathcal W_2\big(\mathcal L_t[X^\alpha],\mathcal E_t[Y^\alpha]\big)^2\right].
\end{equation}
This estimate holds for every fixed admissible feedback. When $\alpha$ ranges over the globally Lipschitz class $\mathcal A_\gamma^\vartheta$, the uniform trajectory comparison of Lemma \ref{lem: uniform trajectory comparison} gives the uniform bound used in Section \ref{subsec: finite horizon equilibrium comparison}.

\section{Supplementary Material: Proofs}\label{sec: proofs}

\begin{proof}[\textbf{Proof of Proposition \ref{prop: estimate trajectories}}]
Let $\mu$ denote the unique invariant probability measure of the McKean--Vlasov dynamics \eqref{MkV SDE}, and introduce the auxiliary frozen process
\[
    \dd Z_t=b(Z_t,\mu)\,\dd t+\sigma(Z_t,\mu)\,\dd W_t,\qquad Z_0=z,
\]
driven by the same Brownian motion as $X$ and $Y$.

The estimates corresponding to \cite[(4.1) and (4.3)]{du2023empirical} yield, for every admissible $\varepsilon$,
\[
    \mathbb E|X_t-Z_t|^2\leq C_1t^{-\varepsilon},\qquad \mathbb E|Y_t-Z_t|^2\leq C_2t^{-\varepsilon},\qquad t\geq1.
\]
Therefore,
\(\;
    \mathbb E|X_t-Y_t|^2\leq2\mathbb E|X_t-Z_t|^2+2\mathbb E|Y_t-Z_t|^2\leq Ct^{-\varepsilon},
\)\;
which proves the first estimate in \eqref{eq: estimate trajectories}. For the second estimate, the triangle inequality for $\mathcal W_2$ gives
\[
    \mathcal W_2\big(\mathcal L_t[X],\mathcal E_t[Y]\big)^2\leq2\mathcal W_2\big(\mathcal L_t[X],\mu\big)^2+2\mathcal W_2\big(\mathcal E_t[Y],\mu\big)^2.
\]
Taking expectations and applying \eqref{eq: conv law to mu} and \eqref{eq: conv path to mu} yields the desired bound.
\end{proof}

\begin{proof}[\textbf{Proof of Lemma \ref{lem: short time moments}}]
Fix $T_0>0$. Throughout the proof, $C$ denotes a positive constant that may change from line to line but depends only on the structural constants in Assumptions \textnormal{(A)}.

We first consider the McKean--Vlasov process $X$. Applying It\^o's formula to $|X_t|^2$, using standard localization, and then applying the Lyapunov inequality \eqref{eq: assumption 1}, we obtain
\[
    \mathbb E|X_t|^2\leq|z|^2+\int_0^t\left[-\mathfrak c\,\mathbb E|X_s|^2+K\big(1+\rho+\mathbb E|X_s|^2\big)\right]\dd s.
\]
Dropping the nonpositive term and applying Gr\"onwall's inequality gives
\[
    \sup_{0\leq t\leq T_0}\mathbb E|X_t|^2\leq\big(|z|^2+K(1+\rho)T_0\big)e^{KT_0}.
\]

For the self-interacting process $Y$, the occupation-measure identity and Tonelli's theorem yield
\[
    \mathbb E\big[\mathcal E_s[Y](|\cdot|^2)\big]=\frac1s\int_0^s\mathbb E|Y_r|^2\,\dd r,\qquad s>0.
\]
Applying It\^o's formula and \eqref{eq: assumption 1}, again with standard localization, we obtain
\[
    \mathbb E|Y_t|^2\leq|z|^2+K(1+\rho)t+K\int_0^t\frac1s\int_0^s\mathbb E|Y_r|^2\,\dd r\,\dd s.
\]
Set $M(t):=\sup_{0\leq r\leq t}\mathbb E|Y_r|^2$. Since $t^{-1}\int_0^t\mathbb E|Y_r|^2\,\dd r\leq M(t)$, it follows that
\[
    M(t)\leq|z|^2+K(1+\rho)t+K\int_0^tM(s)\,\dd s.
\]
Gr\"onwall's inequality therefore gives
\[
    \sup_{0\leq t\leq T_0}\mathbb E|Y_t|^2\leq\big(|z|^2+K(1+\rho)T_0\big)e^{KT_0}.
\]
Together with the estimate for $X$, this proves \eqref{eq: short time moments}.

It remains to establish the integrability assertion. By Cauchy--Schwarz and \eqref{eq: short time moments}, for $0\leq t\leq1$,
\[
    \mathbb E|X_t-Y_t|\leq\big(2\mathbb E|X_t|^2+2\mathbb E|Y_t|^2\big)^{1/2}\leq C_{1,z}.
\]
Moreover, the triangle inequality  for $\mathcal W_2$ applied through $\delta_0$ gives, for $t>0$,
\[
    \mathcal W_2\big(\mathcal L_t[X],\mathcal E_t[Y]\big)^2\leq2\mathbb E|X_t|^2+\frac2t\int_0^t|Y_s|^2\,\dd s.
\]
Taking expectations, applying Tonelli's theorem, and using \eqref{eq: short time moments}, we obtain
\[
    \mathbb E\left[\mathcal W_2\big(\mathcal L_t[X],\mathcal E_t[Y]\big)^2\right]\leq2\mathbb E|X_t|^2+\frac2t\int_0^t\mathbb E|Y_s|^2\,\dd s\leq C_{1,z},\qquad 0<t\leq1.
\]
A final application of Cauchy--Schwarz yields
\[
    \int_0^1\mathbb E\left[|X_t-Y_t|+\mathcal W_2\big(\mathcal L_t[X],\mathcal E_t[Y]\big)\right]\dd t<+\infty.
\]
\end{proof}

\begin{proof}[\textbf{Proof of Lemma \ref{lem: conditions ok}}]
Fix $\alpha\in\mathcal A_\gamma^\vartheta$ and write
\[
    \bar b(x,\mu):=
    b_0(x,\mu)+b_1(x,\alpha(x,\mu)),\qquad \bar\sigma(x,\mu):=
    \sigma_0(x,\mu)+\sigma_1(x,\alpha(x,\mu)).
\]
We verify successively the regularity and the two inequalities in Assumptions (A).

\noindent\textit{\textbf{Step 1} (Continuity and boundedness on bounded sets).} 
Since $\alpha\in\mathcal A_\gamma^\vartheta$, \eqref{eq: lip feedback} implies that $\alpha$ is continuous on $\mathbb R^d\times\mathscr P_2(\mathbb R^d)$. By Assumption (C)(i), $b_0$ and $\sigma_0$ are continuous, while Assumption (C)(ii) implies that $b_1$ and $\sigma_1$ are continuous. Hence $\bar b$ and $\bar\sigma$ are continuous.

Moreover, let $D\subset\mathbb R^d\times\mathscr P_2(\mathbb R^d)$ be bounded. Since $\mathds A$ is compact, the set $\{(x,\alpha(x,\mu)):(x,\mu)\in D\}$ is bounded in $\mathbb R^d\times\mathds A$. Thus, using Assumption (C)(i) and the linear-growth estimates \eqref{eq: b1 sigma1 growth}, both $\bar b$ and $\bar\sigma$ are bounded on $D$.

\noindent\textit{\textbf{Step 2} (Verification of \eqref{eq: assumption 1}).} 
Using $\|A+B\|^2\leq2\|A\|^2+2\|B\|^2$, we obtain
\begin{align*}
    2\langle\bar b(x,\mu),x\rangle+(1+\rho_0)\|\bar\sigma(x,\mu)\|^2
    \leq&2\langle b_0(x,\mu),x\rangle+2(1+\rho_0)\|\sigma_0(x,\mu)\|^2\\
    &+2\langle b_1(x,\alpha(x,\mu)),x\rangle\!+2(1+\!\rho_0)\|\sigma_1(x,\alpha(x,\mu))\|^2.
\end{align*}
By \eqref{eq: D coercivity},
\begin{align*}
    2\langle\bar b(x,\mu),x\rangle+(1+\rho_0)\|\bar\sigma(x,\mu)\|^2
    &\leq-\mathfrak c_0|x|^2+K_0\big[1+\rho_0+\mu(|\cdot|^2)\big]\\
    &\;+2\langle b_1(x,\alpha(x,\mu)),x\rangle+2(1+\rho_0)\|\sigma_1(x,\alpha(x,\mu))\|^2.
\end{align*}
Recall from \eqref{eq: b1 sigma1 growth} that $|b_1(x,a)|\leq C_b|x|+M_b$ and $\|\sigma_1(x,a)\|\leq C_\sigma|x|+M_\sigma$. Therefore,
\[
    2\langle b_1(x,\alpha(x,\mu)),x\rangle\leq2C_b|x|^2+2M_b|x|\leq(2C_b+1)|x|^2+M_b^2.
\]
Similarly,
\[
    2(1+\rho_0)\|\sigma_1(x,\alpha(x,\mu))\|^2\leq4(1+\rho_0)C_\sigma^2|x|^2+4(1+\rho_0)M_\sigma^2.
\]
It follows that
\[
    2\langle\bar b(x,\mu),x\rangle+(1+\rho_0)\|\bar\sigma(x,\mu)\|^2\leq-\mathfrak c|x|^2+\widetilde K\big[1+\rho_0+\mu(|\cdot|^2)\big],
\]
for a suitable constant $\widetilde K>0$, where $\mathfrak c$ is defined in \eqref{eq: resulting abc}. The second inequality in \eqref{eq: condition} ensures that $\mathfrak c>0$. Hence \eqref{eq: assumption 1} is satisfied with $\rho=\rho_0$.

\noindent\textit{\textbf{Step 3} (Verification of \eqref{eq: assumption 2}).} 
Let $x,y\in\mathbb R^d$ and $\mu,\nu\in\mathscr P_2(\mathbb R^d)$, and set $a:=\alpha(x,\mu)$ and $\tilde a:=\alpha(y,\nu)$. Using again $\|A+B\|^2\leq2\|A\|^2+2\|B\|^2$, we obtain
\begin{align*}
    &2\langle\bar b(x,\mu)-\bar b(y,\nu),x-y\rangle+\|\bar\sigma(x,\mu)-\bar\sigma(y,\nu)\|^2\\
    &\quad\leq2\langle b_0(x,\mu)-b_0(y,\nu),x-y\rangle+2\|\sigma_0(x,\mu)-\sigma_0(y,\nu)\|^2\\
    &\qquad+2\langle b_1(x,a)-b_1(y,\tilde a),x-y\rangle+2\|\sigma_1(x,a)-\sigma_1(y,\tilde a)\|^2.
\end{align*}
By \eqref{eq: D dissipativity} and \eqref{eq: D control lip},
\begin{align*}
    &2\langle\bar b(x,\mu)-\bar b(y,\nu),x-y\rangle+\|\bar\sigma(x,\mu)-\bar\sigma(y,\nu)\|^2\\
    &\quad\leq-\mathfrak a_0|x-y|^2+\mathfrak b_0\mathcal W_2(\mu,\nu)^2+2C_b\big(|x-y|+|a-\tilde a|\big)|x-y|\\
    &\qquad+2C_\sigma^2\big(|x-y|+|a-\tilde a|\big)^2.
\end{align*}
Since $\alpha\in\mathcal A_\gamma^\vartheta$, \eqref{eq: lip feedback} gives $|a-\tilde a|\leq\gamma|x-y|+\vartheta\mathcal W_2(\mu,\nu)$. Consequently,
\begin{align*}
    2C_b\big(|x-y|+|a-\tilde a|\big)|x-y|
    &\leq2C_b(1+\gamma)|x-y|^2+2C_b\vartheta|x-y|\mathcal W_2(\mu,\nu)\\
    &\leq\big(2C_b(1+\gamma)+C_b\vartheta\big)|x-y|^2+C_b\vartheta\mathcal W_2(\mu,\nu)^2,
\end{align*}
while
\begin{align*}
    2C_\sigma^2\big(|x-y|+|a-\tilde a|\big)^2
    &\leq2C_\sigma^2\big((1+\gamma)|x-y|+\vartheta\mathcal W_2(\mu,\nu)\big)^2\\
    &\leq4C_\sigma^2(1+\gamma)^2|x-y|^2+4C_\sigma^2\vartheta^2\mathcal W_2(\mu,\nu)^2.
\end{align*}
Combining the preceding estimates yields
\[
    2\langle\bar b(x,\mu)-\bar b(y,\nu),x-y\rangle+\|\bar\sigma(x,\mu)-\bar\sigma(y,\nu)\|^2\leq-\mathfrak a|x-y|^2+\mathfrak b\,\mathcal W_2(\mu,\nu)^2,
\]
where $\mathfrak a$ and $\mathfrak b$ are given by \eqref{eq: resulting abc}. Finally,
\[
    \mathfrak a-\mathfrak b=\mathfrak a_0-\mathfrak b_0-2C_b(1+\gamma+\vartheta)-4C_\sigma^2\big((1+\gamma)^2+\vartheta^2\big),
\]
which is strictly positive by the first inequality in \eqref{eq: condition}. Hence $\mathfrak a>\mathfrak b\geq0$, and \eqref{eq: assumption 2} follows.

Thus $(\bar b,\bar\sigma)=(b^\alpha,\sigma^\alpha)$ satisfy all parts of Assumptions (A). The well-posedness of the corresponding McKean--Vlasov and self-interacting dynamics then follows from the results recalled in Section \ref{sec:path}.
\end{proof}

\begin{proof}[\textbf{Proof of Lemma \ref{lem: extra conditions feedback}}]
Fix $\alpha\in\mathcal A_\gamma^\vartheta$. We first verify the Lipschitz condition on the diffusion coefficient. By \eqref{eq: b sigma separable},
\begin{align*}
    \|\sigma^\alpha(x,\mu)-\sigma^\alpha(y,\nu)\|^2
    &=\big\|\sigma_0(x,\mu)-\sigma_0(y,\nu)+\sigma_1(x,\alpha(x,\mu))-\sigma_1(y,\alpha(y,\nu))\big\|^2\\
    &\leq2\|\sigma_0(x,\mu)-\sigma_0(y,\nu)\|^2+2\|\sigma_1(x,\alpha(x,\mu))-\sigma_1(y,\alpha(y,\nu))\|^2.
\end{align*}
Using \eqref{extra assumptions b0 sigma0}, \eqref{eq: D control lip}, and \eqref{eq: lip feedback}, we obtain
\begin{align*}
    \|\sigma^\alpha(x,\mu)-\sigma^\alpha(y,\nu)\|^2
    &\leq2K_{\sigma,0}\big(|x-y|^2+\mathcal W_2(\mu,\nu)^2\big)+2C_\sigma^2\big(|x-y|+|\alpha(x,\mu)-\alpha(y,\nu)|\big)^2\\
    &\leq2K_{\sigma,0}\big(|x-y|^2+\mathcal W_2(\mu,\nu)^2\big)+2C_\sigma^2\big((1+\gamma)|x-y|+\vartheta\mathcal W_2(\mu,\nu)\big)^2\\
    &\leq\big(2K_{\sigma,0}+4C_\sigma^2(1+\gamma)^2\big)|x-y|^2+\big(2K_{\sigma,0}+4C_\sigma^2\vartheta^2\big)\mathcal W_2(\mu,\nu)^2.
\end{align*}
Thus the first inequality in \eqref{eq: extra feedback estimates} holds, for instance, with
\begin{equation}\label{eq: K sigma feedback}
    \widetilde K_\sigma:=\max\big\{2K_{\sigma,0}+4C_\sigma^2(1+\gamma)^2,\,2K_{\sigma,0}+4C_\sigma^2\vartheta^2\big\}.
\end{equation}

We next verify the growth condition on the drift at the origin. Using \eqref{eq: b sigma separable},
\[
    |b^\alpha(0,\mu)|^2=|b_0(0,\mu)+b_1(0,\alpha(0,\mu))|^2\leq2|b_0(0,\mu)|^2+2|b_1(0,\alpha(0,\mu))|^2.
\]
By \eqref{extra assumptions b0 sigma0}, we have $|b_0(0,\mu)|^2\leq K_{b,0}\big[1+\mu(|\cdot|^2)\big]$. Moreover, by compactness of $\mathds A$ and continuity of $b_1$, the constant $M_b$ defined in \eqref{eq: Mb Msigma} is finite. Hence
\[
    |b^\alpha(0,\mu)|^2\leq2K_{b,0}\big[1+\mu(|\cdot|^2)\big]+2M_b^2\leq\widetilde K_b\big[1+\mu(|\cdot|^2)\big],
\]
where one may take $\widetilde K_b:=2K_{b,0}+2M_b^2$. This proves the second inequality in \eqref{eq: extra feedback estimates}.

The constants $\widetilde K_\sigma$ and $\widetilde K_b$ depend only on $\gamma,\vartheta$ and on the structural constants in Assumption (C), and are therefore independent of the particular feedback $\alpha\in\mathcal A_\gamma^\vartheta$. Thus \eqref{extra assumptions} holds uniformly over $\mathcal A_\gamma^\vartheta$.
\end{proof}

\begin{proof}[\textbf{Proof of Proposition \ref{prop: lip subclass admissible}}]
Let $\alpha\in\mathcal A_\gamma^\vartheta$. By Lemma \ref{lem: conditions ok},  feedback-composed coefficients $(b^\alpha,\sigma^\alpha)$ satisfy Assumptions (A), so condition \textnormal{(i)} in Definition \ref{def: admissible controls} holds. By Lemma \ref{lem: extra conditions feedback}, they also satisfy \eqref{extra assumptions}, which verifies condition \textnormal{(ii)}.

It remains to verify the regularity condition on the feedback-composed running cost. Let $x,y\in\mathbb R^d$ and $\mu,\eta,\nu\in\mathscr P_2(\mathbb R^d)$. By Assumption (B),
\[
    \big|f(x,\mu,\alpha(x,\eta))-f(y,\mu,\alpha(y,\nu))\big|\leq C_f\big(|x-y|+|\alpha(x,\eta)-\alpha(y,\nu)|\big).
\]
Since $\alpha\in\mathcal A_\gamma^\vartheta$, \eqref{eq: lip feedback} gives $|\alpha(x,\eta)-\alpha(y,\nu)|\leq\gamma|x-y|+\vartheta\mathcal W_2(\eta,\nu)$. Therefore,
\begin{align*}
    \big|f(x,\mu,\alpha(x,\eta))-f(y,\mu,\alpha(y,\nu))\big|
    &\leq C_f(1+\gamma)|x-y|+C_f\vartheta\mathcal W_2(\eta,\nu)\\
    &\leq C_{\gamma,\vartheta}\big(|x-y|+\mathcal W_2(\eta,\nu)\big),
\end{align*}
where one may take $C_{\gamma,\vartheta}:=C_f\max\{1+\gamma,\vartheta\}$. Thus condition \textnormal{(iii)} in Definition \ref{def: admissible controls} is also satisfied. Since $\alpha\in\mathcal A_\gamma^\vartheta$ was arbitrary, $\mathcal A_\gamma^\vartheta\subseteq\mathcal A_{\mathrm{adm}}$.

Finally, the constants obtained in Lemmas \ref{lem: conditions ok} and \ref{lem: extra conditions feedback} depend only on the fixed parameters $\gamma,\vartheta$ and the structural constants in Assumption (C). Likewise, $C_{\gamma,\vartheta}$ depends only on $C_f,\gamma$, and $\vartheta$. Consequently, all constants entering the defining properties of $\mathcal A_{\mathrm{adm}}$ can be chosen uniformly over $\alpha\in\mathcal A_\gamma^\vartheta$.
\end{proof}

\begin{proof}[\textbf{Proof of Corollary \ref{cor: feedback trajectory estimates}}]
Let $\alpha\in\mathcal A_\gamma^\vartheta$. By Proposition \ref{prop: lip subclass admissible}, $\alpha\in\mathcal A_{\mathrm{adm}}$. Hence, by Definition \ref{def: admissible controls}, the feedback-composed coefficients $(b^\alpha,\sigma^\alpha)$ satisfy Assumptions (A) and the additional conditions \eqref{extra assumptions}, with the structural parameters $\mathfrak a,\mathfrak b,\rho_0$ common to all feedbacks in $\mathcal A_\gamma^\vartheta$. Proposition \ref{prop: estimate trajectories} therefore applies to $(X^\alpha,Y^\alpha)$ for every $\varepsilon$ in the stated range and yields \eqref{eq: feedback trajectory estimates}.
\end{proof}

\begin{proof}[\textbf{Proof of Lemma \ref{lem: localized composed lip}}]
Fix $r>R$, and let $\gamma_r$ be the local Lipschitz constant of $\alpha$ on $B_r$ given by \eqref{eq: local lip state}. We first consider $B^\alpha(x,\mu)=\chi(x)\beta(\alpha(x,\mu))$. We distinguish two cases.

\noindent\textit{\textbf{Case 1:} Both $x$ and $y$ belong to $B_r$.}

By the triangle inequality and the Lipschitz continuity of $\chi$ and $\beta$,
\begin{align*}
    |B^\alpha(x,\mu)-B^\alpha(y,\nu)|
    &\leq|\chi(x)-\chi(y)|\,|\beta(\alpha(x,\mu))|+|\chi(y)|\,|\beta(\alpha(x,\mu))-\beta(\alpha(y,\nu))|\\
    &\leq L_\chi M_\beta|x-y|+L_\beta|\alpha(x,\mu)-\alpha(y,\nu)|.
\end{align*}
Using \eqref{eq: local joint lip}, we obtain
\[
    |B^\alpha(x,\mu)-B^\alpha(y,\nu)|\leq(L_\chi M_\beta+L_\beta\gamma_r)|x-y|+L_\beta\vartheta\mathcal W_2(\mu,\nu).
\]

\noindent\textit{\textbf{Case 2:} At least one of $x,y$ lies outside $B_r$.}

Without loss of generality, assume $y\notin B_r$. Since $r>R$, we have $\chi(y)=0$ and therefore $B^\alpha(y,\nu)=0$. If also $|x|\geq R$, then $\chi(x)=0$ and the desired estimate is immediate.

It remains to consider $|x|<R$. Choose $\bar x$ such that $|\bar x|=R$ and $|x-\bar x|=R-|x|$. Since $\chi(\bar x)=0$ and $\chi$ is Lipschitz,
\[
    |\chi(x)|=|\chi(x)-\chi(\bar x)|\leq L_\chi(R-|x|).
\]
Moreover, $|x-y|\geq|y|-|x|\geq R-|x|$. Consequently,
\[
    |B^\alpha(x,\mu)-B^\alpha(y,\nu)|=|\chi(x)\beta(\alpha(x,\mu))|\leq L_\chi M_\beta|x-y|.
\]
Thus the required global estimate also holds in this case.

Combining the two cases gives
\[
    |B^\alpha(x,\mu)-B^\alpha(y,\nu)|\leq\widehat C_{b,\alpha}\big(|x-y|+\mathcal W_2(\mu,\nu)\big),
\]
where one may take $\widehat C_{b,\alpha}:=L_\chi M_\beta+L_\beta(\gamma_r+\vartheta)$.

The proof for $\Sigma^\alpha$ is identical. Indeed, applying the same argument with $M_\varsigma$ and $L_\varsigma$ in place of $M_\beta$ and $L_\beta$, respectively, yields
\[
    \|\Sigma^\alpha(x,\mu)-\Sigma^\alpha(y,\nu)\|\leq\widehat C_{\sigma,\alpha}\big(|x-y|+\mathcal W_2(\mu,\nu)\big),
\]
where one may take $\widehat C_{\sigma,\alpha}:=L_\chi M_\varsigma+L_\varsigma(\gamma_r+\vartheta)$. This proves \eqref{eq: localized composed lip}.
\end{proof}

\begin{proof}[\textbf{Proof of Proposition \ref{prop: local subclass admissible}}]
Fix $\alpha\in\mathcal A_{\mathrm{loc}}^\vartheta$. By \eqref{eq: localized coefficients} and \eqref{eq: localized composed terms}, the feedback-composed coefficients are
\[
    b^\alpha(x,\mu)=b_0(x,\mu)+B^\alpha(x,\mu),\qquad \sigma^\alpha(x,\mu)=\sigma_0(x,\mu)+\Sigma^\alpha(x,\mu).
\]
We verify the three conditions in Definition \ref{def: admissible controls}.

\noindent\textit{\textbf{Step 1} (Verification of Assumptions (A)).} 
By Lemma \ref{lem: localized composed lip}, the maps $B^\alpha$ and $\Sigma^\alpha$ are globally Lipschitz, hence continuous. Since $(b_0,\sigma_0)$ are continuous and bounded on bounded sets, the same holds for $(b^\alpha,\sigma^\alpha)$.\\
We first verify the coercivity condition \eqref{eq: assumption 1}. Since $|B^\alpha(x,\mu)|\leq M_\beta$ and $\|\Sigma^\alpha(x,\mu)\|\leq M_\varsigma$, we obtain
\begin{align*}
    2\langle b^\alpha(x,\mu),x\rangle+(1+\rho_0)\|\sigma^\alpha(x,\mu)\|^2
    &\leq2\langle b_0(x,\mu),x\rangle+2(1+\rho_0)\|\sigma_0(x,\mu)\|^2\\
    &\quad+2\langle B^\alpha(x,\mu),x\rangle+2(1+\rho_0)\|\Sigma^\alpha(x,\mu)\|^2.
\end{align*}
Using \eqref{eq: D coercivity} and $2uv\leq u^2+v^2$, we deduce
\begin{align*}
    2\langle b^\alpha(x,\mu),x\rangle+(1+\rho_0)\|\sigma^\alpha(x,\mu)\|^2
    &\leq-\mathfrak c_0|x|^2+K_0\big[1+\rho_0+\mu(|\cdot|^2)\big]\\
    &\quad+|x|^2+M_\beta^2+2(1+\rho_0)M_\varsigma^2.
\end{align*}
Hence, for a suitable constant $\widehat K>0$,
\[
    2\langle b^\alpha(x,\mu),x\rangle+(1+\rho_0)\|\sigma^\alpha(x,\mu)\|^2\leq-(\mathfrak c_0-1)|x|^2+\widehat K\big[1+\rho_0+\mu(|\cdot|^2)\big].
\]
Since $\mathfrak c_0>1$, condition \eqref{eq: assumption 1} holds with $\rho=\rho_0$ and $\mathfrak c:=\mathfrak c_0-1>0$.

We next verify the confluence condition \eqref{eq: assumption 2}. Let $x,y\in\mathbb R^d$ and $\mu,\nu\in\mathscr P_2(\mathbb R^d)$. Using $\|A+B\|^2\leq2\|A\|^2+2\|B\|^2$,
\begin{align*}
    &2\langle b^\alpha(x,\mu)-b^\alpha(y,\nu),x-y\rangle+\|\sigma^\alpha(x,\mu)-\sigma^\alpha(y,\nu)\|^2\\
    &\quad\leq2\langle b_0(x,\mu)-b_0(y,\nu),x-y\rangle+2\|\sigma_0(x,\mu)-\sigma_0(y,\nu)\|^2\\
    &\qquad+2\langle B^\alpha(x,\mu)-B^\alpha(y,\nu),x-y\rangle+2\|\Sigma^\alpha(x,\mu)-\Sigma^\alpha(y,\nu)\|^2.
\end{align*}
By \eqref{eq: D dissipativity} and Lemma \ref{lem: localized composed lip},
\begin{align*}
    &2\langle b^\alpha(x,\mu)-b^\alpha(y,\nu),x-y\rangle+\|\sigma^\alpha(x,\mu)-\sigma^\alpha(y,\nu)\|^2\\
    &\quad\leq-\mathfrak a_0|x-y|^2+\mathfrak b_0\mathcal W_2(\mu,\nu)^2+2\widehat C_{b,\alpha}\big(|x-y|+\mathcal W_2(\mu,\nu)\big)|x-y|\\
    &\qquad+2\widehat C_{\sigma,\alpha}^2\big(|x-y|+\mathcal W_2(\mu,\nu)\big)^2.
\end{align*}
Using $2uv\leq u^2+v^2$ and $(u+v)^2\leq2u^2+2v^2$, we have
\[
    2\widehat C_{b,\alpha}\big(|x-y|+\mathcal W_2(\mu,\nu)\big)|x-y|\leq3\widehat C_{b,\alpha}|x-y|^2+\widehat C_{b,\alpha}\mathcal W_2(\mu,\nu)^2,
\]
and
\[
    2\widehat C_{\sigma,\alpha}^2\big(|x-y|+\mathcal W_2(\mu,\nu)\big)^2\leq4\widehat C_{\sigma,\alpha}^2|x-y|^2+4\widehat C_{\sigma,\alpha}^2\mathcal W_2(\mu,\nu)^2.
\]
Consequently,
\[
    2\langle b^\alpha(x,\mu)-b^\alpha(y,\nu),x-y\rangle+\|\sigma^\alpha(x,\mu)-\sigma^\alpha(y,\nu)\|^2\leq-\mathfrak a_\alpha|x-y|^2+\mathfrak b_\alpha\mathcal W_2(\mu,\nu)^2,
\]
where
\begin{equation}\label{eq: local resulting ab}
    \mathfrak a_\alpha:=\mathfrak a_0-3\widehat C_{b,\alpha}-4\widehat C_{\sigma,\alpha}^2,\qquad \mathfrak b_\alpha:=\mathfrak b_0+\widehat C_{b,\alpha}+4\widehat C_{\sigma,\alpha}^2.
\end{equation}
Therefore,
\[
    \mathfrak a_\alpha-\mathfrak b_\alpha=\mathfrak a_0-\mathfrak b_0-4\widehat C_{b,\alpha}-8\widehat C_{\sigma,\alpha}^2,
\]
which is strictly positive by \eqref{eq: localized compatibility}. Hence $\mathfrak a_\alpha>\mathfrak b_\alpha\geq0$, and \eqref{eq: assumption 2} follows. Thus $(b^\alpha,\sigma^\alpha)$ satisfy Assumptions (A).

\noindent\textit{\textbf{Step 2} (Verification of \eqref{extra assumptions}).} 
By \eqref{extra assumptions b0 sigma0} and Lemma \ref{lem: localized composed lip},
\begin{align*}
    \|\sigma^\alpha(x,\mu)-\sigma^\alpha(y,\nu)\|^2
    &\leq2\|\sigma_0(x,\mu)-\sigma_0(y,\nu)\|^2+2\|\Sigma^\alpha(x,\mu)-\Sigma^\alpha(y,\nu)\|^2\\
    &\leq2K_{\sigma,0}\big(|x-y|^2+\mathcal W_2(\mu,\nu)^2\big)+2\widehat C_{\sigma,\alpha}^2\big(|x-y|+\mathcal W_2(\mu,\nu)\big)^2\\
    &\leq K_\alpha\big(|x-y|^2+\mathcal W_2(\mu,\nu)^2\big),
\end{align*}
for a suitable constant $K_\alpha>0$. Moreover,
\begin{align*}
    |b^\alpha(0,\mu)|^2
    &\leq2|b_0(0,\mu)|^2+2|B^\alpha(0,\mu)|^2\\
    &\leq2K_{b,0}\big[1+\mu(|\cdot|^2)\big]+2M_\beta^2\leq K'_\alpha\big[1+\mu(|\cdot|^2)\big],
\end{align*}
for a suitable constant $K'_\alpha>0$. Thus condition \textnormal{(ii)} in Definition \ref{def: admissible controls} is satisfied.

\noindent\textit{\textbf{Step 3} (Verification of the running-cost condition).} 
By \eqref{eq: localized running cost},
\begin{align*}
    &\big|f(x,\mu,\alpha(x,\eta))-f(y,\mu,\alpha(y,\nu))\big|\\
    &\quad\leq|f_0(x,\mu)-f_0(y,\mu)|+\big|\chi(x)\ell(\alpha(x,\eta))-\chi(y)\ell(\alpha(y,\nu))\big|.
\end{align*}
The first term is bounded by $L_{f,0}|x-y|$ by \eqref{eq: f0 lip}. For the second term, the argument of Lemma \ref{lem: localized composed lip} applies with $\ell$ in place of $\beta$. Since $\ell$ is Lipschitz and bounded on the compact set $\mathds A$, there exists a constant $\widehat C_{f,\alpha}\geq0$ such that
\[
    \big|\chi(x)\ell(\alpha(x,\eta))-\chi(y)\ell(\alpha(y,\nu))\big|\leq\widehat C_{f,\alpha}\big(|x-y|+\mathcal W_2(\eta,\nu)\big).
\]
Consequently,
\[
    \big|f(x,\mu,\alpha(x,\eta))-f(y,\mu,\alpha(y,\nu))\big|\leq C_\alpha\big(|x-y|+\mathcal W_2(\eta,\nu)\big),
\]
where one may take $C_\alpha:=L_{f,0}+\widehat C_{f,\alpha}>0$, independently of the frozen measure $\mu$. This verifies condition \textnormal{(iii)} in Definition \ref{def: admissible controls}.

All three defining conditions of $\mathcal A_{\mathrm{adm}}$ are therefore satisfied, and hence $\alpha\in\mathcal A_{\mathrm{adm}}$.
\end{proof}

\begin{proof}[\textbf{Proof of Lemma \ref{lem: uniform frozen moments}}]
Fix $\alpha\in\mathcal A_\gamma^\vartheta$. Let $X^\alpha$ be the corresponding McKean--Vlasov process starting from $z$. Applying It\^o's formula to $|X_t^\alpha|^2$, with standard localization, and using \eqref{eq: uniform lyapunov longtime}, we obtain
\[
    \frac{\dd}{\dd t}\mathbb E|X_t^\alpha|^2\leq-\mathfrak c\,\mathbb E|X_t^\alpha|^2+K\big[1+\rho+\mathcal L_t[X^\alpha](|\cdot|^2)\big].
\]
Since $\mathcal L_t[X^\alpha](|\cdot|^2)=\mathbb E|X_t^\alpha|^2$, this gives
\[
    \frac{\dd}{\dd t}\mathbb E|X_t^\alpha|^2\leq-(\mathfrak c-K)\mathbb E|X_t^\alpha|^2+K(1+\rho).
\]
Gr\"onwall's inequality yields
\[
    \mathbb E|X_t^\alpha|^2\leq e^{-(\mathfrak c-K)t}|z|^2+\frac{K(1+\rho)}{\mathfrak c-K}\big(1-e^{-(\mathfrak c-K)t}\big).
\]
By Theorem \ref{thm: conv law to mu}, $\mathcal L_t[X^\alpha]\to\mu^\alpha$ in $\mathcal W_2$. Since this implies convergence of second moments, we conclude that
\[
    \sup_{\alpha\in\mathcal A_\gamma^\vartheta}\mu^\alpha(|\cdot|^2)\leq\frac{K(1+\rho)}{\mathfrak c-K}<+\infty.
\]

We next establish the higher-moment estimate for the frozen dynamics. Set $p:=2+\rho$ and $V(x):=(1+|x|^2)^{p/2}$. Let $\mathscr L^\alpha$ denote the generator of the frozen process $Z^\alpha$. A direct computation gives
\[
\begin{aligned}
    \mathscr L^\alpha V(x)
    &=p(1+|x|^2)^{p/2-1}\langle x,b^\alpha(x,\mu^\alpha)\rangle+\frac p2(1+|x|^2)^{p/2-1}\|\sigma^\alpha(x,\mu^\alpha)\|^2\\
    &\quad+\frac{p(p-2)}2(1+|x|^2)^{p/2-2}|(\sigma^\alpha(x,\mu^\alpha))^\top x|^2.
\end{aligned}
\]
Using $|(\sigma^\alpha(x,\mu^\alpha))^\top x|^2\leq(1+|x|^2)\|\sigma^\alpha(x,\mu^\alpha)\|^2$ and $p-1=1+\rho$, we obtain
\[
    \mathscr L^\alpha V(x)\leq\frac p2(1+|x|^2)^{\rho/2}\big[2\langle b^\alpha(x,\mu^\alpha),x\rangle+(1+\rho)\|\sigma^\alpha(x,\mu^\alpha)\|^2\big].
\]
Applying \eqref{eq: uniform lyapunov longtime} with $\eta=\mu^\alpha$ and using the uniform invariant second-moment bound gives
\[
    \mathscr L^\alpha V(x)\leq-\frac{p\mathfrak c}{2}|x|^2(1+|x|^2)^{\rho/2}+C(1+|x|^2)^{\rho/2},
\]
where $C$ is independent of $\alpha$.

Since $|x|^2(1+|x|^2)^{\rho/2}=V(x)-V(x)^{\rho/(2+\rho)}$, we infer that
\[
    \mathscr L^\alpha V(x)\leq-\frac{p\mathfrak c}{2}V(x)+C'V(x)^{\rho/(2+\rho)}.
\]
Since $\theta:=\rho/(2+\rho)\in(0,1)$, Young's inequality gives $u^\theta\leq\delta u+C_{\delta,\theta}$ for every $\delta>0$ and $u\geq0$. Choosing $\delta=p\mathfrak c/(4C')$ and applying this inequality with $u=V(x)$, we obtain
\[
C'V(x)^\theta\leq\frac{p\mathfrak c}{4}V(x)+C.
\]
Consequently,
\[
\mathscr L^\alpha V(x)\leq-\frac{p\mathfrak c}{4}V(x)+C=-cV(x)+C,
\]
where $c:=p\mathfrak c/4>0$. Both constants are independent of $\alpha$, since $\mathfrak c$ and $C'$ are uniform over $\mathcal A_\gamma^\vartheta$.

To justify the application of It\^o's formula to the unbounded function $V$, let $\tau_n^\alpha:=\inf\{t\geq0:|Z_t^\alpha|\geq n\}$. Applying It\^o's formula to $e^{ct}V(Z_t^\alpha)$ up to $\tau_n^\alpha$, taking expectations, and using the preceding Foster--Lyapunov inequality, we obtain
\[
    \mathbb E\big[e^{c(t\wedge\tau_n^\alpha)}V(Z_{t\wedge\tau_n^\alpha}^\alpha)\big]\leq V(z)+\frac Cc\big(e^{ct}-1\big).
\]
Since the frozen dynamics are nonexplosive, $\tau_n^\alpha\to+\infty$ almost surely. Fatou's lemma therefore gives
\[
    \mathbb E[V(Z_t^\alpha)]\leq e^{-ct}V(z)+\frac Cc(1-e^{-ct}),\qquad t\geq0.
\]
The constants $c,C$ are independent of $\alpha\in\mathcal A_\gamma^\vartheta$. Finally, $|x|^{2+\rho}\leq V(x)$, and hence
\[
    \sup_{\alpha\in\mathcal A_\gamma^\vartheta}\sup_{t\geq0}\mathbb E|Z_t^\alpha|^{2+\rho}<+\infty.
\]
This completes the proof.
\end{proof}

\begin{proof}[\textbf{Proof of Lemma \ref{lem: uniform frozen empirical}}]
Fix $\varepsilon\in(0,1)$ and $\alpha\in\mathcal A_\gamma^\vartheta$. Write $Z=Z^\alpha$, $\mu=\mu^\alpha$, and let $(P_t)_{t\geq0}$ be the Markov semigroup of the frozen diffusion. We follow \cite[Proposition 3.1]{du2023empirical}, indicating which bounds are uniform over the feedback class and specifying the adjustment needed in the covariance estimate.

\smallskip
\noindent\textit{\textbf{Step 1} (Uniform input estimates). }
By \eqref{uniform moment} and Lemma \ref{lem: uniform frozen moments}, the quantities
\[
\begin{aligned}
M_2&:=\sup_{\alpha\in\mathcal A_\gamma^\vartheta}\mu^\alpha(|\cdot|^2)<\infty,\qquad \text{and} \qquad 
M_{2+\rho}:=\sup_{\alpha\in\mathcal A_\gamma^\vartheta}\sup_{s\geq0}\mathbb E|Z_s^\alpha|^{2+\rho}<\infty
\end{aligned}
\]
are finite. Moreover, \eqref{eq: uniform frozen contraction} and invariance of $\mu$ imply
\[
\mathcal W_2(P_t^*\nu,\mu)^2\leq e^{-\mathfrak a t}\mathcal W_2(\nu,\mu)^2,\qquad \nu\in\mathscr P_2(\mathbb R^d),
\]
with the same $\mathfrak a>0$ for every $\alpha$. In particular,
\[
\mathcal W_2(\delta_z,\mu)^2\leq2|z|^2+2M_2,\qquad \sup_{s\geq0}\mathbb E|Z_s|^2\leq M_{2+\rho}^{2/(2+\rho)}.
\]
These replace the process-dependent inputs in \cite[Lemma 3.2]{du2023empirical}. Below, $C$ may change from line to line but depends only on $d$, $\rho$, $\mathfrak a$, $z$, $M_2$, $M_{2+\rho}$, and the fixed mollifier; $C_\varepsilon$ may also depend on $\varepsilon$. Neither depends on $\alpha$, $t$, $r$, or $R$.

\smallskip
\noindent\textit{\textbf{Step 2} (Mollification and averaged-law estimate).}
Fix a nonnegative $\varphi\in C_c^\infty(\mathbb R^d)$ with integral one and support in the unit ball, and put $\varphi_r(x):=r^{-d}\varphi(x/r)$. The mollifier-dependent constants involve only $\|\varphi\|_\infty$, $\|\nabla\varphi\|_\infty$, and $\int_{\mathbb R^d}|x|^2\varphi(x)\,\dd x$. The convolution and decomposition in \cite[Eqs. (3.2)--(3.3)]{du2023empirical} therefore give mollification errors bounded by $Cr^2$, uniformly in $\alpha$.

For the ordinary occupation measure the sampling weight is Lebesgue measure on $[0,1]$, which satisfies the condition $\Pi_2(\varepsilon)$ of \cite{du2023empirical} for every $\varepsilon\in(0,1)$: indeed,
\[
\int_0^1u^{-\varepsilon}\,\dd u=\frac1{1-\varepsilon},\qquad \int_0^1\int_0^1|u-v|^{-\varepsilon}\,\dd u\,\dd v=\frac2{(1-\varepsilon)(2-\varepsilon)}.
\]
The averaged-law term in \cite[Eq. (3.4)]{du2023empirical} is bounded by
\[
\begin{aligned}
\mathcal W_2\left(\frac1t\int_0^t\mathcal L_s[Z] \,\dd s,\mu\right)^2&\leq\frac1t\int_0^t\mathcal W_2(\mathcal L_s[Z],\mu)^2\,\dd s\\
&\leq\frac{\mathcal W_2(\delta_z,\mu)^2}{t}\int_0^t e^{-\mathfrak a s}\,\dd s\leq\frac{2|z|^2+2M_2}{\mathfrak a t}.
\end{aligned}
\]
Thus the decay rate and its constant are uniform over $\alpha$.

\smallskip
\noindent\textit{\textbf{Step 3} (Density coupling and tail estimate).}
To identify the terms in \cite[Eqs. (3.5)--(3.7)]{du2023empirical}, let
\[
p_{t,r}(y):=\frac1t\int_0^t\varphi_r(y-Z_s)\,\dd s,\qquad q_{t,r}(y):=\mathbb E[p_{t,r}(y)].
\]
These are the densities of the mollified occupation measure and of its expectation. By \cite[Lemma 3.3]{du2023empirical}, the expected squared Wasserstein distance between these mollified measures is bounded by $C_d(J_1+J_2)$, where
\[
J_1:=\int_{|y|>R}|y|^2\mathbb E|p_{t,r}(y)-q_{t,r}(y)|\,\dd y,\qquad J_2:=\int_{|y|\leq R}|y|^2\mathbb E|p_{t,r}(y)-q_{t,r}(y)|\,\dd y.
\]
The support of $\varphi_r$ in the ball of radius $r$ and the bound $\sup_s\mathbb E|Z_s|^{2+\rho}\leq M_{2+\rho}$ allow the tail argument of \cite[Eq. (3.5)]{du2023empirical} to be used unchanged, giving
\[
J_1\leq C\big[r^2+(R-r)^{-\rho}\big],\qquad 0<r<R.
\]
Here $C$ depends on $\rho$, $M_{2+\rho}$, and the fixed mollifier, but not on $\alpha$.

\smallskip
\noindent\textit{\textbf{Step 4} (Covariance estimate).}
We follow \cite[Eqs. (3.6)--(3.7)]{du2023empirical}. The relevant process-dependent data are the common contraction rate $\mathfrak a$ in \eqref{eq: uniform frozen contraction} and the bound $\sup_s\mathbb E|Z_s|^2\leq M_{2+\rho}^{2/(2+\rho)}$ above. The only adjustment concerns the dependence on $r$: the fixed rescaling satisfies
\[
\|\varphi_r\|_\infty\leq\|\varphi\|_\infty r^{-d},\qquad \|\nabla\varphi_r\|_\infty\leq\|\nabla\varphi\|_\infty r^{-d-1}.
\]
Set $g_y(x):=\varphi_r(y-x)$ and $H_s(y):=g_y(Z_s)-\mathbb E[g_y(Z_s)]$. By the Markov property and synchronous contraction,
\[
\operatorname{Lip}(P_vg_y)\leq\|\nabla\varphi\|_\infty r^{-d-1}e^{-\mathfrak a v/2}.
\]
Combining this with $\|H_s(y)\|_\infty\leq2\|\varphi\|_\infty r^{-d}$ and the uniform second-moment bound yields, for $0\leq s\leq u$,
\[
\left|\mathbb E[H_s(y)H_u(y)]\right|\leq Cr^{-2d}\min\left\{1,r^{-1}e^{-\mathfrak a(u-s)/2}\right\}.
\]
Indeed, for $0\leq s\leq u$, the Markov property gives
\[
\mathbb E[H_u(y)\mid\mathcal F_s]=P_{u-s}g_y(Z_s)-\mathbb E[P_{u-s}g_y(Z_s)].
\]
Let $\widetilde Z_s$ be an independent copy of $Z_s$. Since $|H_s(y)|\leq2\|\varphi\|_\infty r^{-d}$, we obtain
\[
\begin{aligned}
|\mathbb E[H_s(y)H_u(y)]|&\leq2\|\varphi\|_\infty r^{-d}\,\mathbb E\left|P_{u-s}g_y(Z_s)-\mathbb E[P_{u-s}g_y(Z_s)]\right|\\
&\leq2\|\varphi\|_\infty r^{-d}\operatorname{Lip}(P_{u-s}g_y)\,\mathbb E|Z_s-\widetilde Z_s|\\
&\leq Cr^{-2d-1}e^{-\mathfrak a(u-s)/2},
\end{aligned}
\]
where we used $\mathbb E|Z_s-\widetilde Z_s|\leq2(\mathbb E|Z_s|^2)^{1/2}\leq2M_{2+\rho}^{1/(2+\rho)}$. On the other hand, the boundedness of $H_s(y)$ gives the trivial estimate $|\mathbb E[H_s(y)H_u(y)]|\leq4\|\varphi\|_\infty^2r^{-2d}$. Combining the two bounds yields
\[
|\mathbb E[H_s(y)H_u(y)]|\leq Cr^{-2d}\min\left\{1,r^{-1}e^{-\mathfrak a(u-s)/2}\right\}.
\]
The constant in this estimate depends only on $\|\varphi\|_\infty$, $\|\nabla\varphi\|_\infty$, $M_{2+\rho}^{1/(2+\rho)}$, and $\mathfrak a$, not on $\alpha$. For $0<r<1$, integration over $[0,t]^2$ gives
\[
\mathbb E\left|\frac1t\int_0^t H_s(y)\,\dd s\right|^2\leq\frac{Cr^{-2d}}t\int_0^\infty\min\left\{1,r^{-1}e^{-\mathfrak a v/2}\right\}\,\dd v\leq\frac{Cr^{-2d}}t\big(1+\log(1/r)\big).
\]
Applying Cauchy--Schwarz to $J_2$ and using $\int_{|y|\leq R}|y|^2\,\dd y\leq C_dR^{d+2}$, as in \cite[Eq. (3.6)]{du2023empirical}, we conclude that
\[
J_2\leq CR^{d+2}r^{-d}\left(\frac{1+\log(1/r)}{t}\right)^{1/2},\qquad 0<r<1.
\]
This is the covariance estimate used in place of the bound leading to \cite[Eq. (3.7)]{du2023empirical}; the additional logarithm arises from retaining the $r^{-d-1}$ gradient scaling.

\smallskip
\noindent\textit{\textbf{Step 5} (Choice of the parameters and conclusion).}
We now choose the mollification parameters needed for our result. Recall $\kappa:=\rho/((d+2)(\rho+2))$ and, for $t>1$, set
\[
r=t^{-\frac{\rho\varepsilon}{2(d+2)(\rho+2)}},\qquad R=2t^{\frac{\varepsilon}{(d+2)(\rho+2)}}.
\]
Then $0<r<1<R$, $r^2=t^{-\kappa\varepsilon}$, $(R-r)^{-\rho}\leq t^{-\kappa\varepsilon}$, and
\[
R^{d+2}r^{-d}t^{-1/2}=2^{d+2}t^{-\kappa\varepsilon-(1-\varepsilon)/2},\qquad \log(1/r)=\frac{\kappa\varepsilon}{2}\log t.
\]
Consequently,
\[
J_2\leq Ct^{-\kappa\varepsilon}t^{-(1-\varepsilon)/2}\left(1+\frac{\kappa\varepsilon}{2}\log t\right)^{1/2}\leq C_\varepsilon t^{-\kappa\varepsilon},\qquad t>1,
\]
because $\varepsilon<1$. The mollification, averaged-law, and tail bounds are also $O(t^{-\kappa\varepsilon})$, since $t^{-1}\leq t^{-\kappa\varepsilon}$ for $t\geq1$. Combining the estimates as in the conclusion of \cite[Proposition 3.1]{du2023empirical} yields
\[
\mathbb E\left[\mathcal W_2\big(\mathcal E_t[Z^\alpha],\mu^\alpha\big)^2\right]\leq C_{z,\gamma,\vartheta,\varepsilon}t^{-\kappa\varepsilon},\qquad t>1.
\]
At $t=1$, the same estimate follows after enlarging the constant, since the uniform second-moment bounds give $\mathbb E[\mathcal W_2(\mathcal E_1[Z^\alpha],\mu^\alpha)^2]\leq2M_{2+\rho}^{2/(2+\rho)}+2M_2$. This proves \eqref{eq: uniform frozen empirical}, with a constant independent of $\alpha\in\mathcal A_\gamma^\vartheta$.
\end{proof}

\begin{proof}[\textbf{Proof of Lemma \ref{lem: uniform trajectory comparison}}]
Fix $\bar\varepsilon$ satisfying \eqref{eq: choice uniform comparison exponent} and choose $q$ with $\bar\varepsilon<q<\kappa$. For each $\alpha\in\mathcal A_\gamma^\vartheta$, let $\mu^\alpha$ be the invariant measure and couple $X^\alpha$, $Y^\alpha$, and the frozen diffusion $Z^\alpha$ from \eqref{eq: frozen alpha process} with the same Brownian motion and initial condition $z$. We follow the comparison method of \cite[proof of Lemma 4.1, p. 2505]{du2023empirical}, explaining its specialization to Lebesgue sampling and checking that every comparison constant is independent of the feedback.

\smallskip
\noindent\textit{\textbf{Step 1} (Uniform data and the frozen empirical estimate).}
By the first part of \eqref{uniform moment} and Lemma \ref{lem: uniform frozen moments}, define the finite constants
\[
M_2:=\sup_{\alpha\in\mathcal A_\gamma^\vartheta}\mu^\alpha(|\cdot|^2),\qquad M_{2+\rho}:=\sup_{\alpha\in\mathcal A_\gamma^\vartheta}\sup_{s\geq0}\mathbb E|Z_s^\alpha|^{2+\rho},\qquad K_z:=2|z|^2+2M_2.
\]
The contraction \eqref{eq: uniform frozen contraction} gives $\mathcal W_2(\mathcal L_t[Z^\alpha],\mu^\alpha)^2\leq K_ze^{-\mathfrak a t}$, uniformly in $\alpha$. For $t\geq1$, Lemma \ref{lem: uniform frozen empirical}, applied with $\varepsilon=q/\kappa\in(0,1)$, gives the rate $t^{-q}$. For $0<t\leq1$, the triangle inequality via $\delta_0$ gives
\[
\mathbb E\left[\mathcal W_2\big(\mathcal E_t[Z^\alpha],\mu^\alpha\big)^2\right]\leq\frac2t\int_0^t\mathbb E|Z_s^\alpha|^2\,\dd s+2\mu^\alpha(|\cdot|^2)\leq2M_{2+\rho}^{2/(2+\rho)}+2M_2.
\]
Hence the uniform constant
\[
F_q:=\sup_{\alpha\in\mathcal A_\gamma^\vartheta}\sup_{t>0}\frac{\mathbb E[\mathcal W_2(\mathcal E_t[Z^\alpha],\mu^\alpha)^2]}{1\wedge t^{-q}}<\infty
\]
is well defined, and \eqref{eq: frozen empirical used comparison} holds with $C_{z,\gamma,\vartheta,q}=F_q$. The process-dependent input constants below are thus the common dissipation coefficients $\mathfrak a>\mathfrak b\geq0$ from Assumption \textnormal{(A)}, $K_z$, and $F_q$; which all do not depend on $\alpha$.

If $\mathfrak b=0$, the synchronous-coupling estimates from Assumption \textnormal{(A)} imply $X^\alpha=Y^\alpha=Z^\alpha$ almost surely. Frozen contraction and \eqref{eq: frozen empirical used comparison} then yield \eqref{eq: uniform trajectory comparison}. We therefore assume $\mathfrak b>0$ and fix $\delta>0$ such that $\widehat{\mathfrak b}:=(1+\delta)\mathfrak b<\mathfrak a$ and $\bar\varepsilon<1-\widehat{\mathfrak b}/\mathfrak a$. The choice depends only on $\mathfrak a$, $\mathfrak b$, and $\bar\varepsilon$.

\smallskip
\noindent\textit{\textbf{Step 2} (Law-dependent versus frozen dynamics).}
Set $h_{\texttt L}^\alpha(t):=\mathbb E|X_t^\alpha-Z_t^\alpha|^2$. It\^o's formula, Assumption \textnormal{(A)}, and  $\mathcal W_2(\mathcal L_t[X^\alpha],\mu^\alpha)^2\leq(1+\delta)h_{\texttt L}^\alpha(t)+(1+\delta^{-1})K_ze^{-\mathfrak a t}$ yield
\[
\frac{\dd}{\dd t}h_{\texttt L}^\alpha(t)\leq-(\mathfrak a-\widehat{\mathfrak b})h_{\texttt L}^\alpha(t)+A_ze^{-\mathfrak a t},\qquad A_z:=\mathfrak b(1+\delta^{-1})K_z.
\]
Since $h_{\texttt L}^\alpha(0)=0$, integration gives, with $\lambda:=\mathfrak a-\widehat{\mathfrak b}>0$,
\[
h_{\texttt L}^\alpha(t)\leq A_z\int_0^t e^{-\lambda(t-s)}e^{-\mathfrak a s}\,\dd s\leq\frac{A_z}{\widehat{\mathfrak b}}e^{-\lambda t},\qquad t\geq0.
\]
This proves \eqref{eq: uniform X Z estimate} and \eqref{eq: uniform X Z polynomial}, with constants depending on $\mathfrak a$, $\mathfrak b$, $\delta$, $K_z$, and $\bar\varepsilon$, but not on $\alpha$.

\smallskip
\noindent\textit{\textbf{Step 3} (Self-interacting versus frozen dynamics).}
Set $h_{\texttt E}^\alpha(t):=\mathbb E|Y_t^\alpha-Z_t^\alpha|^2$. By It\^o's formula and Assumption \textnormal{(A)},
\[
\frac{\dd}{\dd t}h_{\texttt E}^\alpha(t)\leq-\mathfrak a h_{\texttt E}^\alpha(t)+\mathfrak b\,\mathbb E\left[\mathcal W_2\big(\mathcal E_t[Y^\alpha],\mu^\alpha\big)^2\right].
\]
The pathwise occupation-measure coupling and the triangle inequality give
\[
\mathcal W_2\big(\mathcal E_t[Y^\alpha],\mathcal E_t[Z^\alpha]\big)^2\leq\frac1t\int_0^t|Y_s^\alpha-Z_s^\alpha|^2\,\dd s
\]
and
\[
\begin{aligned}
\mathcal W_2\big(\mathcal E_t[Y^\alpha],\mu^\alpha\big)^2&\leq(1+\delta)\mathcal W_2\big(\mathcal E_t[Y^\alpha],\mathcal E_t[Z^\alpha]\big)^2\;+(1+\delta^{-1})\mathcal W_2\big(\mathcal E_t[Z^\alpha],\mu^\alpha\big)^2.
\end{aligned}
\]
Taking expectations and using \eqref{eq: frozen empirical used comparison} gives the scalar inequality \eqref{eq: scalar inequality Y Z} with the  uniform  constant
\[
M_q:=\mathfrak b(1+\delta^{-1})F_q,\qquad \frac{\dd}{\dd t}h_{\texttt E}^\alpha(t)\leq-\mathfrak a h_{\texttt E}^\alpha(t)+\widehat{\mathfrak b}\frac1t\int_0^t h_{\texttt E}^\alpha(s)\,\dd s+M_q(1\wedge t^{-q}).
\]
These differential inequalities may be interpreted almost everywhere; the expected squared differences are locally absolutely continuous, which suffices for the comparison argument.

\smallskip
\noindent\textit{\textbf{Step 4} (Scalar comparison and its constant).}
We now specialize \cite[Lemma 4.1 and its proof, p. 2505]{du2023empirical}. In that lemma, the dissipation coefficient, interaction coefficient, forcing exponent, and forcing amplitude correspond respectively to $\mathfrak a$, $\widehat{\mathfrak b}$, $q$, and $M_q$ above. Its time-sampling family is here $\pi_t(\dd u)=\dd u$, independent of $t$. Choose $\bar\varepsilon<r<1-\widehat{\mathfrak b}/\mathfrak a$, and put $p:=r\wedge q>\bar\varepsilon$. The required integrability condition for the polynomial comparison function is
\[
\int_0^1u^{-r}\,\dd u=\frac1{1-r}<\frac{\mathfrak a}{\widehat{\mathfrak b}},\qquad \Delta:=\mathfrak a-\frac{\widehat{\mathfrak b}}{1-p}>0.
\]
Thus the sampling family belongs to $\Pi_1(r)$ in the notation of \cite{du2023empirical}, and the method of its Lemma 4.1 applies with decay exponent $p$. To make the uniformity of the comparison constant transparent, we record the short supersolution calculation used there.

Write $h:=h_{\texttt E}^\alpha$ and $M:=M_q$. The variation-of-constants formula, $h(0)=0$, and $0\leq1\wedge t^{-q}\leq1$ imply, for $H_T:=\sup_{0\leq s\leq T}h(s)$,
\[
H_T\leq\frac{\widehat{\mathfrak b}H_T+M}{\mathfrak a},\qquad \text{hence}\qquad \sup_{t\geq0}h(t)\leq H:=\frac{M}{\mathfrak a-\widehat{\mathfrak b}}.
\]
Set $T_0:=\max\{1,2p/\Delta\}$, choose $A:=\max\{1,2M/\Delta,HT_0^p\}$, and define $w(t):=At^{-p}$ for $t>0$. Because $p<1$, the polynomial comparison computation of \cite[proof of Lemma 4.1]{du2023empirical} reads
\[
w'(t)+\mathfrak a w(t)-\frac{\widehat{\mathfrak b}}t\int_0^t w(s)\,\dd s=At^{-p}\left(\Delta-\frac pt\right)\geq\frac{A\Delta}{2}t^{-p}\geq M(1\wedge t^{-q}),\qquad t\geq T_0.
\]
Also $h(t)\leq H\leq w(t)$ on $(0,T_0]$. For completeness, the same comparison step as in \cite{du2023empirical} can be expressed by setting $v:=(h-w)^+$. Variation of constants and the fact that $h\leq w$ on $(0,T_0]$ imply
\[
v(t)\leq\widehat{\mathfrak b}\int_{T_0}^t\frac1s\int_{T_0}^s v(u)\,\dd u\,\dd s,\qquad t\geq T_0.
\]
Gronwall's inequality applied to $\sup_{T_0\leq u\leq t}v(u)$ gives $v\equiv0$. Therefore $h(t)\leq At^{-p}$ for $t>0$, and the bound $h\leq H$ yields
\[
h_{\texttt E}^\alpha(t)\leq\max\{A,H\}(1\wedge t^{-p}),\qquad t\geq0.
\]
This is \eqref{eq: uniform Y Z estimate} with $r_0=p$. Crucially, $H$, $\Delta$, $T_0$, and $A$ depend only on $\mathfrak a$, $\mathfrak b$, $\delta$, $r$, $q$, and $F_q$, all fixed independently of $\alpha$. Since $p>\bar\varepsilon$, \eqref{eq: uniform Y Z desired exponent} follows.

\smallskip
\noindent\textit{\textbf{Step 5} (Conclusion for states and measures).}
By the triangle inequality and the two state estimates \eqref{eq: uniform X Z polynomial} and \eqref{eq: uniform Y Z desired exponent},
\[
\mathbb E|X_t^\alpha-Y_t^\alpha|^2\leq2h_{\texttt L}^\alpha(t)+2h_{\texttt E}^\alpha(t)\leq Ct^{-\bar\varepsilon},\qquad t\geq1.
\]
For the measure estimate, synchronous coupling and frozen contraction yield
\[
\mathcal W_2\big(\mathcal L_t[X^\alpha],\mu^\alpha\big)^2\leq2h_{\texttt L}^\alpha(t)+2K_ze^{-\mathfrak a t}\leq Ct^{-\bar\varepsilon},\qquad t\geq1.
\]
Since $p<1$, the occupation-measure coupling and the scalar bound imply
\[
\mathbb E\left[\mathcal W_2\big(\mathcal E_t[Y^\alpha],\mathcal E_t[Z^\alpha]\big)^2\right]\leq\frac1t\int_0^t h_{\texttt E}^\alpha(s)\,\dd s\leq Ct^{-p},\qquad t\geq1.
\]
Combining this with \eqref{eq: frozen empirical used comparison} gives $\mathbb E[\mathcal W_2(\mathcal E_t[Y^\alpha],\mu^\alpha)^2]\leq C(t^{-p}+t^{-q})\leq Ct^{-\bar\varepsilon}$. A final application of the triangle inequality through $\mu^\alpha$ yields the second estimate in \eqref{eq: uniform trajectory comparison}. The constants throughout depend on $z$, $\gamma$, $\vartheta$, $\bar\varepsilon$, and the common structural bounds, but not on $\alpha$.
\end{proof}

\begin{proof}[\textbf{Proof of Lemma \ref{lem: continuity occupation measure map}}]
Fix $y,\widetilde y\in C([0,t];\mathbb R^d)$. Since both paths are continuous on a compact interval, their occupation measures have finite second moments. Define
\[
    \pi_{y,\widetilde y}:=\frac1t\int_0^t\delta_{(y_s,\widetilde y_s)}\,\dd s\in\mathscr P_2(\mathbb R^d\times\mathbb R^d).
\]
For every bounded Borel function $\varphi:\mathbb R^d\to\mathbb R$,
\[
    \int_{\mathbb R^d\times\mathbb R^d}\varphi(x)\,\pi_{y,\widetilde y}(\dd x,\dd\widetilde x)=\frac1t\int_0^t\varphi(y_s)\,\dd s=\int_{\mathbb R^d}\varphi(x)\,\mathfrak E_t(y)(\dd x).
\]
Similarly, the second marginal of $\pi_{y,\widetilde y}$ is $\mathfrak E_t(\widetilde y)$. Thus $\pi_{y,\widetilde y}\in\Pi\big(\mathfrak E_t(y),\mathfrak E_t(\widetilde y)\big)$, and the definition of $\mathcal W_2$ yields
\[
    \mathcal W_2\big(\mathfrak E_t(y),\mathfrak E_t(\widetilde y)\big)^2\leq\int_{\mathbb R^d\times\mathbb R^d}|x-\widetilde x|^2\,\pi_{y,\widetilde y}(\dd x,\dd\widetilde x)=\frac1t\int_0^t|y_s-\widetilde y_s|^2\,\dd s.
\]
The second inequality in \eqref{eq: occupation measure map Lipschitz} follows immediately from the definition of the uniform norm. This proves the result.
\end{proof}

\begin{proof}[\textbf{Proof of Lemma \ref{lem: Wasserstein averaged random measure}}]
We first verify that \eqref{eq: mean random probability measure} defines a probability measure. Since $\boldsymbol\nu$ is almost surely a probability measure, $\overline{\boldsymbol\nu}(\varnothing)=0$ and $\overline{\boldsymbol\nu}(\mathbb R^d)=1$.  If $(A_n)_{n\geq1}$ are pairwise disjoint Borel sets, then \vspace*{-0.5em}
\[
    \boldsymbol\nu\left(\bigcup_{n\geq1}A_n\right)=\sum_{n\geq1}\boldsymbol\nu(A_n)\qquad\text{a.s.}
\]
Taking expectations and applying the monotone convergence theorem yields
\[
    \overline{\boldsymbol\nu}\left(\bigcup_{n\geq1}A_n\right)=\sum_{n\geq1}\mathbb E[\boldsymbol\nu(A_n)]=\sum_{n\geq1}\overline{\boldsymbol\nu}(A_n).
\]
Thus $\overline{\boldsymbol\nu}\in\mathscr P(\mathbb R^d)$.

By definition, for every Borel set $A$,
\[
    \int_{\mathbb R^d}\mathbf 1_A(x)\,\overline{\boldsymbol\nu}(\dd x)=\mathbb E\left[\int_{\mathbb R^d}\mathbf 1_A(x)\,\boldsymbol\nu(\dd x)\right].
\]
By linearity, the same identity holds for nonnegative simple functions. Approximating a nonnegative Borel function by an increasing sequence of nonnegative simple functions and applying monotone convergence, we obtain  \vspace*{-0.5em}
\begin{equation}\label{eq: integration mean random measure}
    \int_{\mathbb R^d}h(x)\,\overline{\boldsymbol\nu}(\dd x)=\mathbb E\left[\int_{\mathbb R^d}h(x)\,\boldsymbol\nu(\dd x)\right]
\end{equation}
for every nonnegative Borel function $h:\mathbb R^d\to[0,+\infty]$. By considering positive and negative parts, the identity also holds for every bounded Borel function, and more generally whenever the corresponding absolute integrals are finite.

Taking $h(x)=|x|^2$ and using \eqref{eq: random measure second moment}, we obtain $\overline{\boldsymbol\nu}(|\cdot|^2)=\mathbb E\big[\boldsymbol\nu(|\cdot|^2)\big]<+\infty$, so $\overline{\boldsymbol\nu}\in\mathscr P_2(\mathbb R^d)$.

Fix $\mu\in\mathscr P_2(\mathbb R^d)$. The map $\eta\mapsto\mathcal W_2(\mu,\eta)^2$ is continuous on $\mathscr P_2(\mathbb R^d)$, so $\mathcal W_2(\mu,\boldsymbol\nu)^2$ is measurable. Moreover, the triangle inequality gives
\[
    \mathcal W_2(\mu,\boldsymbol\nu)^2\leq2\mathcal W_2(\mu,\delta_0)^2+2\mathcal W_2(\boldsymbol\nu,\delta_0)^2=2\mu(|\cdot|^2)+2\boldsymbol\nu(|\cdot|^2).
\]
Therefore,
\begin{equation}\label{eq: random Wasserstein integrability}
    \mathbb E\big[\mathcal W_2(\mu,\boldsymbol\nu)^2\big]\leq2\mu(|\cdot|^2)+2\mathbb E\big[\boldsymbol\nu(|\cdot|^2)\big]<+\infty.
\end{equation}

We now apply Kantorovich duality. Since $\mathbb R^d$ is Polish and the cost $c(x,y)=|x-y|^2$ is continuous and nonnegative, the Kantorovich duality theorem \cite[Theorem 1.42]{santambrogio2015optimal} gives, for every $\eta\in\mathscr P_2(\mathbb R^d)$,
\begin{equation}\label{eq: Wasserstein Kantorovich duality appendix}
    \mathcal W_2(\mu,\eta)^2=\sup_{(\varphi,\psi)\in\mathscr D}\left\{\int_{\mathbb R^d}\varphi(x)\,\mu(\dd x)+\int_{\mathbb R^d}\psi(y)\,\eta(\dd y)\right\},
\end{equation}
where
\[
    \mathscr D:=\left\{(\varphi,\psi)\in C_b(\mathbb R^d)\times C_b(\mathbb R^d):\ \varphi(x)+\psi(y)\leq|x-y|^2,\quad x,y\in\mathbb R^d\right\}.
\]

Fix $(\varphi,\psi)\in\mathscr D$. Since $\varphi$ and $\psi$ are bounded, the integration identity \eqref{eq: integration mean random measure} gives
\[
    \int_{\mathbb R^d}\varphi(x)\,\mu(\dd x)+\int_{\mathbb R^d}\psi(y)\,\overline{\boldsymbol\nu}(\dd y)=\mathbb E\left[\int_{\mathbb R^d}\varphi(x)\,\mu(\dd x)+\int_{\mathbb R^d}\psi(y)\,\boldsymbol\nu(\dd y)\right].
\]
For every realization of $\boldsymbol\nu$, the pair $(\varphi,\psi)$ is admissible in the dual problem \eqref{eq: Wasserstein Kantorovich duality appendix}. Hence,
\[
    \int_{\mathbb R^d}\varphi(x)\,\mu(\dd x)+\int_{\mathbb R^d}\psi(y)\,\boldsymbol\nu(\dd y)\leq\mathcal W_2(\mu,\boldsymbol\nu)^2\qquad\text{a.s.}
\]
Taking expectations and using \eqref{eq: random Wasserstein integrability}, we obtain
\[
    \int_{\mathbb R^d}\varphi(x)\,\mu(\dd x)+\int_{\mathbb R^d}\psi(y)\,\overline{\boldsymbol\nu}(\dd y)\leq\mathbb E\big[\mathcal W_2(\mu,\boldsymbol\nu)^2\big].
\]
Since this holds for every $(\varphi,\psi)\in\mathscr D$, taking the supremum and applying \eqref{eq: Wasserstein Kantorovich duality appendix} with $\eta=\overline{\boldsymbol\nu}$ proves
\[
    \mathcal W_2\big(\mu,\overline{\boldsymbol\nu}\big)^2\leq\mathbb E\big[\mathcal W_2(\mu,\boldsymbol\nu)^2\big].
\]
\end{proof}

\bibliographystyle{plain}
\bibliography{bibliography}

\end{document}